\documentclass[reqno]{amsart}
\usepackage{amscd,amsmath,amsfonts}
\usepackage[all,cmtip]{xy}
\usepackage[dvips]{graphicx}

\newcommand*{\GL}{{\mathrm{GL}}}
\newcommand*{\oo}{{\mathsf{o}}}
\newcommand*{\ra}{{\rightarrow}}
\newcommand*{\R}{{\mathbb R}}
\newcommand*{\CC}{{\mathbb{C}}}

\newcommand*{\Gr}{{\mathrm{Gr}}}
\newcommand*{\Hb}{{\mathbb H}}
\newcommand*{\Ub}{{\mathbb U}}
\newcommand*{\Vb}{{\mathbb V}}

\newcommand*{\ib}{{\mathfrak i}}
\newcommand*{\V}{U} 
\newcommand*{\U}{V}
\newcommand*{\Rey}{\mathrm{R_e}}
\newcommand*{\Ro}{\mathrm{R_o}}
\newcommand*{\Ek}{\mathrm{E_k}}

\theoremstyle{plain}
\newtheorem{theorem}{Theorem}
\newtheorem{corollary}{Corollary}
\newtheorem{lemma}{Lemma}
\theoremstyle{definition}

\theoremstyle{remark}

\title[Grassmannian shooting]{Grassmannian spectral shooting}
\author[Ledoux, Malham and Th\"ummler]{Veerle Ledoux$^1$, 
Simon J.A. Malham$^2$ and Vera Th\"ummler$^3$}
\address{$^1$Vakgroep Toegepaste Wiskunde en Informatica, Ghent University,
Krijgslaan, 281-S9, B-9000 Gent, Belgium\\
$^2$Department of Mathematics,
Heriot-Watt University, Edinburgh EH14 4AS, UK \\
$^3$Fakult\"at f\"ur Mathematik, Universit\"at Bielefeld,
33501 Bielefeld, Germany} 

\date{6th July 2009}

\begin{document}
\begin{abstract}
We present a new numerical method for computing the 
pure-point spectrum associated with the linear
stability of coherent structures. In the context
of the Evans function shooting and matching approach,
all the relevant information is carried by the flow
projected onto the underlying Grassmann manifold. 
We show how to numerically construct this projected 
flow in a stable and robust manner. In particular, 
the method avoids representation singularities by, in practice, 
choosing the best coordinate patch representation for the flow
as it evolves. The method is analytic in the spectral parameter 
and of complexity bounded by the order of the spectral problem 
cubed. For large systems it represents a competitive 
method to those recently developed that are 
based on continuous orthogonalization. We demonstrate this 
by comparing the two methods in three applications: 
Boussinesq solitary waves, autocatalytic travelling waves
and Ekman boundary layer.
\end{abstract}

\maketitle

\section{Introduction}
We introduce a new numerical method for solving high order
linear spectral problems by shooting and matching. The numerical
construction of pure-point spectra is important in determining 
the linear stability of coherent structures. Examples of such
structures are: ground and higher excited states of molecules in 
quantum chemistry (Johnson~\cite{Johnson}, Hutson~\cite{Hu},
Gray and Manopoulous~\cite{GM}, Manopoulous and Gray~\cite{MG},
Chou and Wyatt~\cite{CW1,CW2}, Ledoux~\cite{L},
Ledoux, Van Daele and Vanden Berghe~\cite{LVV}, Ixaru~\cite{Ixtalk});
nonlinear travelling fronts in reaction-diffusion
such as autocatalysis or combustion (Billingham and Needham~\cite{BN}, 
Metcalf, Merkin and Scott~\cite{MMS}, Doelman, Gardner and Kaper~\cite{DGK}, 
Terman~\cite{T}, Gubernov, Mercer, Sidhu and Weber~\cite{GMSW}); 
nerve impulses (Alexander, Gardner and Jones~\cite{AGJ}), 
neural waves (Coombes and Owen~\cite{CO}); solitary waves or steady flows over 
compliant surfaces (Pego and Weinstein~\cite{PW}, 
Alexander and Sachs~\cite{AS}, Chang, Demekhin and Kopelevich~\cite{CDK},
Kapitula and Sandstede~\cite{KS},
Bridges, Derks and Gottwald~\cite{BDG},
Allen~\cite{A}, Allen and Bridges~\cite{AB});
laser pulses (Swinton and Elgin~\cite{SE});
nonlinear waves along elastic rods (Lafortune and Lega~\cite{LL});
ionization fronts (Derks, Ebert and Meulenbroek~\cite{DEM}) or 
spiral waves (Sandstede and Scheel~\cite{SandScheel}).

For such problems the matching condition is typically a 
discriminant known as the Evans function
(Evans~\cite{Evans}, Alexander, Gardner and Jones~\cite{AGJ}) 
or miss-distance function (Pryce~\cite{Pryce}, Greenberg and Marletta~\cite{GM}).
It measures the degree of (possibly transversal) intersection of
the stable and unstable solution subspaces satisfying the longitudinally 
separated far field boundary data. The stable subspace decays
in the direction of wave propagation whilst the unstable subspace
decays in the opposite direction. Equivalently, the Evans function 
is the determinant of the set of solution vectors spanning both
subspaces. With this end-goal matching condition in mind, 
the problem boils down to how to numerically construct  
the solution subspaces in a robust fashion, as well as where to
match longitudinally. This is especially difficult for large
scale problems. These might either emerge from high order systems,
or more specifically, we envisage the stability of nonlinear 
travelling waves with multi-dimensional structure; for example
wrinkled fronts travelling in a fixed longitudinal direction.
The transverse structural information can be projected onto
finite dimensional transverse basis, generating a large linear 
spectral problem posed on the one-dimensional longitudinal 
coefficient function set (see 
Ledoux, Malham, Niesen and Th\"ummler~\cite{LMNT}).
Hitherto such large problems could not be solved by 
shooting and matching and were resolved by projecting the 
whole problem onto a finite basis and solving the resulting
large algebraic eigenvalue problem. However recently, in
the context of the Evans function, Humpherys and Zumbrun~\cite{HZ}
proposed continuous orthogonalization as a viable approach 
to help make large scale problems amenable to shooting and matching.
Here we provide our own answer.

We propose the \emph{new} Grassmann Gaussian elimination method (GGEM) 
which resolves several numerical problems all in one, in particular it:
\begin{enumerate}
\item Evolves the solution along the underlying Grassmann
manifold avoiding representation singularities.
\item Retains analyticity in the spectral parameter.
\item Allows for matching anywhere in the longitudinal 
computational domain. 
\item Has polynomial complexity, operations are 
of the order of the size of the system cubed.
\item Naturally evolves the solution in what in practice
is the optimal coordinate representation (generated
by optimal partial pivoting). 
\end{enumerate}
For each property what is new, expected, proved, 
numerically observed, and its context?

First, evolving the solution subspaces, considered as curves in 
the Grassmann manifold is not new for autonomous problems---see 
Hermann and Martin~\cite{HMc,HM1,HM2,HM3,HM4}, Martin and Hermann~\cite{MH}
Brockett and Byrnes~\cite{BB}, Shayman~\cite{Sh},
Rosenthal~\cite{Ro1}, Ravi, Rosenthal and Wang~\cite{RRW},
Zelikin~\cite{Ze}, Abou--Kandil, Freiling, Ionescu and Jank~\cite{Abu}
and Bittanti, Laub and Willems~\cite{BLW}.
Using Riccati systems to solve nonautonomous spectral problems 
is also not new---see Johnson~\cite{Johnson}, 
Hutson~\cite{Hu}, Pryce~\cite{Pryce}, Manopoulous and Gray~\cite{MG},
Gray and Manopoulous~\cite{GrM} and Chou and Wyatt~\cite{CW1,CW2}.
Here Riccati systems correspond to the flow of the linear
spectral problem projected onto the Grassmann manifold 
with a fixed coordinate patch representation---see Schneider~\cite{Schneider}.
Also see Schiff and Shnider~\cite{SS} and Chou and Wyatt~\cite{CW2}
who use this connection to integrate Riccati systems through singularities.

That we have a numerical method that avoids representation singularities 
for nonautonomous systems is \emph{new}. The idea is as follows.
Given data that lies in a suitable coordinate patch of the Grassmann
manifold, pullback to the Stiefel manifold. Evolve the solution
one steplength along the Stiefel manifold either directly using
a Runge--Kutta method or a Lie group method. Then project onto a suitable and
possibly different coordinate patch of the Grassmann manifold using
optimal Gaussian elimination. In this last step, the practical method
we propose picks a quasi-optimal coordinate patch
in which to best represent the solution in the Grassmann manifold
(see below and Section~\ref{sec:NM}). 

Note that on first inspection the Stiefel manifold is the 
direct natural setting for the stable and unstable solution 
subspaces. Afterall in each case we have to construct the 
full set of solutions to a large system of differential equations
satisfying the correct respective asymptotic boundary conditions 
in the far field.
Each solution set represents a curve in the Stiefel manifold
of dimension commensurate with the size of the solution set.
That the spectral problems are linear means that all the
relevant spectral information can be reconstructed from the flow in
the corresponding Grassmann manifold. That the matching condition
is determinental, means that we only need the Grassmann flow
information---see Martin and Hermann~\cite{MH} and Brockett and
Byrnes~\cite{BB} where this reduction was first considered for autonomous
linear control problems  (in practice we will also need to retain
a complex scalar field to ensure analytic dependence on parameters).
This reduction is crucial because long-range integration along
the Stiefel manifold has been problematic (due to multiple distinct 
exponential growth and decay rates) and one of the 
simplest resolutions in the Evans function context was to use 
Pl\"ucker coordinates---whilst ignoring the Pl\"ucker relations
(more on these below).

Second, retaining analyticity away from the essential spectrum is
\emph{standard} for any shooting method; we prove analyticity 
in Section~\ref{sec:SP}. This allows for a global search 
for eigenvalues in that region by numerically computing
the change in argument of the Evans function round any 
closed contour. Invoking the argument principle, this integer value 
represents the number of zeros of the Evans function, and hence the number 
of eigenvalues counting multiplicity, inside the contour; 
see Brockett and Byrnes~\cite{BB}, 
Alexander, Gardner and Jones~\cite{AGJ} and Ying and Katz~\cite{YK}.

Third, that our method allows matching anywhere in the longitudinal
domain is \emph{new}. We provide substantive numerical evidence.
Previously, other than in trivial cases, most numerical practioners
used the common-sense rule of thumb of integrating the spectral
problem from both ends of the longitudinal domain and matching
at a point roughly centered on the front (which is assumed to be 
localized). When solving the linear problem with Pl\"ucker coordinates, 
it was first important to rescale for the far field spatial
behaviour to neutralize its total exponential growth. Integrating
from the far field initial conditions (a subset of the spatial eigenvectors),
the solution remains roughly constant until the coefficient 
matrix starts to reveal its nonautonmous character due to the 
integration step impinging on the front. Accuracy is retained whilst 
integrating through the front, but thereafter the problem becomes
stiff. The issue is that the numerical methods cannot resolve the
simultaneous exponential growth and decay character associated with
the other far-end stable and unstable subspaces. 

Fourth, having polynomial complexity is \emph{essential} and
we provide here a \emph{new} alternate. 
After Humpherys and Zumbrun~\cite{HZ} introduced their 
continuous orthogonalization method in this context,
which also has polynomial complexity, any new numerical 
spectral shooting method should have this property and 
also be competitive. Previous successful methods used  
Pl\"ucker coordinates, ignoring the quadratic Pl\"ucker relations 
(see Section~\ref{sec:GS}). They solved the flow for the 
corresponding linear vector field in the higher dimensional 
Pl\"ucker embedding space. Details of this Pl\"ucker coordinate or
compound matrix approach can be found in, for example, 
Alexander and Sachs~\cite{AS}, Brin~\cite{Br,Brp}
and Allen and Bridges~\cite{AB}. Unfortunately the number of 
Pl\"ucker coordinates typically grows exponentially with the 
order of the original system, and so this approach cannot be 
used for medium to large order systems. 
However the continuous orthogonalization method
of Humpherys and Zumbrun, and our method, are especially suited to
large scale problems. 

Fifth, our method for dynamic \emph{practical optimal} coordinate 
representation is \emph{new}. Given data on the Stiefel manifold,
for example generated by advancing the solution one steplength
along the Stiefel manifold, how can we project down onto the 
Grassmann manifold using the best representation patch possible?
The idea is as follows. The natural map projection from the 
Stiefel to the Grassmann manifolds is a linear fractional 
map (Milnor and Stasheff~\cite{MS}; Martin and Hermann~\cite{MH}). 
This map represents the action of equivalencing by transformations
whose rank matches that of the Stiefel manifold (this takes us from
the space of frames to the space of spaces spanning those frames).
The Stiefel manifold has a non-square matrix representation.
Projection onto the Grassmann manifold corresponds to equivalencing
by a full rank submatrix of the non-square Stiefel matrix---this renders 
the corresponding submatrix as the identity matrix. We are 
free to choose which submatrix to equivalence by, each distinct choice
corresponds to the matrix representation of a coordinate patch on
the Grassmann manifold. We can use Gaussian elimination,
via elementary column operations, to render any given full rank submatrix 
of the Stiefel matrix as the identity matrix. The key is to 
try to pick the full rank submatrix which has the largest determinant---this
corresponds to choosing the Grassmannian patch that gives the best representation
for the projection from the Stiefel to Grassmann manifold. 
Ideally we would check the size of every full rank submatrix of the 
Stiefel matrix and equivalence by the one with the largest determinant.
However this is an NP problem (equivalent to using the Pl\"ucker coordinates
described above). We provide a practical solution of 
polynomial complexity. The method maximises the pivot used at each step of
the Gaussian elimination process. In the current context, we 
call it quasi-optimal Gaussian elimination. 

Our paper is organised as follows. In Section~\ref{sec:GS}, 
we provide a tailored review of Grassmann manifolds 
and their representation. We then show, in Section~\ref{sec:VF} 
how flows generated by linear vector fields on the Stiefel manifold, 
produce a natural flow on the underlying Grassmann manifold that
is decoupled from the flow through the remaining fibres. We 
subsequently show how this leads to using Riccati systems to
resolve spectra, but that singularities that develop in the 
Riccati flows present spectral matching problems.
In Section~\ref{sec:PGI} we introduce two new practical approaches 
to avoiding these representation singularities. One is the idea behind
our main method, the Grassmann Gaussian elimination method.
The other is a modification of the Riccati approach that changes
the coordinate patch when deemed necessary. Also in this section
we show the connection between the Riccati and  
continuous orthogonalization approaches. 
We present our proposed Grassmann Gaussian elimination
method fully in Section~\ref{sec:NM}, including details 
of how in practice to choose the quasi-optimal 
Grassmannian coordinate representation patch.
We review the Evans function in Section~\ref{sec:SP} and discuss
further simple practical numerical refinements that retain analyticity 
and prevent potential numerical overflow.
We then implement and compare all the competing numerical 
methods in Section~\ref{sec:ES} in three distinct applications.
Finally in Section~\ref{sec:conclu} we conclude and 
present future directions.

\section{Review: Grassmann manifolds}\label{sec:GS}

\subsection{Grassmann and Stiefel manifolds}
A $k$-frame is a $k$-tuple of $k\leq n$ linearly independent
vectors in $\CC^n$. The \emph{Stiefel manifold} $\Vb(n,k)$ of $k$-frames 
is the open subset of $\CC^{n\times k}$ of all $k$-frames 
centred at the origin. The set of $k$ dimensional subspaces of $\CC^n$
form a complex manifold $\mathrm{Gr}(n,k)$ called
the \emph{Grassmann manifold} of $k$-planes in $\CC^n$
(see Steenrod~\cite[p.~35]{S} or Griffiths and Harris~\cite[p.~193]{GH}).

The fibre bundle $\pi\colon\Vb(n,k)\ra\Gr(n,k)$ is a 
\emph{principle fibre bundle}. For each $y$ in the 
base space $\Gr(n,k)$, the inverse image $\pi^{-1}(y)$ is homeomorphic 
to the fibre space $\GL(k)$ which is a Lie group---see 
Montgomery~\cite[p.~151]{M}. 
The projection map $\pi$ is the natural quotient map 
sending each $k$-frame centered at the origin to the 
$k$-plane it spans---see Milnor and Stasheff~\cite[p.~56]{MS}.

\subsection{Representation}
Following the exposition in Griffiths and Harris~\cite{GH},  
any $k$-plane in $\CC^n$ can be represented by an 
$n\times k$ matrix of rank $k$, say $Y\in\CC^{n\times k}$. 
Any two such matrices $Y$ and $Y'$ represent the 
same $k$-plane element of $\Gr(n,k)$ if and only if 
$Y'=Yu$ for some $u\in\mathrm{GL}(k)$ (the $k$-dimensional
subspace elements are invariant to rank $k$ closed transformations
mapping $k$-planes to $k$-planes).

Let $\ib=\{i_1,\ldots,i_k\}\subset\{1,\ldots,n\}$ denote a multi-index
of cardinality $k$. Let $Y_{\ib^\circ}\subset\CC^n$ denote the $(n-k)$-plane
in $\CC^n$ spanned by the vectors $\{e_j\colon j\not\in\ib\}$ and
\begin{equation*}
\Ub_\ib=\bigl\{Y\in\Gr(n,k)\colon Y\cap Y_{\ib^\circ}=\{0\}\bigr\}\,.
\end{equation*}
In other words, $\Ub_\ib$ is the set of $k$-planes $Y\in\Gr(n,k)$
such that the $k\times k$ submatrix of one, and hence any,
matrix representation of $Y$ is nonsingular (representing a
coordinate patch labelled by $\ib$).

Any element of $\Ub_\ib$ has a unique matrix representation 
$y_{\ib^\circ}$ whose $\ib$th $k\times k$ submatrix is the identity matrix. 
For example, if $\ib=\{1,\ldots,k\}$ then any element of 
$\Ub_{\{1,\ldots,k\}}$ can be uniquely represented by a matrix 
of the form
\begin{equation*}
y_{\ib^\circ}=\begin{pmatrix} 1 & 0 & \cdots & 0 \\  
 0 & 1 & \cdots & 0 \\  
\vdots & \vdots & \ddots & \vdots \\ 
 0 & 0 & \cdots & 1 \\ 
\hat y_{k+1,1} & \hat y_{k+1,2} & \cdots & \hat y_{k+1,k}  \\
\hat y_{k+2,1} & \hat y_{k+2,2} & \cdots & \hat y_{k+2,k} \\
 \vdots & \vdots & \ddots & \vdots \\
\hat  y_{n,1} & \hat y_{n,2} &\cdots & \hat y_{n,k} 
\end{pmatrix}\,,
\end{equation*}
where $\hat y_{i,j}\in\CC$ for $i=k+1,\ldots,n$ and $j=1,\ldots,k$.
Conversely, a $n\times k$ matrix of this form represents 
a $k$-plane in $\Ub_\ib$. Each coordinate patch $\Ub_\ib$ 
is an open, dense subset of $\Gr(n,k)$ and the union of
all such patches covers $\Gr(n,k)$.
For each $\ib$, there is a bijective map 
$\varphi_\ib\colon\Ub_\ib\rightarrow \CC^{(n-k)k}$
given by
\begin{equation*}
\varphi_\ib\colon y_{\ib^\circ}\mapsto\hat y\,.
\end{equation*}
Each $\varphi_\ib$ is thus a local coordinate chart
for the coordinate patch $\Ub_\ib$ of $\Gr(n,k)$. 
For all $\ib,\ib'$, if $Y\in\Ub_\ib\cap\Ub_{\ib'}$ and $u_{\ib,\ib'}$
is the $\ib'$th $k\times k$ submatrix of $y_{\ib^\circ}$, then 
$y_{(\ib')^\circ}=y_{\ib^\circ}(u_{\ib,\ib'})^{-1}$. 
Since $u_{\ib,\ib'}$ represents the transformation between representative patchs
and depends holomorphically on $y_{\ib^\circ}$, we deduce 
$\varphi_\ib\circ\varphi_{\ib'}^{-1}$ is holomorphic. 
Note that $\Gr(n,k)$ has a structure of a complex manifold
(see Griffiths and Harris~\cite[p.~194]{GH}).
Further the unitary group $\Ub(n)$ acts continuously and surjectively
on $\Gr(n,k)$. Hence $\Gr(n,k)$ is compact and connected.
Lastly, the general linear group $\GL(n)$ acts transitively on
$\Gr(n,k)$ and it is a homogeneous manifold isomorphic 
to $\GL(n)/\GL(n-k)\times\GL(k)$ 
(see Chern~\cite[p.~65]{C} or Warner~\cite[p.~130]{W}).

\subsection{Pl\"ucker embedding} 
There is a natural map, the \emph{Pl\"ucker map},
\begin{equation*}
p\colon \Gr(n,k)\ra\mathbb P\bigl(\textstyle{\bigwedge^k\CC^n}\bigr)
\end{equation*}
that sends each $k$-plane with basis $Y=[Y_1 \ldots Y_k]$ 
to $Y_1\wedge\ldots\wedge Y_k$; here $\mathbb P^n$
denotes the complex projective space of dimension $n$.
See Griffiths and Harris~\cite{GH} or Coskun~\cite{Coskun}
for more details. If we change the basis, the basis for
the image changes by the determinant of the transformation
matrix. Hence the map is a point in 
$\mathbb P\bigl(\bigwedge^k\CC^n\bigr)$. 
We can recover $Y$ from its image $Y_1\wedge\ldots\wedge Y_k$ as
the set of all vectors $v$ such that $v\wedge Y_1\wedge\ldots\wedge Y_k=0$.
Further, a point of $\mathbb P\bigl(\bigwedge^k\CC^n\bigr)$ is in
the image of $p$ if and only if its representation as a linear
combination of the basis elements of $\bigwedge^k\CC^n$, 
consisting of all possible distinct wedge products
of a $k$-dimensional basis in $\CC^n$, is completely decomposable.
Hence the image of $p$ is a subvariety of $\mathbb P\bigl(\bigwedge^k\CC^n\bigr)$
of completely decomposable elements. It can also be realized as follows. 
A natural coordinatization of $\mathbb P\bigl(\bigwedge^k\CC^n\bigr)$ is
through the determinants of all the $k\times k$ submatrices
of $Y$, normalized by a chosen minor characterized by an index $\ib$,
hence
$\mathbb P\bigl(\bigwedge^k\CC^n\bigr)\cong\mathbb P^{\left(\substack{n\\k}\right)-1}$.
These minor determinants---the Pl\"ucker coordinates---are not all 
independent, indeed, they
satisfy quadratic relations known as the \emph{Pl\"ucker relations}
(which may themselves not all be independent). 
The image of the Pl\"ucker map $p$ is thus the subspace
of $\mathbb P^{\left(\substack{n\\k}\right)-1}$ cut out by 
the quadratic Pl\"ucker relations.

\section{Grassmannian flows} \label{sec:VF}

\subsection{Tangent space decomposition}
Recall that we can consider the Stiefel manifold 
as a principle fibre bundle $\pi\colon\Vb(n,k)\to\Gr(n,k)$. 
Our goal here is to characterize the induced decomposition 
of the tangent space $\mathrm{T}_Y\Vb(n,k)$ for 
$Y\in\Vb(n,k)$.
We can decompose the tangent space $\mathrm{T}_Y\Vb(n,k)$
into horizontal and vertical subspaces 
(see, for example, Montgomery~\cite[p.~149]{M})
\begin{equation*}
\mathrm{T}_Y\Vb(n,k)=\Hb_Y\oplus\Vb_Y\,.
\end{equation*}
The horizontal subspace $\Hb_Y$ is 
associated with the tangent space of the Grassmannian base space, 
while the vertical subspace $\Vb_Y$ is associated with
the fibres homeomorphic to $\GL(k)$.
Let us choose the coordinate patch representation $\Ub_\ib$ for $\Gr(n,k)$
for some $\ib=\{i_1,\ldots,i_k\}\subset\{1,\ldots,n\}$.
Let $P_{\ib}$ denote the projection matrix of size $n\times n$ that 
contains zeros everywhere except at positions $(i_l,i_l)$ for 
$l=1,\ldots,k$ where it contains ones.
Note that one can additively decompose any given tangent vector
$V=P_{\ib}V+P_{\ib^\circ}V$. Hence we have
\begin{align*}
\Hb_Y&=\bigl\{P_{\ib^\circ}V \colon 
V\in\mathrm{T}_Y\Vb(n,k)\bigr\}\cong\CC^{(n-k)k}\,,\\
\Vb_Y&=\bigl\{P_{\ib}V \colon 
V\in\mathrm{T}_Y\Vb(n,k)\bigr\}\cong\mathfrak{gl}(k)\,.
\end{align*}

\subsection{Fibre bundle flow} 
Suppose we are given a vector field on the Stiefel manifold
\begin{equation*}
V(x,Y)=A(x,Y)\, Y, 
\end{equation*}
for any $(x,Y)\in\R\times\Vb(n,k)$ where 
$A\in\mathfrak{gl}(n)$. Fixing a 
coordinate patch $\Ub_\ib$ for $\Gr(n,k)$
for some $\ib=\{i_1,\ldots,i_k\}$ we can always decompose
$Y\in\Vb(n,k)$ into 
\begin{equation*}
Y=y_{\ib^\circ}u,
\end{equation*}
where $u\in\GL(k)$.
Let $a$, $b$, $c$ and $d$ denote
the $\ib\times\ib$, $\ib\times \ib^\circ$, $\ib^\circ\times \ib$
and $\ib^\circ\times \ib^\circ$ submatrices of $A$, respectively.

\begin{theorem}\label{th:coupled}
The flow governed by $V(x,Y)$ generates a coupled flow in
the base space $\Gr(n,k)$ and fibres $\GL(k)$: if $Y=y_{\ib^\circ}u$
for a given fixed $\ib$, the flows in the coordinate chart
variables $\hat y=\varphi_\ib\circ y_{\ib^\circ}$, and
rank $k$ transformations $u$, are 
\begin{equation*}
\hat y'=c+d\,\hat y-\hat y(a+b\,\hat y)
\qquad\text{and}\qquad
u'=(a+b\,\hat y)\,u,
\end{equation*}
where we now think of $a$, $b$, $c$ and $d$ 
as functions of $x$, $u$ and $\hat y$. 
\end{theorem}

\begin{proof}
Using that $Y=y_{\ib^\circ}u$ the ordinary differential
system $Y'=V(x,Y)$ becomes
\begin{equation*}
y_{\ib^\circ}'u+y_{\ib^\circ}u'
=(A_{\ib}+A_{\ib^\circ}\hat y)\,u,
\end{equation*}
where $A_{\ib}$ represents the submatrix obtained 
by restricting the matrix $A(x,y_{\ib^\circ}u)$ to its $\ib$th columns. 
Applying the projections $P_{\ib^\circ}$ and $P_{\ib}$ to both sides
of this equation, respectively generates the equations for 
$\hat y$ and $u$ shown. Note that $\hat y$ is the projection
of $y_{\ib^\circ}$ onto its $\ib^\circ$th rows, as well as its image
under the coordinate chart $\varphi_{\ib}$.
\end{proof}

\subsection{Riccati flow} 
A natural decoupling of the flow on the base space 
$\Gr(n,k)$ from the flow on the fibres $\GL(k)$ occurs
when the vector field $V$ is linear, i.e.\ when
\begin{equation*}
V(x,Y)=A(x)\,Y,
\end{equation*}
The following corollary is immediate from Theorem~\ref{th:coupled}. 

\begin{corollary}\label{cor:linVF}
If the vector field $V$ is linear so that $A=A(x)$ only, 
then:
\begin{enumerate}
\item The flow on the base space $\Gr(n,k)$
decouples from the flow evolving through the fibres 
$\GL(k)$---the flow in the fibres is slaved to that
in the base space;
\item For a fixed coordinate patch $\Ub_\ib$ index by $\ib$,
in the coordinate chart variables 
$\hat y=\varphi_\ib\circ y_{\ib^\circ}$,
the flow is governed by the Riccati differential system:
\begin{equation*}
\hat y'=c(x)+d(x)\hat y-\hat ya(x)-\hat yb(x)\hat y.
\end{equation*}
\end{enumerate}
\end{corollary}

Suppose we are required to determine the flow generated by 
a linear nonautonomous vector field defined on the 
Stiefel manifold $\Vb(n,k)$. The first conclusion in the corollary
implies that all the relevant information is carried in the
flow in the Grassmann manifold $\Gr(n,k)$, and the flow through
the fibres $\GL(k)$ can be completely determined a-posteriori 
from the Grassmannian flow. The second conclusion suggests
that if we fix a coordinate patch, then the flow in the 
Grassmannian can be determined from the solution to the 
Riccati system for $\hat y$. If required, we can solve
the differential system for $u$ in Theorem~\ref{th:coupled} 
to determine $Y$, thereby completely resolving the 
flow generated by the linear vector field $V$ on the Stiefel manifold.
Provided $\hat y$ remains finite, this approach in fact works.

The problem is that, though $Y=y_{\ib^\circ}u$ must be globally finite as 
it is generated by a linear vector field (with globally smooth coefficients), 
the Riccati solution $\hat y$ itself can become singular in a finite integration interval. 
Of course simultaneously the determinant of $u$ itself becomes zero.
The solution on the Grassmannian does not become singular.
The issue is representation. (Note that the flow on $\GL(k)$ is linear 
but rank is not preserved because its coefficients depend on the Riccati flow.)

The Riccati flow is a flow in a given fixed coordinate chart 
indexed by $\ib$, which is chosen at the start of integration.
Given an initial element in $\Vb(n,k)$, we pick a (good) coordinate patch
$\Ub_\ib$ for $\Gr(n,k)$, this fixes the Grassmannian representation
$y_{\ib^\circ}$. Projecting onto the $\ib^\circ$th rows of
$y_{\ib^\circ}$, or equivalently looking at the image under the
coordinate chart $\varphi_\ib$, generates $\hat y\in\CC^{(n-k)k}$.
The Riccati flow is the flow in the Euclidean chart image space
$\CC^{(n-k)k}$. Each coordinate patch $\Ub_\ib$ is dense in
$\Gr(n,k)$. Therefore in numerical computations, the Riccati
solution $\hat y$ is likely to leave and return to 
the patch $\Ub_\ib$ across any discrete integration step that
staddles a representation singularity (generating a large but 
finite solution $\hat y$ the other side). 

With this in mind, Schiff and Shnider~\cite{SS} suggested 
the following method that integrates through singularities 
in the Riccati flow. Fix a coordinate patch $\Ub_\ib$ with index $\ib$. 
The general linear group $\GL(n)$ acts transitively on $\Vb(n,k)$
(and also $\Gr(n,k)$): the left Lie group action 
$\Lambda\colon\GL(n)\times\Vb(n,k)\rightarrow\Vb(n,k)$ is defined by 
$\Lambda\colon(S,Y)\mapsto S\,Y$. For $Y_0\in\Vb(n,k)$ we set 
\begin{equation*}
\Lambda_{Y_0}\colon S\mapsto S\,Y_0.
\end{equation*}
The M\"obius Lie group action 
$\mu_{\hat y_0}\colon\GL(n)\ra\CC^{(n-k)k}$
is defined by
$\mu_{\hat y_0}\colon S\mapsto\varphi_{\ib}\circ\pi_\ib\circ
\Lambda_{\varphi_{\ib}^{-1}\circ\hat y_0}\circ S$,
where $\pi_\ib\colon\Vb(n,k)\ra\Ub_\ib$ is the quotient map 
$\pi_\ib\colon Y\mapsto y_{\ib^\circ}$. 
Explicitly, if $S_{\ib,\ib'}$ represents the $\ib\times\ib'$ submatrix of $S$,
we have:
\begin{equation*}
\mu_{\hat y_0}\colon S\mapsto 
(S_{\ib,\ib}+S_{\ib,\ib^\circ}\hat y_0)
(S_{\ib^\circ,\ib}+S_{\ib^\circ,\ib^\circ}\hat y_0)^{-1}.
\end{equation*}
Thus, given data $\hat y_0$ in the Euclidean chart image space 
$\CC^{(n-k)k}$, pullback to the Lie group $\GL(n)$, 
to the identity element $I_n$, using the M\"obius Lie group 
action map $\mu_{\hat y_0}$. 
Advance the solution across one integration step in the Lie group
generating the element $S\in\GL(n)$. Schiff and Shnider used a
Neumann/Runge--Kutta method to do this, but a Lie group method
could also be used. 
Now push forward to $\hat y=\mu_{\hat y_0}\circ S\in\CC^{(n-k)k}$ 
using the M\"obius Lie group action map. This takes you back to
the chart corresponding to the original coordinate patch $\Ub_\ib$.

This method integrates through singularities in the Riccati flow.  
However, we are still left with another associated practical problem. 
If the linear vector field $V$ depends on a parameter we wish to vary, 
the singularities in the Riccati flow can drift in the domain 
of integration and impinge on the matching position.

\section{Practical Grassmann integration}\label{sec:PGI}
Our goal in this section is construct numerical methods that
integrate the flow associated with the push forward of the linear
vector field $V$ onto the Grassmannian manifold, 
whilst avoiding the representation singularities 
that occur in the Riccati approach. The solution is to change patch 
continuously in some optimal fashion or whenever the coordinate patch 
becomes a poor representation. 

The following diagram helps map out the two new strategies we suggest.
\begin{equation*}
\xymatrixcolsep{3.5pc}
\xymatrix{
\CC^{(n-k)k}\ar[r]^-{\varphi_{\ib}^{-1}}  \ar@{->}[d]
&\Ub_{\ib}\ar[r]^-{\mathrm{id}} \ar@{->}[d]^-{\mathrm{GGEM}}
&\Vb(n,k)\ar[r]^-{(\Lambda_{Y_0})^{\ast}} \ar@{->}[d]
& \GL(n)\ar[r]^-{\log}  \ar@{->}[d]
&\mathfrak{gl}(n) \ar@{->}[d]_-{\mathrm{Magnus}}\\
\CC^{(n-k)k} 
&\ar[l]^-{\varphi_{\ib'}} \Ub_{\ib'}
&\ar[l]^-{\mathrm{QOGE}} \Vb(n,k)
&\ar[l]^-{\Lambda_{Y_0}} \GL(n)
&\ar[l]^-{\exp} \mathfrak{gl}(n) 
} 
\end{equation*}

\subsection{Continuous optimal patch evolution}
The idea behind the \emph{Grassmann Gaussian elimination method} 
(GGEM) is this. Given data in $\Gr(n,k)$ that lies in a given patch 
$Y_0\in\Ub_\ib$ identified by $\ib$, pullback to the Stiefel manifold
using the identity map---note $\Ub_\ib\subset\Vb(n,k)$---so $Y_0\in\Vb(n,k)$.
Advance the solution one integration step along the 
Stiefel manifold $\Vb(n,k)$, using say a classical Runge--Kutta method, 
thus generating the element $Y\in\Vb(n,k)$. 

Now with the next step solution on the Stiefel manifold, we
use \emph{quasi-optimal Gaussian elimination} with partial pivoting (QOGE) to 
decompose $Y=y_{(\ib')^\circ}u$ 
and project onto the coordinate patch $\Ub_{\ib'}$ producing the element 
$y_{(\ib')^\circ}\in\Ub_{\ib'}$. The quasi-optimal Gaussian elimination
process (described in Section~\ref{sec:NM}) picks out a suitable 
coordinate patch $\Ub_{\ib'}$ to represent the solution, 
which may be different than the original patch $\Ub_{\ib}$.

To ensure we numerically remain with the Stiefel manifold,
rather than use a Runge--Kutta method,  
we might use a Lie group method as follows 
(see Munthe--Kaas~\cite{MK}).
Pullback the data $Y_0\in\Vb(n,k)$ 
from the Stiefel manifold to the general linear
group $\GL(n)$, via the action map $\Lambda_{Y_0}$. 
The corresponding element in $\GL(n)$ is naturally
the identity element $I_n$. Subsequently pullback, 
via the exponential map, to the zero element in the 
corresponding Lie algebra, i.e.\ $\oo\in\mathfrak{gl}(n)$. 
Evolve the solution on the Lie algebra to 
$\sigma\in\mathfrak{gl}(n)$ using the Magnus expansion 
(Magnus~\cite{Ma}; also see Iserles, Munthe--Kaas, 
N\o rsett and Zanna~\cite{IMNZ}). 
Pushforward from $\mathfrak{gl}(n)$ to $\GL(n)$, via the 
exponential map, producing the Lie group element 
$S=\exp\sigma\in\GL(n)$. Now pushfoward to the Stiefel
manifold $\Vb(n,k)$ via the Lie group action map $\Lambda_{Y_0}$.
For more details on Lie group methods for Stiefel manifolds,
see Krogstad~\cite{Kr} and Celledoni and Owren~\cite{CeO}.

\subsection{Riccati flow with patch swapping}\label{subsec:RPS}
Since we are interested in constructing the flow generated 
by the push forward of the vector field $V$ onto the 
Grassmann manifold, we can avoid singularities in any
given Riccati system chart flow associated with a given
Grassmannian coordinate patch, by simply changing patch
when the solution representation appears to become poor.
In particular, we could either change the coordinate patch when the: 
\begin{itemize}
\item Norm $\|\hat y\|_{\infty}$ becomes too large 
(the easier and preferred approach we take);
\item Determinant $\det u$ becomes too small (this involves 
constructing $u$ as integration proceeds and we want to avoid 
carrying information unnecessarily).
\end{itemize}
Hence if say $\|\hat y\|_{\infty}$ exceeds a prescribed tolerance
at the end of one step, to change patch, 
we apply the quasi-optimal Gaussian elimination method (QOGE) 
to the matrix $y_{\ib^\circ}=\varphi_{\ib}^{-1}\circ \hat y$.
This identifies a new patch and index $\ib'$ to use for the 
next set of successive steps until
$\|\hat y\|_{\infty}$ becomes too large again, and so forth.

\subsection{Drury--Oja flow}
There is a close connection between the Riccati flow in
Corollary~\ref{cor:linVF} and the continuous orthogonalization 
method of Humpherys and Zumbrun~\cite{HZ}, proved by direct comparison.

\begin{lemma}\label{lem:DO}
If $\hat y\in\CC^{(n-k)k}$ satisfies the Riccati flow 
$\hat y'=c(x)+d(x)\hat y-\hat ya(x)-\hat yb(x)\hat y$ 
and $u\in\Ub(k)$ satisfies 
$u'=-\hat y^\dag(c+d\,\hat y)\,u$, then for any index $\ib$, 
we have that $Q=y_{\ib^\circ}u$
satisfies the Drury--Oja flow: $Q'=(I_n-QQ^\dag)AQ$.
\end{lemma}

Humpherys and Zumbrun~\cite{HZ} derive
this flow by a $QR$-decomposition of the solution 
$Y=QR\in\Vb(n,k)$ to the linear, spectral, globally bounded 
flow $Y'=A(x)\,Y$. 
We can think of this continuous orthogonalization method 
as generating an approximate flow on the Grassmann manifold
whilst evolving the coordinatization (see Edelman, Arias and Smith~\cite{EAS};
Bindel, Demmel and Friedman~\cite{BDF}). 
It is also known as a Drury--Oja flow (see Drury~\cite{D}, 
Oja~\cite{Oja}, Yan, Helmke and Moore~\cite{YHM}, 
Bridges and Reich~\cite{BR}, Dieci and Van Vleck~\cite{DVV}
and Hairer, Lubich and Wanner~\cite[p.~136]{HLW}). 
The determinant of the upper triangular matrix,
$\det R$, grows exponentially. Thus, since it is nonzero in the far field, 
we know it remains nonzero in the whole integration interval $\mathbb R$.
Consequently, $Q=YR^{-1}$ is globally finite on $\mathbb R$,
i.e.\ there are no singularities in the Drury--Oja flow.

There is a natural $\mathfrak{su}(n)$ Lie algebra action 
on the Stiefel manifold of orthonormal $k$-frames, $\Vb_0(n,k)$, 
generated by the map $(\sigma,Q)\mapsto\exp(\sigma)\,Q$,
for $Q\in\Vb_0(n,k)$ and $\sigma\in\mathfrak{su}(n)$.
If $v\circ Q\equiv (I-QQ^\dag)\,A$,
the flow on the Lie algebra $\mathfrak{su}(n)$ that generates
the Drury--Oja flow on $\Vb_0(n,k)$ is governed by
$\sigma'=\mathrm{dexp}_\sigma^{-1}\circ v\bigl(\exp(\sigma)Q_0\bigr)$.
We could use this to construct a numerical method that preserves 
orthonormality of $Q$.

\section{Grassmann Gaussian elimination method}\label{sec:NM}

\subsection{Algorithm}
The \emph{Grassmann Gaussian elimination method} using
a Lie group method on the Stiefel manifold (GGEM-LG), 
proceeds as follows:
\begin{enumerate}
\item Suppose initially we are given data $y_{\ib^\circ_m}(x_m)$
in the coordinate patch $\Ub_{\ib_m}$.
\item Across the integration interval $[x_m,x_{m+1}]$, compute $\sigma_m$ 
using the Magnus expansion---we recommend the fourth order Magnus expansion.
\item Compute $Y_{m+1}=\exp\sigma_m\cdot y_{\ib^\circ_m}(x_m)$.
\item Apply quasi-optimal Gaussian elimination (QOGE)
(outlined below) to $Y_{m+1}$; 
this generates the solution $y_{\ib^\circ_{m+1}}(x_{m+1})$ 
in the coordinate chart $\Ub_{\ib_{m+1}}$, and the rank $k$ matrix
$U_{m+1}$ we have effectively equivalenced by.
\end{enumerate}
Across the integration interval $[x_m,x_{m+1}]$, we could
generate $Y_{m+1}$ from $y_{\ib^\circ_m}(x_m)$ by solving the flow 
on the Stiefel manifold $\Vb(n,k)$ using a classical Runge--Kutta step
(GGEM-RK). Both algorithms are summarized in the following diagram. 
\begin{equation*}
\xymatrixcolsep{3.5pc}
\xymatrix{
y_{\ib^\circ_m}(x_m)\ar[r]^-{\mathrm{id}} \ar@{->}[d]^-{\mathrm{GGEM}}
&y_{\ib^\circ_m}(x_m)
\ar[r]^-{\bigl(\Lambda_{y_{\ib^\circ_m}(x_m)}\bigr)^{\ast}}\ar@{->}[d]^-{\mathrm{RK}}
& I_n\ar[r]^-{\log}  
&\oo \ar@{->}[d]_-{\mathrm{Magnus}}\\
y_{\ib^\circ_{m+1}}(x_{m+1})
&\ar[l]^-{\mathrm{QOGE}} Y_{m+1}
&\ar[l]^-{\Lambda_{y_{\ib^\circ_m}(x_m)}} S_{m}
&\ar[l]^-{\exp} \sigma_{m}
} 
\end{equation*}
How much of $U_{m+1}$ do we retain at each step? 
We address this in Section~\ref{subsec:GGEMmatching}.

\subsection{Quasi-optimal Gaussian elimination}
To project the endpoint solution $Y_{m+1}$ in the Stiefel manifold
onto a quasi-optimal coordinate patch, say $\Ub_{\ib_{m+1}}$, 
we use a Gauss--Jordan approach. This entails using elementary 
column operations and \emph{optimal pivoting}, 
until a specified $k\times k$ submatrix of $Y_{m+1}$ becomes
the identity matrix (we naturally assume $k\geq2$). 
To be more explicit about what we mean
by optimal pivoting, we outline the procedure:
\begin{enumerate}
\item Look for the largest term in magnitude in $Y^{(1)}:=Y_{m+1}$.
Nominate this term as the pivot term and suppose this occurs in row $i_1$, 
column $j_1$. Use elementary column operations with this term, 
$p_{i_1j_1}:=Y_{i_1,j_1}^{(1)}$, 
as the pivot to render all the other $(k-1)$ terms in that row equal 
to zero. Finally, use the elementary column operation of 
scalar multiplication to render the pivot term itself equal to one. 
Let us call the resulting $n\times k$ matrix $Y^{(2)}$; we can write
(see for example Meyer~\cite{Me})
\begin{equation*}
Y^{(2)}=Y^{(1)}\,U^{(1)},
\end{equation*}
where the elementary column operations performed are 
encoded in the elementary matrix $U^{(1)}$. To be precise,
we can write
\begin{equation*}
U^{(1)}=E_1^{(1)}\cdots E_{j_1-1}^{(1)} E_{j_1+1}^{(1)}
\cdots E_{k-1}^{(1)}E_{j_1}^{(1)},
\end{equation*}
where for $\ell\in\{1,\ldots,k\}/\{j_1\}$ the elementary matrix
$E_\ell^{(1)}$ encodes the following elementary column operation
on column $c_\ell$:
\begin{equation*}
c_\ell\to c_\ell-\frac{Y^{(1)}_{i_1,\ell}}{p_{i_1j_1}}c_{j_1}.
\end{equation*}
For $\ell=j_1$ the elementary matrix $E_{j_1}^{(1)}$ encodes the 
elementary column operation $c_{j_1}\to c_{j_1}/p_{i_1j_1}$. 
In summary, the matrix $Y^{(2)}$ has value
one at position $(i_1,j_1)$ and otherwise zeros in row $i_1$;
and we have $\det U^{(1)}=(p_{i_1j_1})^{-1}$.
\item In the $(n-1)\times(k-1)$ submatrix of $Y^{(2)}$ identified by 
excluding row $i_1$ and column $j_1$, look for the largest term in magnitude.
Again nominate this term as the pivot term and suppose this occurs in row $i_2$
(which will be distinct from $i_1$), and column $j_2$ (which is distinct 
from $j_1$). Here $i_2$ and $j_2$ refer to the row and column
relative to the original $n\times k$ matrix $Y^{(2)}$.
Use elementary column operations with this term, $p_{i_2j_2}:=Y_{i_2,j_2}^{(2)}$, 
as the pivot to render the terms in row $i_2$ and columns
$\{1,\ldots,k\}/\{j_1,j_2\}$ of $Y^{(2)}$ equal to zero. 
Again use the elementary column operation of 
scalar multiplication to render the pivot term itself equal to one. 
Let us call the resulting $n\times k$ matrix $Y^{(3)}$; we can write
\begin{equation*}
Y^{(3)}=Y^{(2)}\,U^{(2)}=Y^{(1)}\,U^{(1)}\,U^{(2)},
\end{equation*}
where the elementary column operations performed on $Y^{(2)}$ are 
encoded in the elementary matrix $U^{(2)}$. Again, to be precise,
we can write
\begin{equation*}
U^{(2)}=E_1^{(2)}\cdots E_{j_1-1}^{(2)} E_{j_1+1}^{(2)}
\cdots E_{j_2-1}^{(2)} E_{j_2+1}^{(2)}
\cdots E_{k-1}^{(2)}E_{j_2}^{(2)},
\end{equation*}
where for $\ell\in\{1,\ldots,k\}/\{j_1,j_2\}$ the elementary matrix
$E_\ell^{(2)}$ encodes the following elementary column operation
on column $c_\ell$:
\begin{equation*}
c_\ell\to c_\ell-\frac{Y^{(2)}_{i_2,\ell}}{p_{i_2j_2}}c_{j_2}.
\end{equation*}
Note that in our expression for $U^{(2)}$ above we could have
either $j_1<j_2$ or $j_1>j_2$. For $\ell=j_2$ then $E_{j_2}^{(2)}$ 
encodes the elementary column operation $c_{j_2}\to c_{j_2}/p_{i_2j_2}$. 
In summary, the resulting matrix $Y^{(3)}$ has value
one at positions $(i_1,j_1)$ and $(i_2,j_2)$ and otherwise zeros
in row $i_1$, and zeros in row $i_2$ in columns 
$\{1,\ldots,k\}/\{j_1,j_2\}$; 
and we have $\det U^{(2)}=(p_{i_2j_2})^{-1}$.
\item Continue this process. Focus on the $(n-2)\times(k-2)$ 
submatrix of $Y^{(3)}$ identified by excluding the rows $i_1, i_2$
and columns $j_1, j_2$; look for the largest term in magnitude.
Nominate this term as the pivot term---suppose it occurs in 
row $i_3$ and column $j_3$, relative to the original $n\times k$
matrix $Y^{(3)}$, and so forth. On completing the final $k$th step
in this process, we will have 
\begin{equation*}
Y^{(k)}=Y^{(1)}\,U^{(1)}\,U^{(2)}\,\cdots\,U^{(k)},
\end{equation*}
where
\begin{equation*}
\det\bigl(U^{(1)}\,U^{(2)}\,\cdots\,U^{(k)}\bigr)
=(p_{i_1j_1}p_{i_2j_2}\cdots p_{i_kj_k})^{-1}.
\end{equation*}
The final $n\times k$ matrix $Y^{(k)}$ will have ones 
in positions $(i_\ell,j_\ell)$, for $\ell=1,\ldots,k$.
In row $i_\ell$ it will have zeros in columns 
$\{1,\ldots,k\}/\{j_1,\ldots,j_\ell\}$. We set
\begin{equation*}
\ib_{m+1}:=\{i_1,i_2,\ldots,i_k\}.
\end{equation*}
\item Perform column swaps in $Y^{(k)}$ so that column $j_1$
becomes column $1$, column $j_2$ becomes column $2$, and so forth
so that finally column $j_k$ is forced to be column $k$. These
column swaps can be encoded in the elementary matrix $\Sigma$
with $\det\Sigma=(-1)^{\#\{\mathrm{swaps}\}}$. The resulting matrix
$\tilde Y^{(k)}$ is given by
\begin{equation*}
\tilde Y^{(k)}=Y^{(k)}\,\Sigma.
\end{equation*}
The $\ib_{m+1}\times\{1,\ldots,k\}$ submatrix of $\tilde Y^{(k)}$ 
given by 
\begin{equation*}
L:=\bigl[\tilde Y^{(k)}\bigr]_{\ib_{m+1}\times\{1,\ldots,k\}}
\end{equation*}
is lower triangular with ones on the diagonal. Finally we set
\begin{equation*}
y_{\ib_{m+1}^\circ}:=\tilde Y^{(k)}\,L^{-1}.
\end{equation*}
In practice of course, we do not compute $L^{-1}$, but continue
performing elementary column operations on $\tilde Y^{(k)}$ so that 
the submatrix $L$ becomes the identity matrix, 
thus generating $y_{\ib_{m+1}^\circ}$.
Hence we have effectively performed the decomposition
\begin{equation*}
Y_{m+1}=y_{\ib_{m+1}^\circ}\cdot U_{m+1},
\end{equation*}
where 
$U_{m+1}=L\,\Sigma^{-1}\bigl(U^{(k)}\bigr)^{-1}\cdots\bigl(U^{(1)}\bigr)^{-1}\in\GL(k)$,
and in particular
\begin{equation*}
\det U_{m+1}=(-1)^{\#\{\mathrm{swaps}\}}\,p_{i_1j_1}\cdots p_{i_kj_k}.
\end{equation*}
\end{enumerate}

Note that if $Y_{m+1}\in\Vb(n,k)$ then the entries 
in $y_{\ib^\circ_{m+1}}$ and $U_{m+1}$, produced as a result of this process,
will all be finite. Further if $Y_{m+1}$ depends analytically on a
parameter, then the product $y_{\ib_{m+1}^\circ}\cdot U_{m+1},$ naturally does as well.
Lastly we remark that we could have performed alternative elementary column operations
of the form $c_\ell\to p_{ij}c_\ell-Y_{i,\ell}^{(\cdot)}c_j$ (we do not
include the scalar multipication operations here) with the result
that $\det U_{m+1}=(-1)^{\#\{\mathrm{swaps}\}}\,(p_{i_1j_1})^{2-k}
(p_{i_2j_2})^{3-k}\cdots (p_{i_kj_k})^1$.

\subsection{Complexity}
The complexity of the quasi-optimal Gaussian elimination algorithm,
dominated by the search for the largest elements in the successively
decreasing submatrices of $Y_{m+1}$, is of order $nk^2$.
The method is a practical approach to maximize the determinant of the 
$k\times k$ submatrix removed from $Y_{m+1}$. It will not in 
general choose the submatrix with the largest determinant---hence the label
\emph{quasi-optimal}. This could be achieved by searching through all the 
$k\times k$ submatrices of $Y_{m+1}$, i.e.\ all the Pl\"ucker coordinates, 
but this has complexity of order $n$ choose $k$. An interesting
question here is whether there is an efficient way to use the 
Pl\"ucker relations to reduce this complexity?

\section{Spectral problems}\label{sec:SP}

\subsection{Linear Stiefel flow}
Consider the linear spectral problem on $\R$:
\begin{equation*}
Y'=A(x;\lambda)\,Y 
\end{equation*}
We assume there exists a subdomain $\Omega\subseteq\CC$
containing the right-half complex plane, such that
for $\lambda\in\Omega$ there exists exponential dichotomies 
on $\R^-$ and $\R^+$ with the same Morse index $k$ in each case 
(see Henry~\cite{H} and Sandstede~\cite{S}). 
Let $Y^-(x;\lambda)\in\Vb(n,k)$ denote 
the matrix whose columns are solutions to the spectral problem 
and which span the unstable manifold section at $x\in[-\infty,+\infty)$.
Let $Y^+(x;\lambda)\in\Vb(n,n-k)$ denote the matrix 
whose columns are the solutions which span the stable 
manifold section at $x\in(-\infty,+\infty]$.

\subsection{Matching}
The values of spectral parameter $\lambda\in\Omega$ for which the 
columns of $Y^-$ and columns of $Y^+$ are linearly dependent
on $\R$ are pure-point eigenvalues.  
The \emph{Evans function} $D(\lambda)$ is the measure of the degree
linear dependence between the two basis sets $Y^-$ and $Y^+$, i.e.\  
of the degree of transversal intersection between the 
unstable and stable manifolds (see Alexander, Gardner and Jones~\cite{AGJ}; 
Nii~\cite{Nii}): 
\begin{equation*}
D(\lambda)\equiv 
\mathrm{e}^{-\int_0^x\mathrm{Tr}A(\xi;\lambda)\,\mathrm{d}\xi}\,
\mathrm{det}\bigl(Y^-(x;\lambda)\,\, Y^+(x;\lambda)\bigr).
\end{equation*}
It is analytic in $\Omega$.
In practice we drop the non-zero, scalar exponential 
prefactor and evaluate the Evans function at a matching point $x_\ast$.

There are other matching criteria measuring the degree of intersection
between subspaces that do not use the determinant: for example
computing the angle between subspaces 
as suggested by Bj\"orck and Golub~\cite{BG} or computing
the smallest eigenvalue as suggested by Hutson \cite{Hu} and Ixaru~\cite{Ixtalk}. 
Both these latter techniques might be important for 
large systems when computing the determinant could be
an unstable process, indeed, we investigate them in
this context in Ledoux, Malham, Niesen and Th\"ummler~\cite{LMNT}.
However in both cases the magnitude of a function of the 
spectral parameter is computed, whose zeros correspond to
eigenvalues. Hence we must search for touchdowns to zero 
in the complex parameter spectral plane which can be
problematic. For the examples we consider here, which
are not too large, using the determinant suffices.

\subsection{Initialization}
We construct the $n\times k$ matrix $Y_0^-(\lambda)$ 
whose columns are the $k$ eigenvectors of $A(-\infty;\lambda)$ 
corresponding to eigenvalues with a positive real part
(see Humpherys and Zumbrun~\cite{HZ} and 
also Humpherys, Sandstede and Zumbrun~\cite{HSZ} 
for how to preserve analyticity
with respect to the spectral parameter $\lambda$).
In practice, integration starts at $x=\ell_-$ for some suitable, 
usually negative, value of $\ell_-$.
Analogously we construct the $n\times(n-k)$ matrix 
$Y_0^+(\lambda)$ whose columns are the $n-k$ eigenvectors 
of $A(+\infty;\lambda)$ 
corresponding to eigenvalues with a negative real part.
Again, in practice, we integrate backwards from $x=\ell_+$ 
for some suitable, usually large and positive, value of $\ell_+$. 

\subsection{GGEM matching and analyticity}\label{subsec:GGEMmatching}
To compute $Y^\pm(x_\ast;\lambda)$ we start with $Y_0^\pm(\lambda)$
at $x=\ell_\pm$ and integrate centrally towards $x=x_\ast$. Our
goal in this section is to show that we only need the determinant
of the rank $k$ transformations in the 
GGEM method described at the beginning of Section~\ref{sec:NM}, 
to retain analyticity for the Evans function in $\Omega$.
Since the procedure is the same in both intervals $[\ell_\pm,x_\ast]$
we will describe it for the generic interval $[\ell,x_\ast]$ for
$Y(x;\lambda)\in\Vb(n,k)$ starting with value $Y_0(\ell;\lambda)$ 
at $x=\ell$. Suppose we use $M$ successive computation subintervals
$[x_m,x_{m+1}]$ in $[\ell,x_\ast]$ where $x_m:=\ell+m\,(x_\ast-\ell)/M$.
We label the nodal solution values at $x=x_m$ as $Y_m(\lambda)$.

At the start $x_0=\ell$, perform quasi-optimal Gaussian elimination
(QOGE) on $Y_0(\lambda):=Y_0(\ell;\lambda)$ to obtain the decomposition
\begin{equation*}
Y_0(\lambda)=y_{\ib_0^\circ}(x_0;\lambda)\,U_0(\lambda).
\end{equation*}
As we shall see, we do not need to actually store $U_0(\lambda)$,
but only $\det U_0(\lambda)$. 

Let $S_{m,m+1}(\lambda)$ denote an approximation to the flow-map across
$[x_m,x_{m+1}]$ to the linear system $Y'=A(x;\lambda)\,Y$.
We assume that $S_{m,m+1}(\lambda)$ preserves analytic dependency
on $\lambda$, so that the next step solution value 
$Y_{m+1}(\lambda)=S_{m,m+1}(\lambda)\,Y_m(\lambda)$ analytically
depends on $\lambda$ if $Y_m(\lambda)$ does.
For example, in the case of GGEM-LG then 
$S_{m,m+1}(\lambda)=\exp\sigma_m$ and most straightforward Magnus
based integrators will naturally preserve analyticity with respect 
to $\lambda$. Similarly most simple Runge--Kutta methods used to generate
$S_{m,m+1}(\lambda)$, or directly generate the next step solution value
$Y_{m+1}(\lambda)$, will preserve analyticity.

Our numerical procedure would proceed as follows. Across $[x_0,x_1]$
we have
\begin{align*}
S_{0,1}(\lambda)\,Y_0(\lambda)&=\bigl(S_{0,1}(\lambda)\,
y_{\ib_0^\circ}(x_0;\lambda)\bigr)\,U_0(\lambda)\\
&=y_{\ib_1^\circ}(x_1;\lambda)\,U_1(\lambda)\,U_0(\lambda),
\end{align*}
where in the last step we applied QOGE to the $n\times k$ matrix 
$S_{0,1}(\lambda)\,y_{\ib_0^\circ}(x_0;\lambda)$. Subsequently,
across $[x_1,x_2]$ we have
\begin{align*}
S_{1,2}(\lambda)\,\bigl(y_{\ib_1^\circ}(x_1;\lambda)
\,U_1(\lambda)\,U_0(\lambda)\bigr)&=\bigl(S_{1,2}(\lambda)\,
y_{\ib_1^\circ}(x_0;\lambda)\bigr)\,U_1(\lambda)\,U_0(\lambda)\\
&=y_{\ib_2^\circ}(x_2;\lambda)\,U_2(\lambda)\,U_1(\lambda)\,U_0(\lambda),
\end{align*}
where we applied QOGE to 
$S_{1,2}(\lambda)\,y_{\ib_1^\circ}(x_0;\lambda)$. Repeating this
argument across the subsequent intervals $[x_m,x_{m+1}]$ for
$m=2,\ldots,M-1$, we get the following approximation 
$\hat Y(x_\ast;\lambda)$ to $Y(x_\ast;\lambda)$:
\begin{align*}
\hat Y(x_\ast;\lambda)&=S_{M-1,M}(\lambda)\,\cdots S_{0,1}(\lambda)\,Y_0(\lambda)\\
&=y_{\ib^\circ_M}(x_M;\lambda)\,U_M(\lambda)\,\cdots\,U_0(\lambda)\\
&=y_{\ib^\circ_M}(x_M;\lambda)\,U(\lambda),
\end{align*}
where $U(\lambda):=U_M(\lambda)\,\cdots\,U_0(\lambda)$. Naturally
the product $y_{\ib^\circ_M}(x_M;\lambda)\,U(\lambda)$ depends analytically
on $\lambda$.

Returning to the separate interval calculations on $[\ell_\pm,x_\ast]$,
we see that by the procedure just outlined we can generate the 
solution approximations 
\begin{equation*}
\hat Y_\pm(x_\ast;\lambda)=y_{\ib^\circ_M(\pm)}^\pm(x_\ast;\lambda)\,U^\pm(\lambda).
\end{equation*}
The number of integration steps $M$ can of course be different in each
interval.
Hence, after dropping the exponential prefactor and fixing the matching
point to be $x_\ast$, the Evans function $D(\lambda;x_\ast)$ can
be approximated by $\hat D(\lambda;x_\ast)$ where
\begin{align*}
\hat D(\lambda;x_\ast):=&\;
\det\bigl(\hat Y_-(x_\ast;\lambda) \,\,\hat Y_+(x_\ast;\lambda)\bigr)\\
=&\;\det\Bigl(y_{\ib^\circ_M(-)}^-(x_\ast;\lambda) \,\,
y_{\ib^\circ_M(+)}^+(x_\ast;\lambda)\Bigr)\cdot
\det U^-(\lambda)\cdot\det U^+(\lambda).
\end{align*}
This is an analytic function of $\lambda$. Indeed as we hinted
previously, at each computation step $x=x_m$ we need only
store $y_{\ib^\circ_m}(x_m;\lambda)$ and the value
\begin{equation*}
\prod_{\ell=m}^0\det U_\ell(\lambda)
\end{equation*}
which gets updated at each step by simply multiplying the 
previous step value by the complex scalar determinental factor 
for current step. Hence to preserve analyticity for the
Evans function using GGEM we must generate an approximate 
flow on $\Gr(n,k)\otimes\mathbb C$.

\subsection{Scaled GGEM}\label{subsec:FPR}
The scalar determinental factor just described, that
we update at each step, grows exponentially. 
This would be tempered by
the scalar exponential prefactor in the definition of 
the Evans function. An accurate practical procedure here is 
as follows (to be applied with due care). 
When integrating in the interval $[\ell_-,x_\ast]$,
at each step, divide the scalar determinental factor in GGEM 
by $\exp\bigl((\mu_1^-(\lambda)+\cdots+\mu_k^-(\lambda))h\bigr)$,
where $h$ is the stepsize, and the $\mu_i^-(\lambda)$ are
the (spatial) eigenvalues, with positive real part, of 
$A(-\infty;\lambda)$. When integrating in the interval 
$[\ell_+,x_\ast]$, at each step, divide the scalar factor
by $\exp\bigl(-(\mu_1^+(\lambda)+\cdots+\mu_{n-k}^+(\lambda))h\bigr)$,
where the $\mu_i^+(\lambda)$ are the eigenvalues, 
with negative real part, of $A(+\infty;\lambda)$.

To see that this normalization is appropriate,
we recall the Pl\"ucker coordinates of Section~\ref{sec:GS}. 
After applying the optimal Gaussian elimination algorithm
to $Y_{m+1}$ the $\ib^\circ_{m+1}$th row elements of $y_{\ib^\circ_{m+1}}$
are themselves $(n-k)k$ Pl\"ucker coordinates; normalized
by $\det U_{m+1}$. The $\ib^\circ_{m+1}$th row elements of $y_{\ib^\circ_{m+1}}$
and $\det U_{m+1}$, can be used to reconstruct the remaining 
Pl\"ucker coordinates through the homogeneous, quadratic 
Pl\"ucker relations. Hence the Pl\"ucker coordinates, or complete 
set of $k\times k$ minors, of $Y_{m+1}$ and $y_{\ib^\circ_{m+1}}$ differ
by a factor $\det U_{m+1}$. It is well known that if the original
vector field on $\Vb(n,k)$ is linear, then the
Pl\"ucker coordinates corresponding 
to $Y_{m+1}$ satisfy a (larger) linear system of equations. 
In the left far field the coordinates thus grow exponentially, in fact with 
growth rate $\mu_1^-(\lambda)+\cdots+\mu_k^-(\lambda)$; hence our
recommendation to divide by the exponential factor suggested
(with an analogous argument for the right far field).
See Alexander, Gardner and Jones~\cite{AGJ},
Alexander and Sachs~\cite{AS}, Brin~\cite{Br,Brp}
or Allen and Bridges~\cite{AB} for more details.

\section{Applications}\label{sec:ES}
We present some numerical results for three different applications. 
The three applications reduce to the solution of a system showing 
multiple distinct exponential growth and decay rates in the 
stable and unstable subspaces, respectively.
We show that our approach resolves this numerical obstacle 
successfully and can compete with the continuous orthogonalization method 
of Humpherys and Zumbrun~\cite{HZ}.

\subsection{Algorithms}
We implement six different algorithms as follows.

\emph{(1) Riccati-RK:} Riccati method with fixed coordinatization
with the flow of the Riccati vector field approximated by the 
classical fourth order Runge--Kutta method. We generically
chose the coordinatization labelled by $\ib^-=\{1,\ldots,k\}$ and 
$\ib^+=\{k+1,\ldots,n\}$ for the left-hand and right-hand
intervals, respectively. Hence if $Y_0^-(\lambda)$ and $Y_0^+(\lambda)$
denote the unstable and stable subspaces of $A(-\infty;\lambda)$ 
and $A(+\infty;\lambda)$, respectively, then we set 
\begin{equation*}
\hat y^\pm(\ell_\pm;\lambda)=\bigl(Y_0^\pm(\lambda)\bigr)_{(\ib^\pm)^\circ,\ib^\pm}
\bigl(Y_0^\pm(\lambda)\bigr)^{-1}_{\ib^\pm,\ib^\pm},
\end{equation*}
where $(Y)_{\ib,\ib'}$ denotes the $\ib\times\ib'$ submatrix of $Y$.
We integrate the Riccati equation outlined in Corollary~\ref{cor:linVF}
in the two intervals and evaluate the modified Evans function
\begin{equation*}
D(\lambda;x_\ast)\equiv 
\det\begin{pmatrix} I_k & \hat y^+(x_\ast;\lambda) \\ 
\hat y^-(x_\ast;\lambda) & I_{n-k}\end{pmatrix}.
\end{equation*}
Provided neither Riccati flow becomes singular, this 
Evans function is analytic in the spectral parameter $\lambda$.

\emph{(2) M\"obius--Magnus:} Uses the Schiff and Shnider approach to
integrate through singularities, combined with a Lie group method
to advance the solution on the general linear group, as described 
at the end of Section~\ref{sec:VF}. The same generic fixed coordinate
charts are used as for the Riccati-RK method above.
Over the integration interval $[x_m,x_{m+1}]$ with an equidistant 
mesh stepsize $h$, we advance the solution on the Lie algebra 
using the fourth order Magnus method 
\begin{equation*}
\sigma_m=\tfrac{1}{2}h\bigl(A(x_m^{[1]})+A(x_m^{[2]})\bigr)
-\tfrac{\sqrt{3}}{12}h^2\bigl[A(x_m^{[1]}),A(x_m^{[2]})\bigr],
\end{equation*}
with the two Gauss--Legendre points 
(see Iserles, Marthinsen and N\o rsett~\cite{IMN})
\begin{equation*}
x_m^{[1]}=x_m+(\tfrac{1}{2}-\tfrac{1}{6}\sqrt{3})h
\quad\text{and}\quad 
x_m^{[2]}=x_m+(\tfrac{1}{2}+\tfrac{1}{6}\sqrt{3})h.
\end{equation*}
We then compute the M\"obius map $\hat y_{m+1}=\mu_{\hat y_m}\circ\exp\sigma_m$
to advance the solution in the fixed Grassmannian chart---for the 
left-hand interval $\ib^-=\{1,\ldots,k\}$ while for the right-hand
interval $\ib^+=\{k+1,\ldots,n\}$. We evaluate the same Evans function
as for the Riccati-RK method above.

\emph{(3) GGEM-RK:} Scaled Grassmann Gaussian elimination method, with the 
classical fourth order Runge--Kutta method used to advance the solution on the 
Stiefel manifold, as described in Sections~\ref{sec:NM} and \ref{subsec:FPR}. 
We evaluate the Evans function $\hat D(\lambda;x_\ast)$ in Section~\ref{sec:SP}.

\emph{(4) GGEM-LG:} Same as GGEM-RK but with a fourth order Magnus method 
used to advance the solution on the Stiefel manifold instead, i.e.\ 
$Y_{m+1}=\exp(\sigma_m)\,y_{\ib_m}$ where $\sigma_m$ is generated
as for the M\"obius--Magnus method above. 

\emph{(5) Riccati-QOGE:} Riccati method with coordinate swapping as described 
in Section~\ref{subsec:RPS}. We have chosen to implement the method in the
following form. At each 
integration step we advance the solution on the Stiefel manifold using the 
Magnus method (we could also use a Runge--Kutta method here). 
We apply elementary column operations to the resulting solution matrix
to convert the pre-determined rows indexed by $\ib$ from the previous step
to the identity matrix. Then if $\|\hat y\|_\infty$ is less than or equal
to a tolerance size, we keep this index $\ib$ for the next step. If 
it is greater, we apply QOGE at the end of the next step after advancing the
solution on the Stiefel manifold, thus generating a new index. 
As for GGEM-RK and GGEM-LG, we update the scalar determinental factor at each step
(produced by the elementary column operations with the pre-determined index
or QOGE). We divide the scalar determinental factor by 
the scalar exponential factors, as described for the scaled GGEM method.
We evaluate the same Evans function also.

\emph{(6) CO-RK:} Continuous orthogonalization method 
of Humpherys and Zumbrun with the classical fourth order Runge--Kutta method 
used to advance the solution on the Stiefel manifold of orthonormal frames. 
The initial conditions $Q_0^\pm(\lambda)$ for the $Q$-matrices are obtained by 
QR-factorization of $Y_0^\pm(\lambda)$. From Humpherys and Zumbrun~\cite{HZ},
to ensure analyticity we must also solve the scalar problems
$(\det R^\pm)'=\text{Tr}\bigl(Q^\dag A(x;\lambda)Q
-(Q^\pm_0(\lambda))^\dag A(\pm\infty;\lambda)Q^\pm_0(\lambda)\,x\bigr)\det R^\pm$. 
The Evans function is then given by 
\begin{equation*}
D(\lambda;x_\ast)=
\det R^-(x_\ast;\lambda)\cdot
\det R^+(x_\ast;\lambda)\cdot
\det\bigl(Q^-(x_\ast;\lambda)\,\,Q^+(x_\ast;\lambda)\bigr).
\end{equation*}

\subsection{Boussinesq system}
As the first test system, we consider the   
Boussinesq system studied by 
Humpherys and Zumbrun~\cite{HZ}.
The (good) Boussinesq equation,
expressed in a co-moving frame moving to the right with
wave speed $c$, is given by    
\begin{equation*}
u_{tt}=(1-c^2)\,u_{xx}+2c\,u_{xt}-u_{xxxx}-(u^2)_{xx}.
\end{equation*}
It has solitary wave solutions of the form
$\bar{u}(x)\equiv \tfrac32(1-c^2)\mathrm{sech}^2
\bigl(\tfrac12\sqrt{1-c^2}\,x\bigr)$, 
where $|c|<1$. These waves are stable 
when $1/2<|c|<1$ and unstable when $|c|<1/2$.

If we consider small perturbations about the travelling wave 
$\bar u$ we generate a linear spectral problem 
of the form $Y'=A(x;\lambda)Y$, where
\begin{equation*}
A(x;\lambda)=\left(\begin{matrix}0&1&0&0\\0&0&1&0\\0&0&0&1\\
-\lambda^2-2{\bar u}''&2\lambda c-4{\bar u}'&
(1-c^2)-2{\bar u}&0\end{matrix}\right).
\end{equation*}
When the spectral parameter $\lambda$ lies in the right-half
complex plane the eigenvalues of $A(\pm\infty;\lambda)$
spectrally separate into two growth and two decay modes, 
i.e.\ $k=2$. We used $\ell_\pm=\pm8$ in our experiments.

\begin{figure}
\begin{center}
		\includegraphics[width=0.45\textwidth]{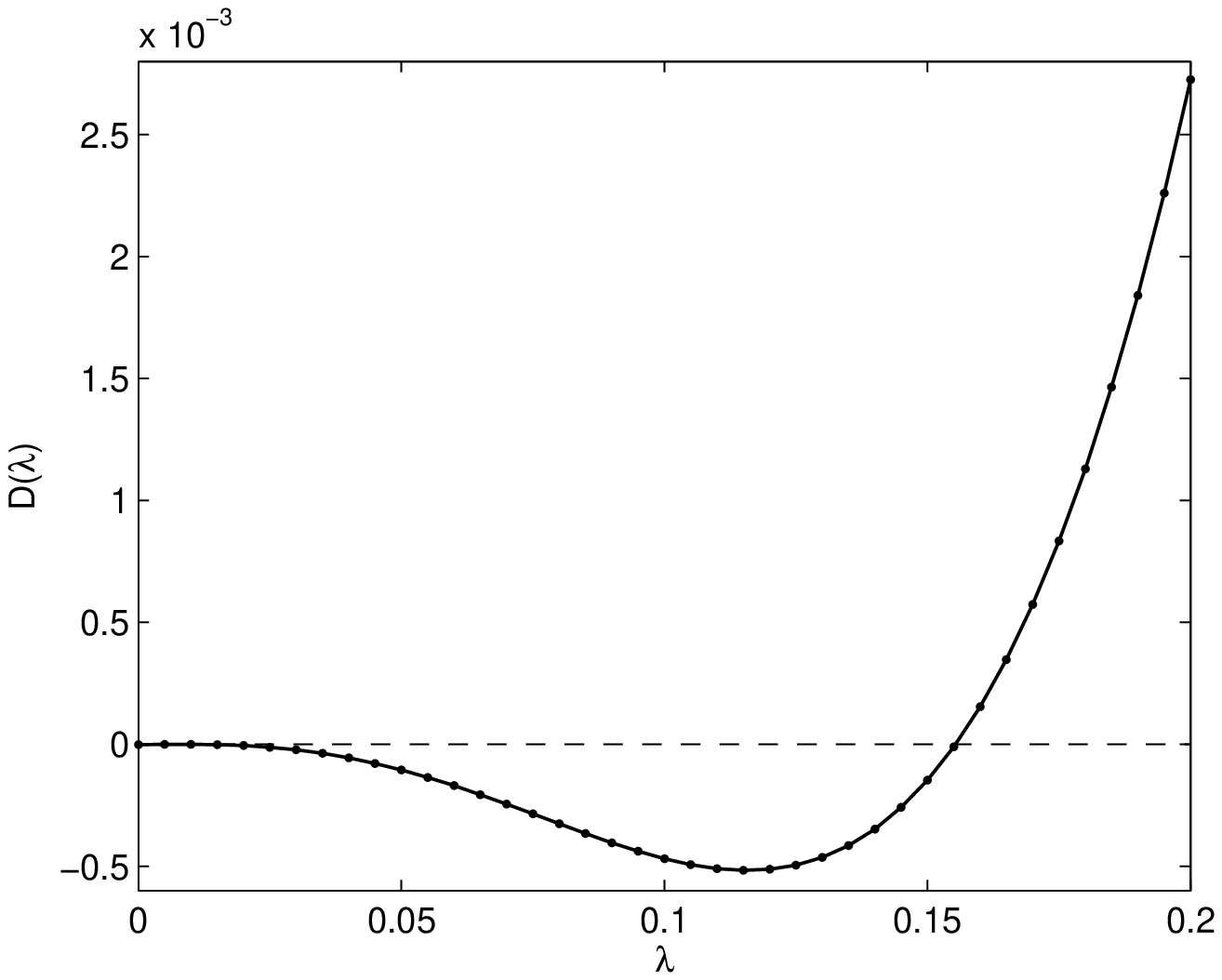}
		\includegraphics[width=0.45\textwidth]{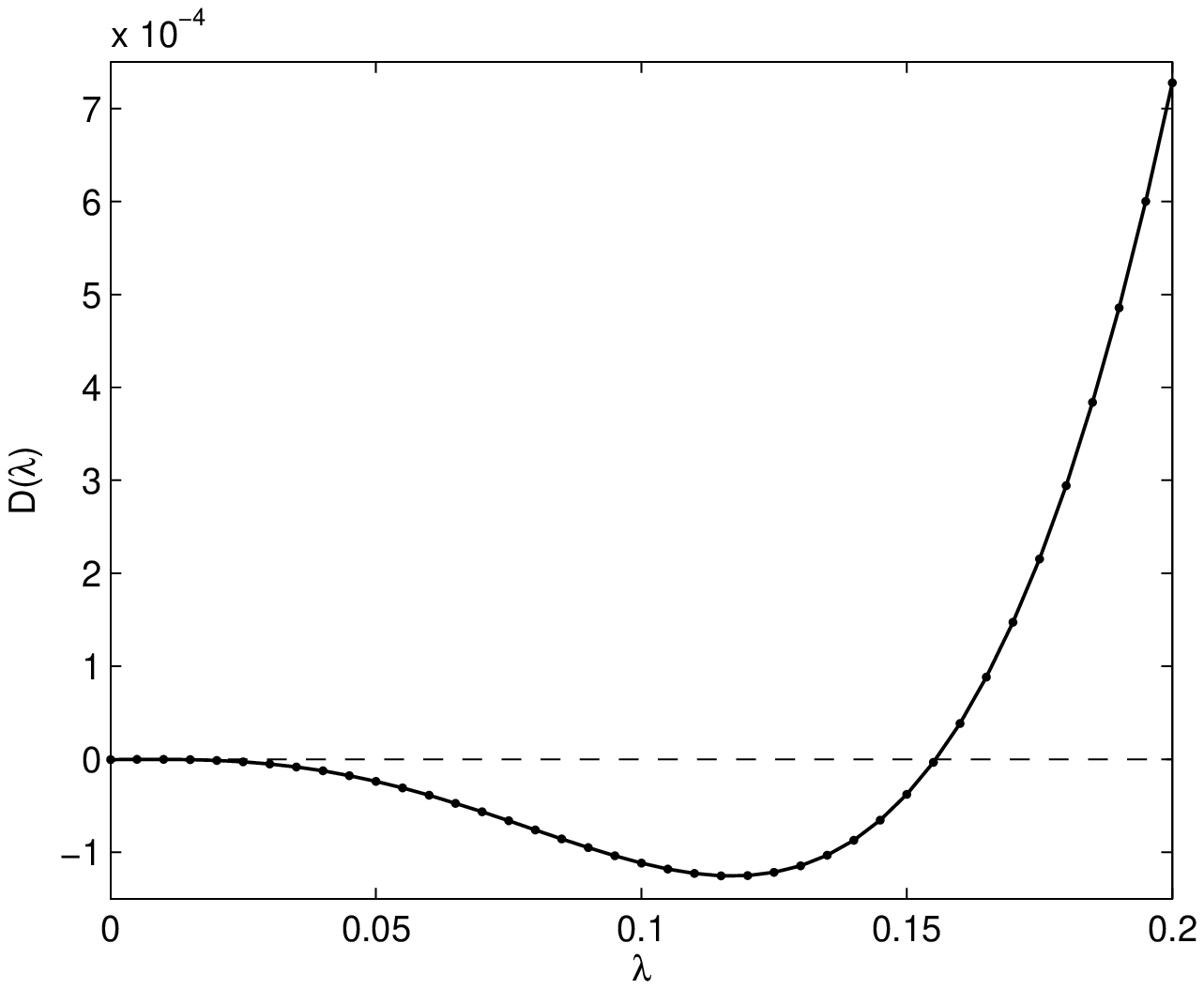}
\end{center}
\caption{The Evans function of the Boussinesq system for the 
unstable pulse having wave speed $c=0.4$. The left plot is generated 
by the Riccati-RK method with fixed coordinate patches identified by 
$\ib^-=\{1,2\}$ over $[-8,0]$ and $\ib^+=\{3,4\}$ over $[0,8]$. 
The right plot shows the result of the CO-RK method.}
\label{fig1} 
\end{figure}
 
\begin{table}
\caption{Zero of the Evans function for the Boussinesq problem 
computed with the Riccati-RK method (with fixed coordinate patches 
identified by $\ib=^-\{1,2\}$ over $[-8,0]$ and $\ib^+=\{3,4\}$ 
over $[0,8]$) and the CO-RK method. 
Here $N$ is the number of (equidistant) steps used in the mesh. }
\centering
		\begin{tabular}{ccc}
		\hline
		$N$&Riccati-RK&CO-RK\\
		\hline
		128&0.15544090&0.15540090\\
		256&0.15543184&0.15542952\\
		512&0.15543143&0.15543129\\
		1024&0.15543141&0.15543140\\
		2048&0.15543141&0.15543141\\	
		4096&0.15543141&0.15543141\\\hline
		\end{tabular}
	\label{tab:convergenceToEigenvalue}
\end{table}

\begin{figure}
\includegraphics[width=9cm,height=5cm]{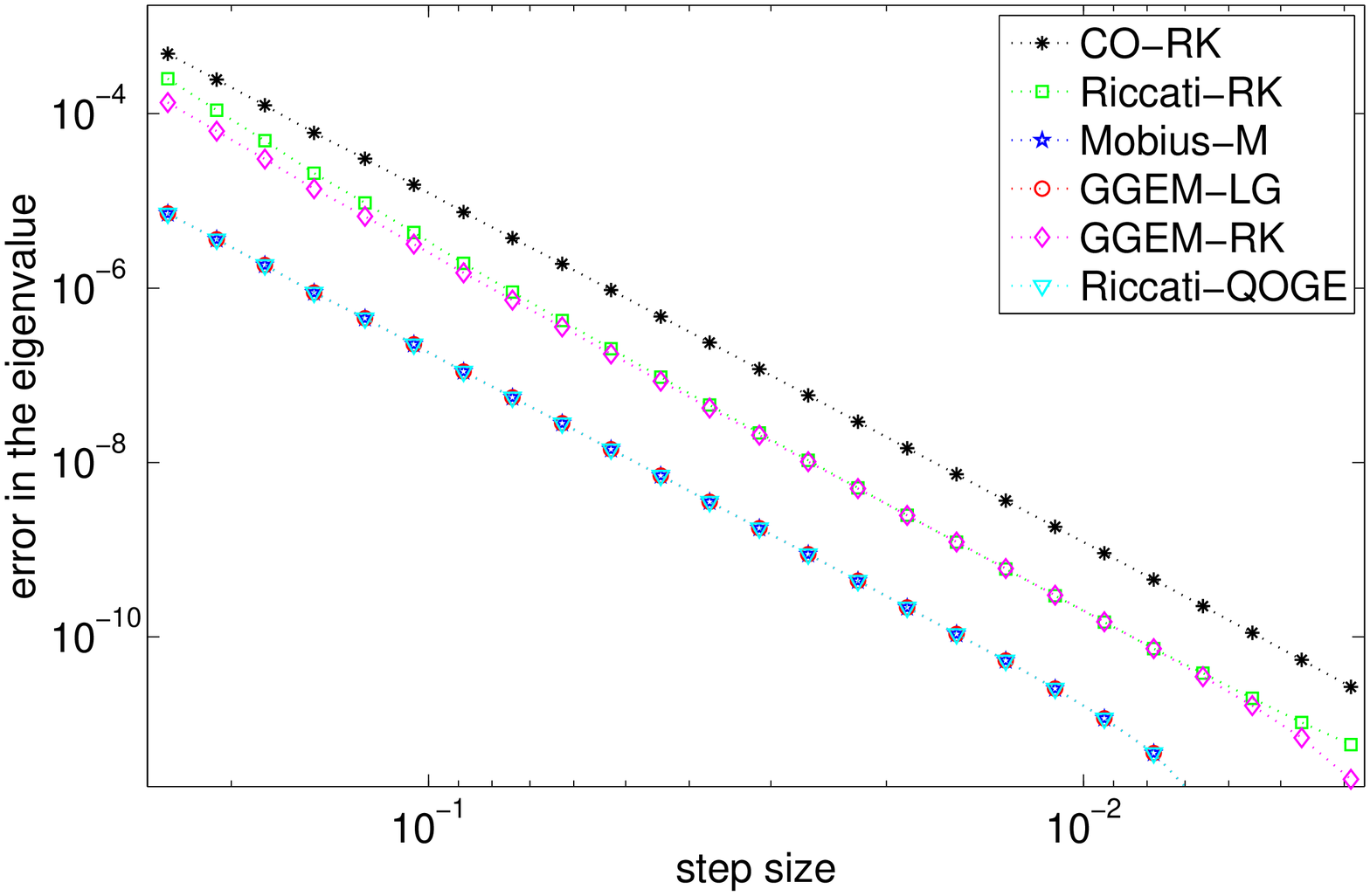}
\includegraphics[width=9cm,height=5cm]{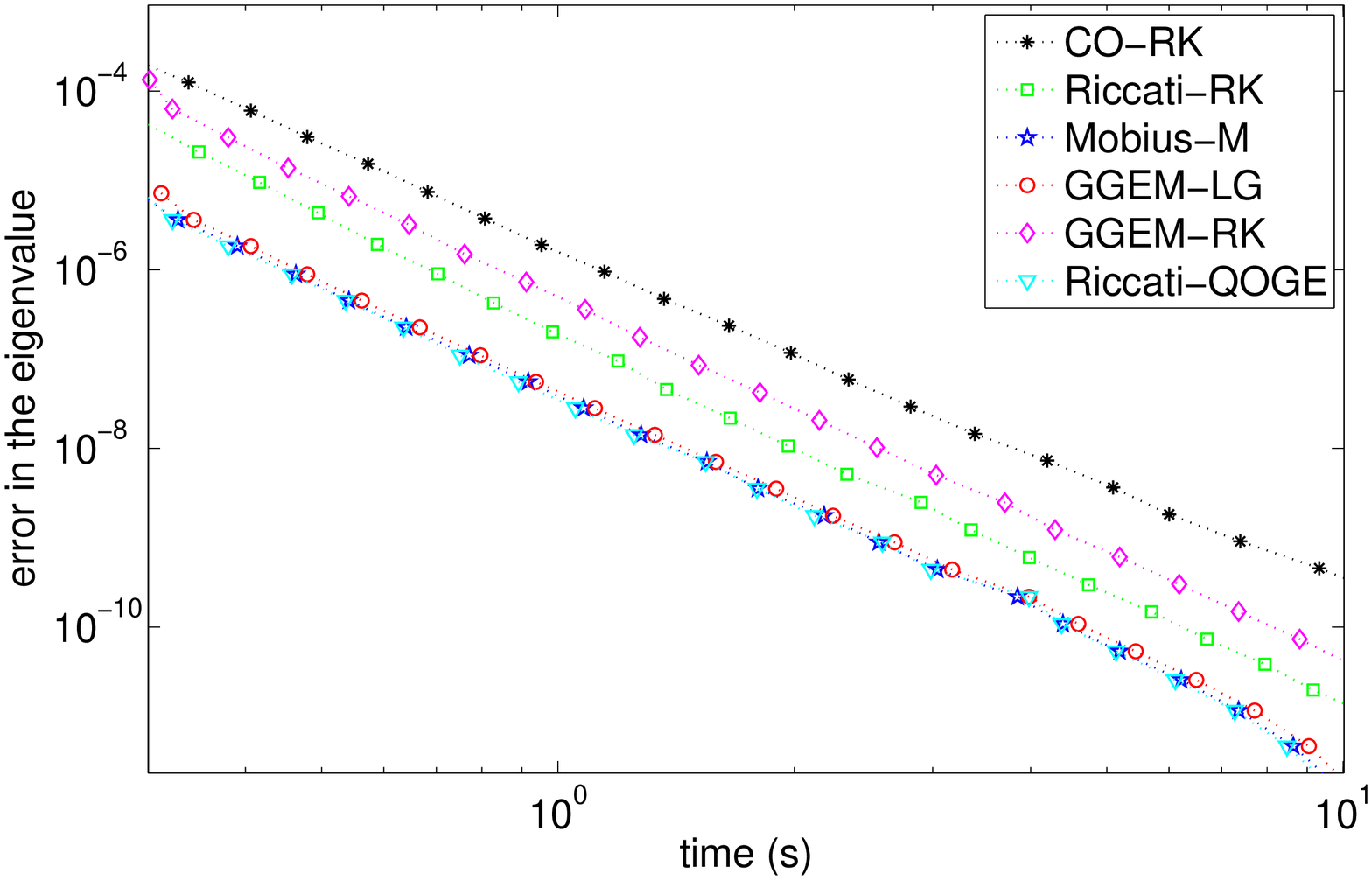}
\caption{Error in the eigenvalue, vs stepsize (upper panel) 
and vs cputime (lower panel), matching at $x_\ast=0$.}
\label{Bouss:matchzero}
\end{figure}

\begin{figure}
\includegraphics[width=9cm,height=5cm]{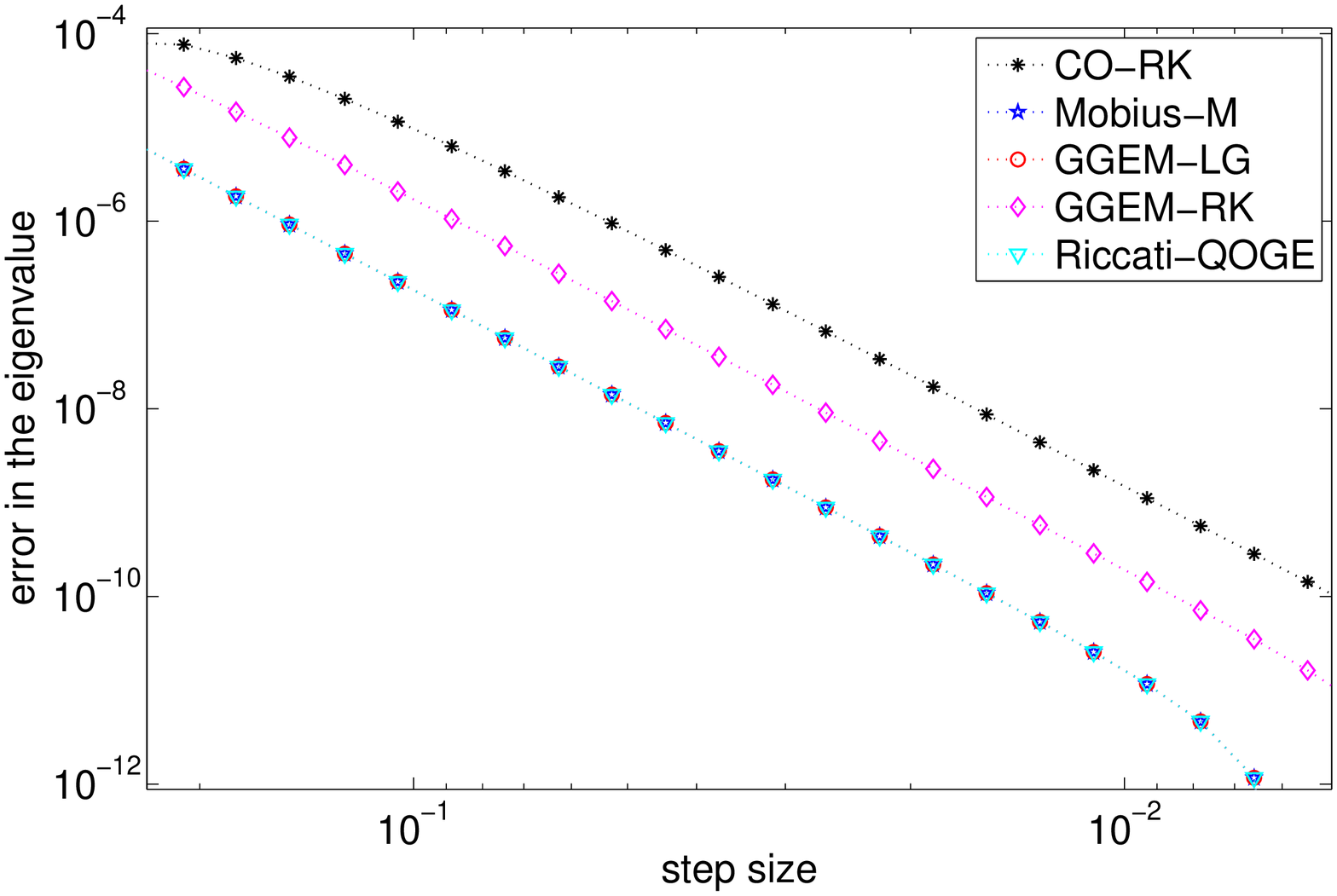}
\includegraphics[width=9cm,height=5cm]{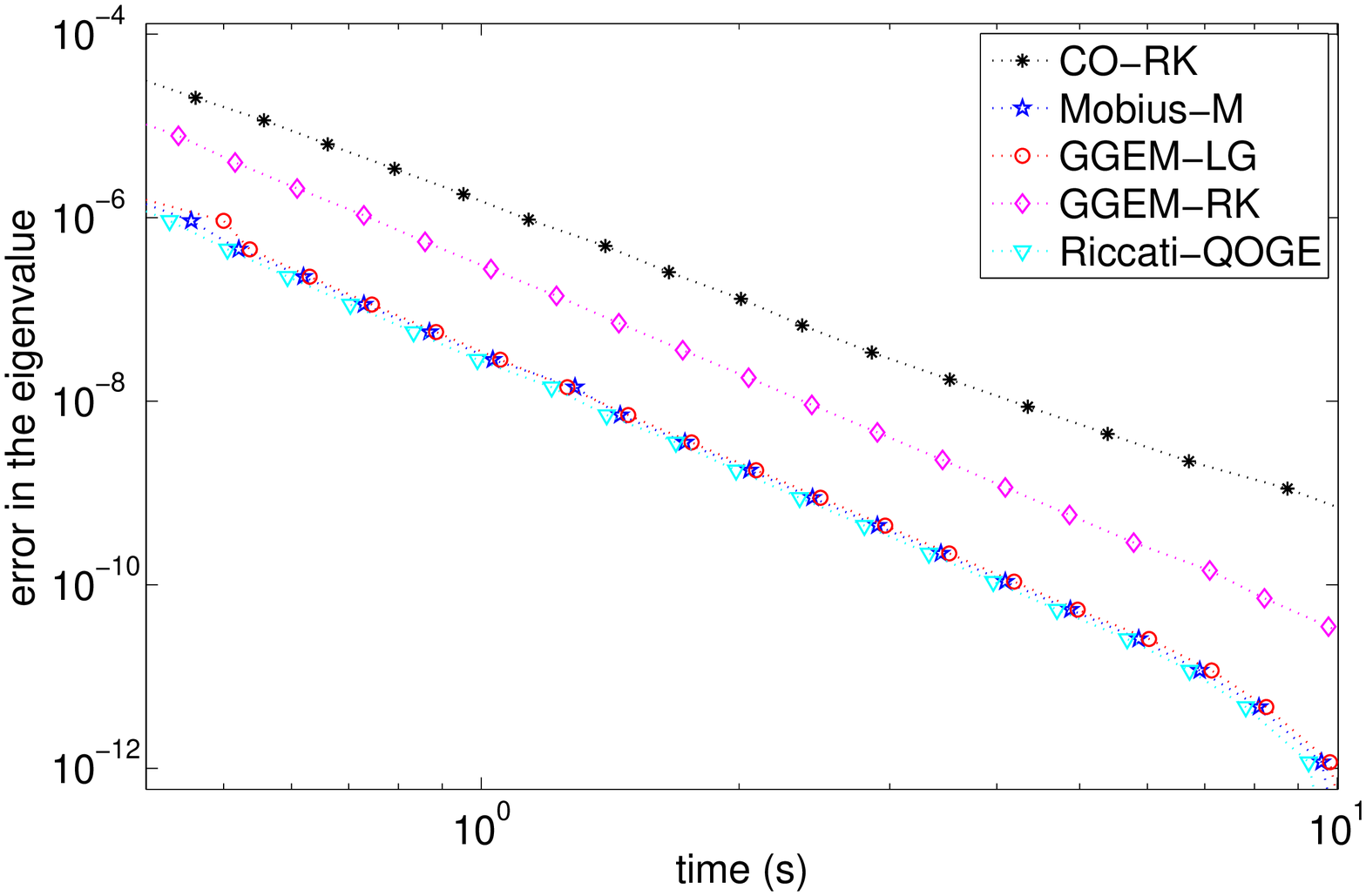}
\caption{Error in the eigenvalue vs stepsize (upper panel)  
and vs cputime (lower panel), matching at $x_\ast=8$.}
\label{Bouss:matchright}
\end{figure}

In Figure~\ref{fig1} we show the Evans function computed     
along the real axis from $\lambda=0$ to $\lambda=0.2$
for the unstable pulse with $c=0.4$. The Riccati-RK (left plot)
and CO-RK (right plot) methods
detect a zero of the Evans function near $\lambda=0.155$, 
indicating an unstable eigenvalue there. 
An accurate value of the eigenvalue can be found by using a standard 
root-finding method. This yields the values 
in Table~\ref{tab:convergenceToEigenvalue}. The Riccati-RK and 
the CO-RK methods both converge to the same 
eigenvalue when the number of steps $N$ increases. As a check, 
the Matlab \texttt{ode45} solver was used with a relative tolerance $10^{-8}$ and
absolute tolerance $10^{-10}$ to integrate the systems, leading to the 
same resulting eigenvalue: $\lambda=0.15543141$ for both methods.

Function evaluation for the Riccati vector field requires 
three matrix-matrix multiplications. 
This is the same number of matrix-matrix multiplications 
needed to evaluate the Drury--Oja vector field. 
However, the matrices in the Drury--Oja vector field 
are $n \times k$ and $n\times(n-k)$, respectively, 
while the matrices in the Riccati vector fields have 
smaller dimension: $(n-k)\times k$. Because of the 
smaller dimension of the systems to be integrated, 
our Riccati approach is faster than the continuous orthogonalization problem.
For example, to construct Figure~\ref{fig1} the 
Evans function was evaluated at 200 distinct $\lambda$ values 
between $\lambda=0$ and $\lambda=0.2$. Using the fourth-order 
Runge-Kutta method with $N=512$ steps, this required 33 seconds 
for the CO-RK method, while the Riccati-RK method 
needed 24 seconds (Matlab-implementation, CPU 2.4GHz).

\begin{figure}
		\includegraphics[width=0.45\textwidth]{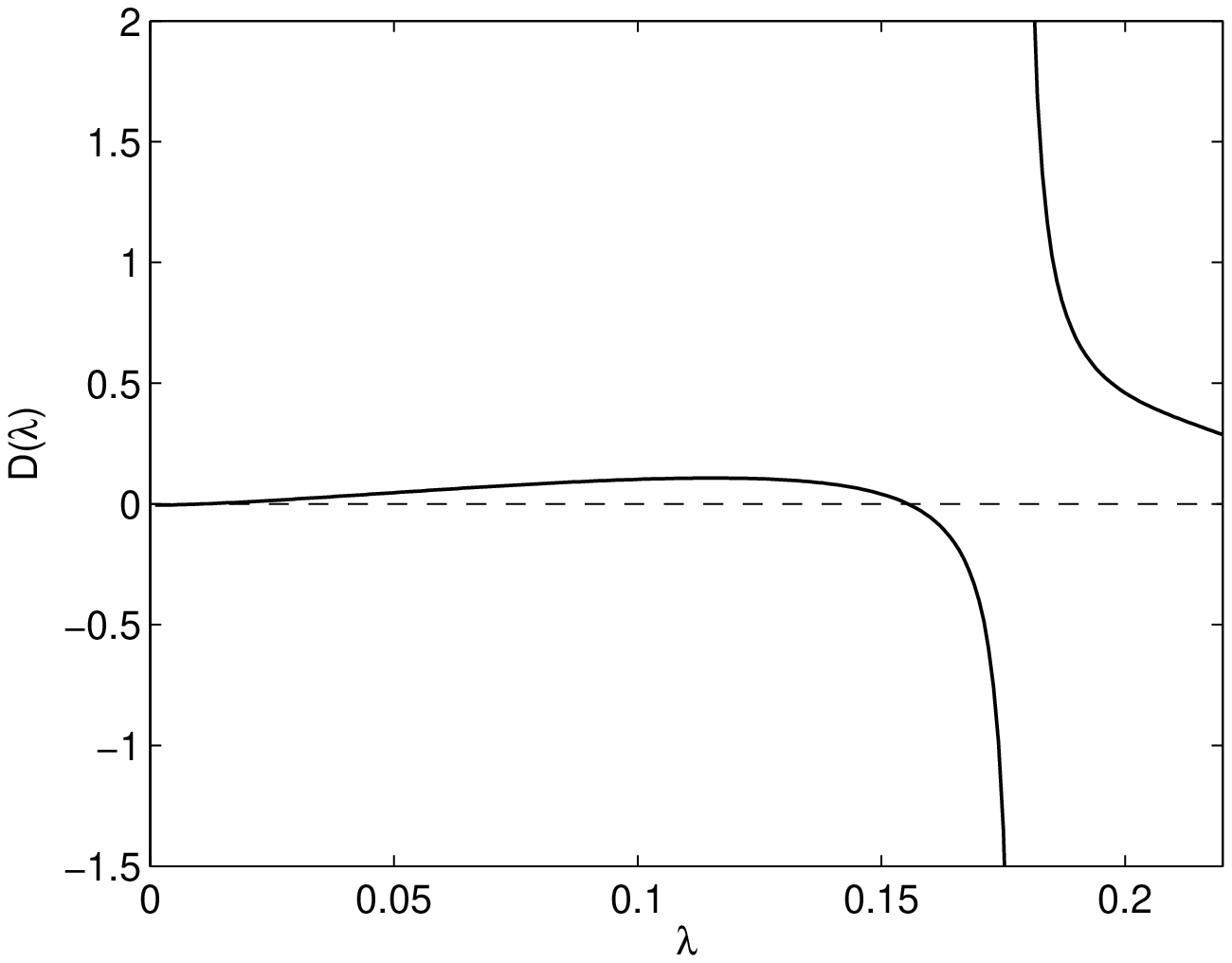}
		\includegraphics[width=0.45\textwidth]{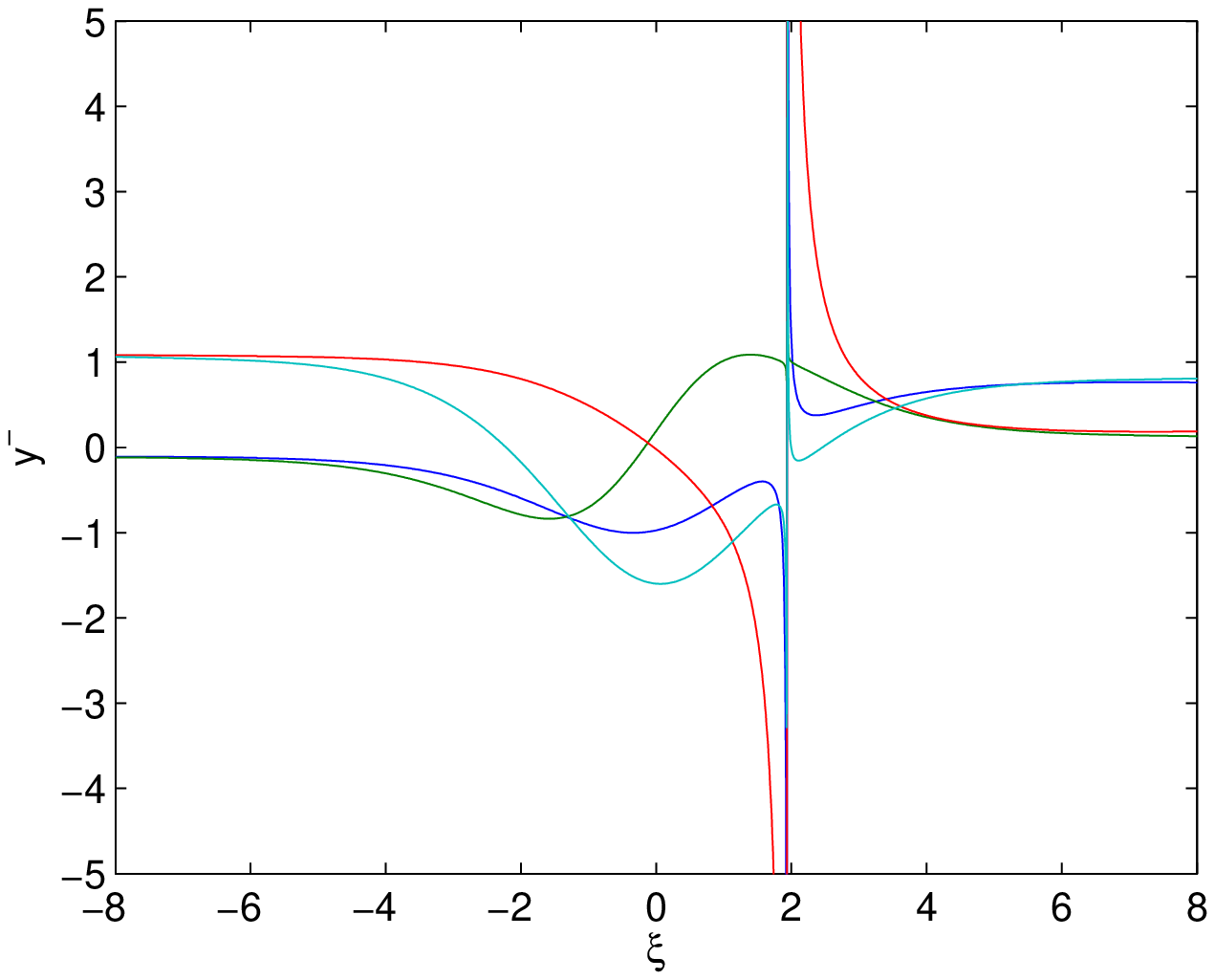}
\caption{The Evans function when the M\"obius--Magnus method is
applied (left panel) with the matching point as $x_\ast=+8$.
The entries of the solution $\hat y^-$,
passing through the singularity when integrating from 
$-8$ to $8$, for $\lambda=0.15543141$ (right panel).}\label{fig3bb} 
\end{figure}

\begin{figure}
		\includegraphics[width=0.45\textwidth]{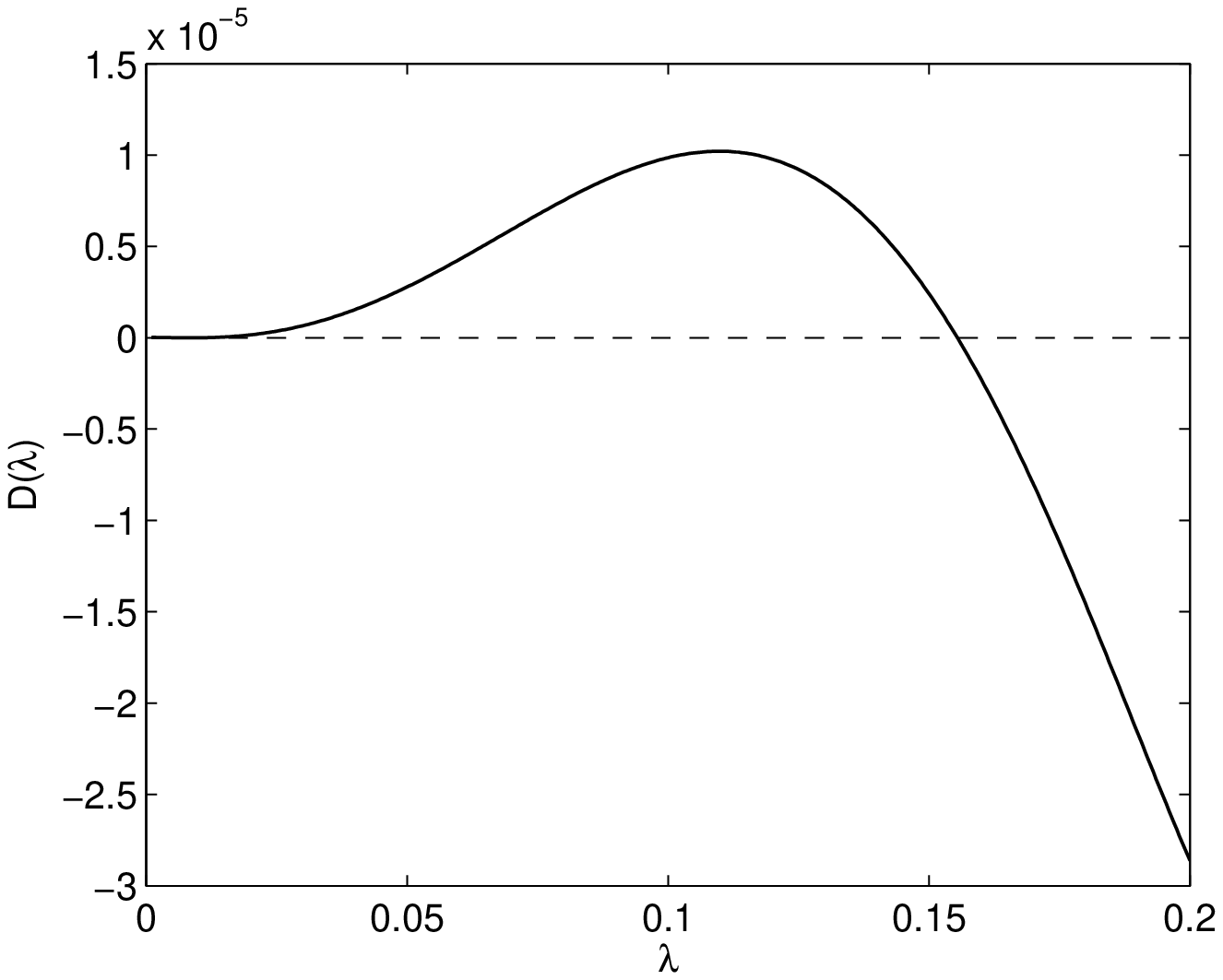}
		\includegraphics[width=0.45\textwidth]{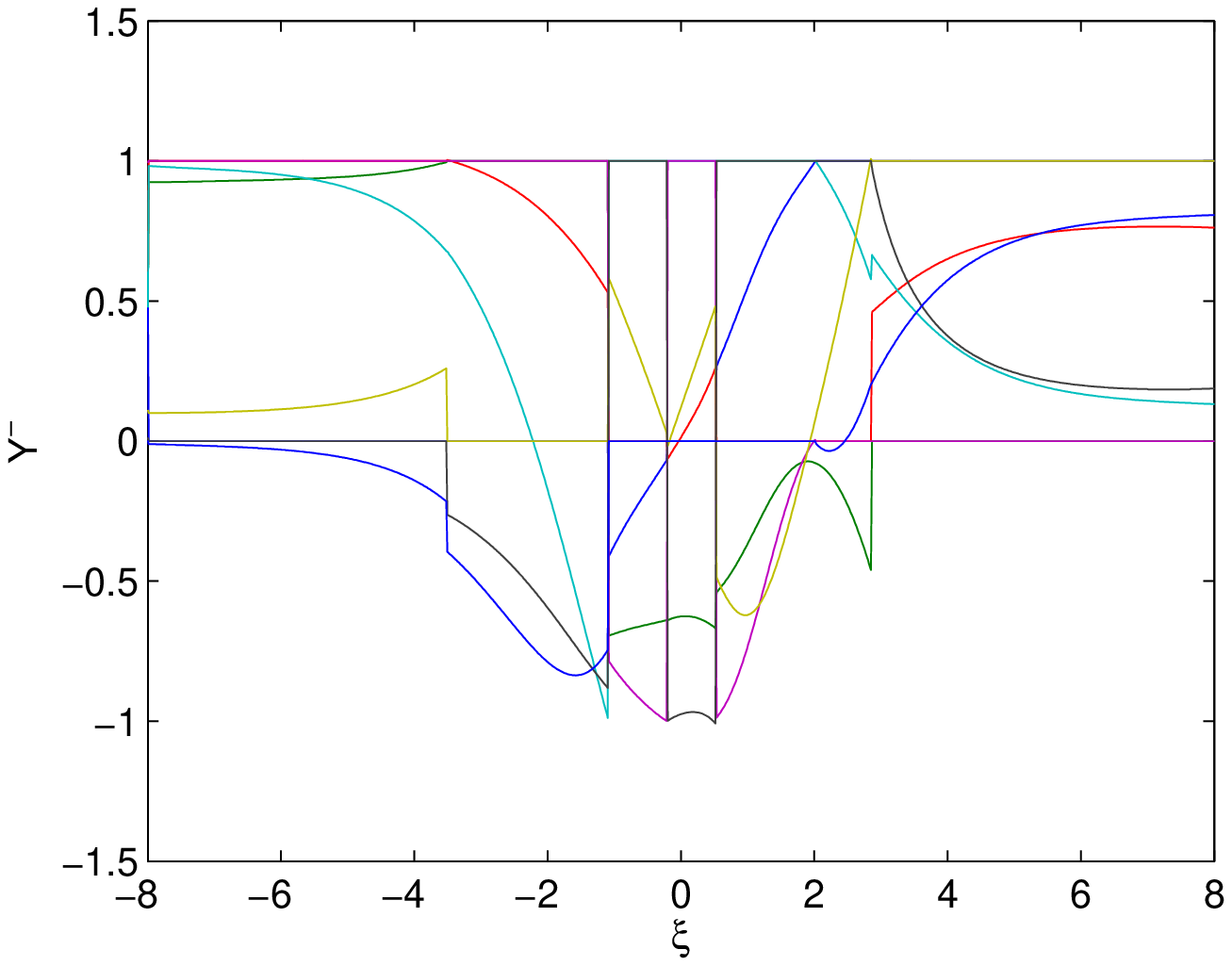}
\caption{The Evans function which results when the 
GGEM-RK scheme is applied over $[-8,8]$, matching point in 
$x_\ast=8$ (left panel).
The entries of $y_{\ib^\circ}$ for $\lambda=0.15543141$ (right panel).}\label{figRGE} 
\end{figure}

\begin{figure}
\includegraphics[width=0.45\textwidth]{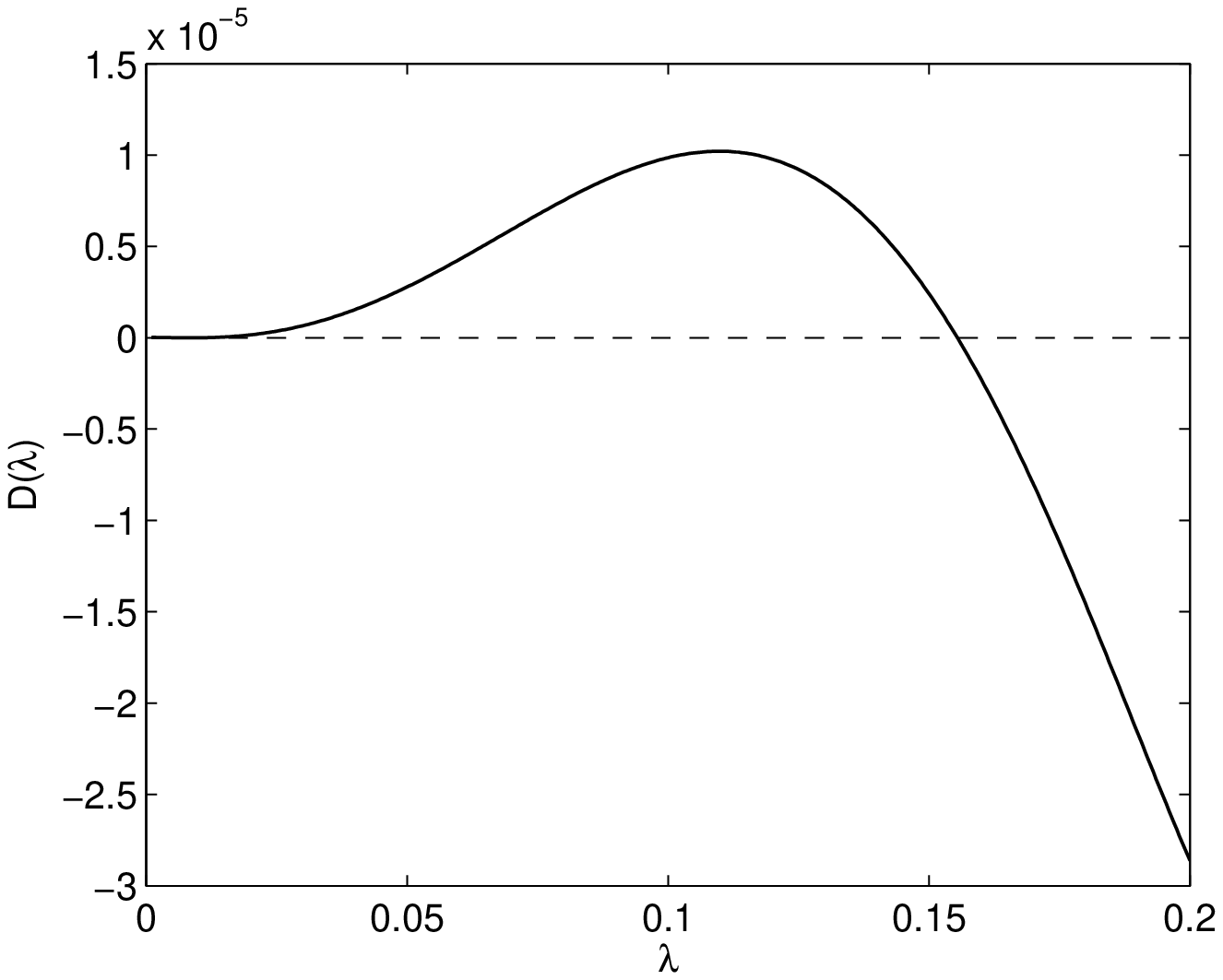}
\includegraphics[width=0.45\textwidth]{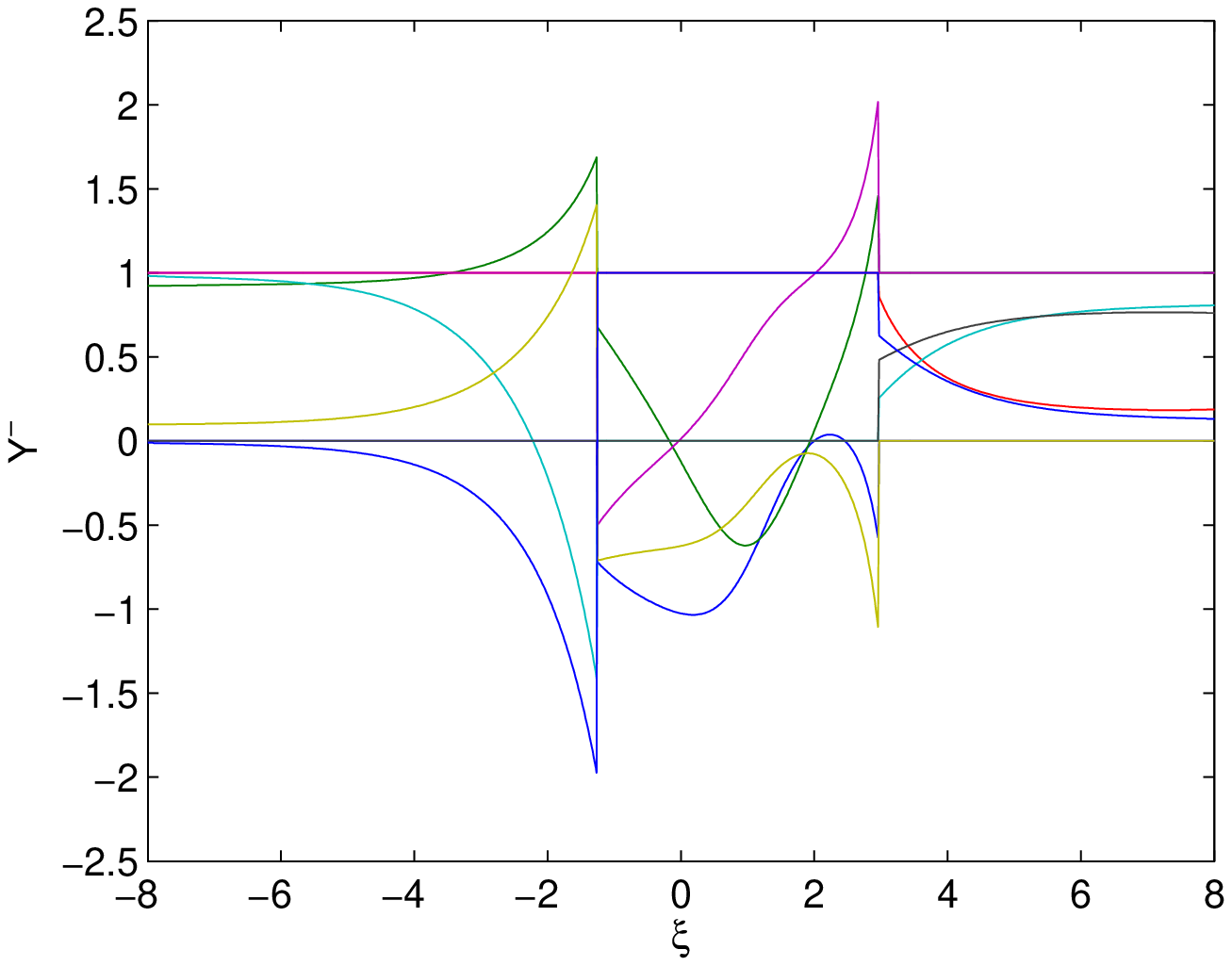}
\caption{The Evans function when the Riccati-QOGE method is 
applied (left panel). The entries of $y_{\ib^\circ}$ are shown for 
$\lambda=0.15543141$ (right panel).
The criterion used for swapping to a new coordinate patch
was $\|\hat y^-\|_{\infty}>2$.
}\label{figRGE2} 
\end{figure}

\begin{figure}
\includegraphics[width=9cm,height=5cm]{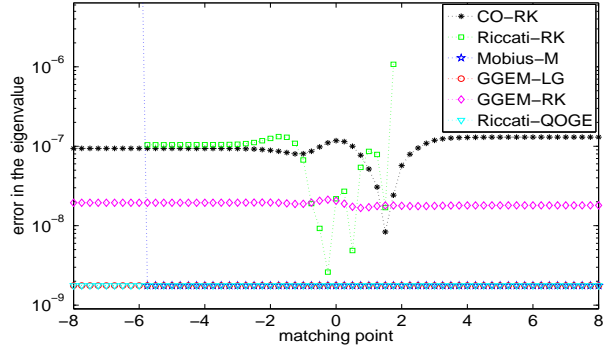}
\caption{Error in the eigenvalue for different choices 
of the matching point. The number of steps in the 
equidistant mesh was $N=512$.}
\label{Bouss:varymatcheval}
\end{figure}

\begin{figure}
\includegraphics[width=9cm,height=5cm]{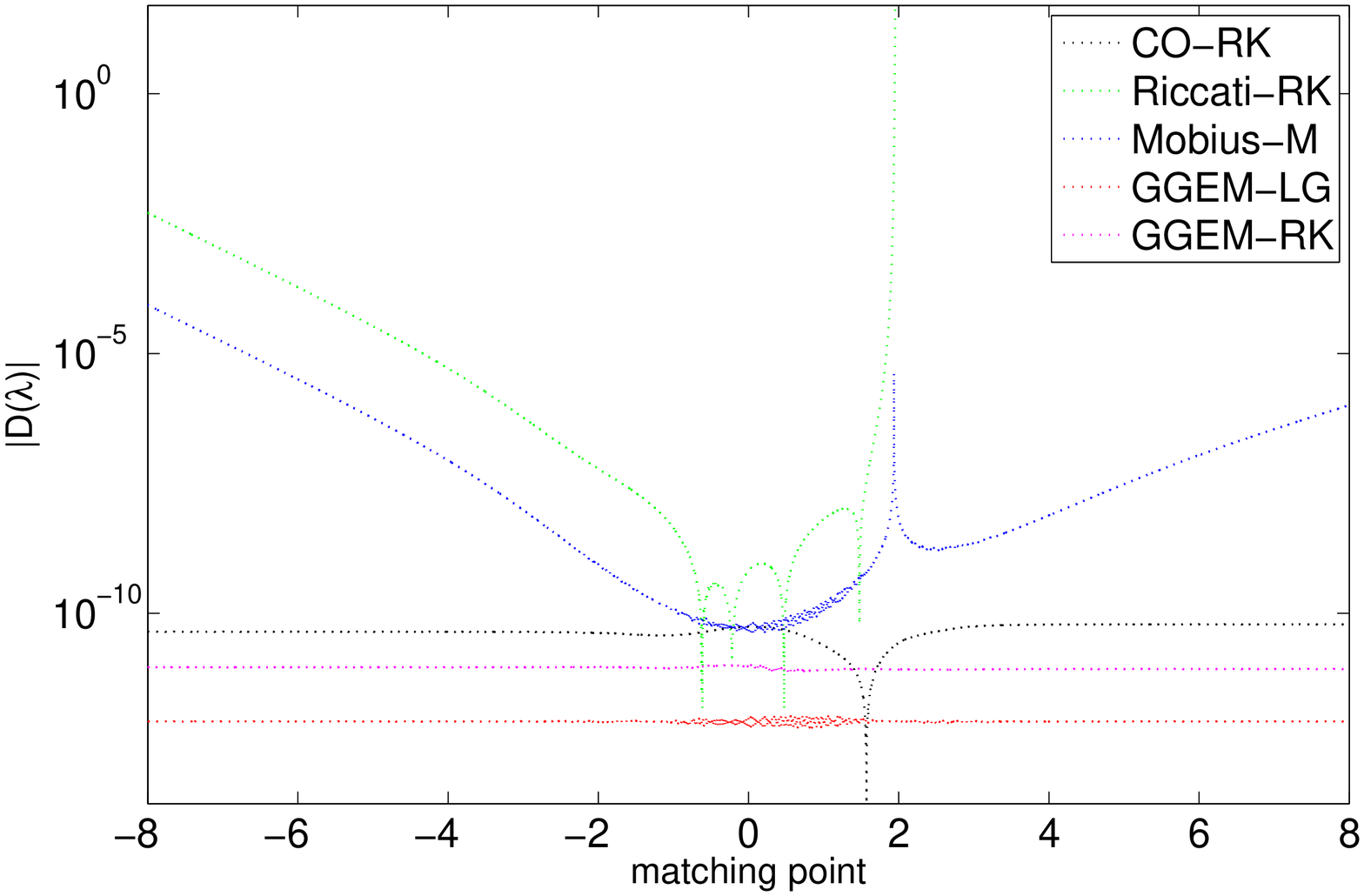}
\caption{$|D(\lambda)|$ for $\lambda$ equal to the eigenvalue 
for different matching points.
The number of steps in the equidistant mesh was $N=512$.}
\label{Bouss:varymatchEvans}
\end{figure}

In Figure~\ref{Bouss:matchzero} we compare the 
error in the eigenvalue and efficiency of computation 
for all six methods, when we match at $x_\ast=0$.
We see that the methods that use the Magnus expansion to
advance the solution on the Stiefel manifold are the most
accurate for a given stepsize. They are also the most efficient,
delivering the best accuracy for given computational effort.
The Riccati-RK method does not suffer from singularities
for the chosen fixed patches when matching at $x_\ast=0$,
at least for the range of values of the spectral parameter
in the vicinity of the eigenvalue (as well as the origin 
and anywhere inbetween). However, if we change the matching
point to $x_\ast=\ell_+=+8$ there are singularities in the 
Riccati-RK solution (as a result of a vanishing determinant 
of $u^-$). In particular, a singularity appears 
around $x=2$ for $\lambda$ equal to the eigenvalue (see Figure~\ref{fig3bb}).
Hence we compare the remaining five methods in 
Figure~\ref{Bouss:matchright} in this case. We see that 
using the M\"obius--Magnus method to integrate through a 
singularity does not introduce loss in accuracy. The resulting 
Evans function can have poles, as seen in Figure \ref{fig3bb},
which appear at $\lambda$-values where the Riccati equation 
has a singularity at the matching point. This means that in some cases 
the matching point should be chosen rather carefully 
in order not to have the poles interfering with the eigenvalue(s).  
When applying GGEM-RK the Evans function is analytic and the choice of the 
matching point is less important. 

Figure~\ref{figRGE} shows the Evans function obtained when the 
GGEM-RK evolves $y_{\ib}$ from $\ell_-=-8$ to $\ell_+=+8$. 
To construct the plot in Figure~\ref{figRGE}, 
the quasi-optimal Gaussian elimination process was applied 
at each step in the integration. However it is clear 
from the right plot in Figure~\ref{figRGE}, that multiple 
successive steps can be integrated in the same 
coordinate patch. For example, between $x=-8$ and $x=-3$ the 
coordinate patch does not change. Performing the whole 
quasi-optimal Gaussian elimination process only when a certain 
criterion is satisfied, reduces the computing time. 
Using the Riccati-QOGE method, we change the 
coordinatization when $\|\hat y^-\|_{\infty}>2$. This 
generates an Evans function very similar to that in
Figure~\ref{figRGE}. As seen in Figure~\ref{figRGE2} 
the quasi-optimal Gaussian elimination process is then performed 
only two times for $\lambda$ equal to the eigenvalue. 

We compare the error in computing the eigenvalue for different
choices of matching point---in fact for $x_\ast$ anywhere in the interval
$[-8,8]$---for all six methods in Figure~\ref{Bouss:varymatcheval}.
We see that the most accurate and robust methods are the GGEM-LG
and Riccati-QOGE methods. Some methods, such as the Riccati-RK method 
as discussed already, break down when singularities impinge on
the matching point---the singularities in the Evans function are observed 
in Figure~\ref{Bouss:varymatchEvans}. Generally we also see 
in Figure~\ref{Bouss:varymatcheval} that the GGEM-LG
and Riccati-QOGE methods outperform the CO-RK method 
in terms of accuracy.

Overall, we observe in this example that when computing the eigenvalue,
those methods based on the Magnus expansion are superior in accuracy
and efficiency. 
Note that for GGEM-RK, the quasi-optimal Gaussian elimination process 
is an additional $nk^2$ operation. However, to ensure analyticity 
for the CO-RK method, there are two additional 
matrix-matrix multiplications in the equations for $\det R^\pm$ 
(operational cost $kn^2$) required at each step.

\subsection{Autocatalytic fronts}
As a second example, we study travelling waves in a model of 
autocatalysis in an infinitely extended medium 
\begin{align*}
	u_t &=\delta u_{xx}+cu_x-uv^m,\\
	v_t &=v_{xx}+cv_x+uv^m.
\end{align*}
Here $u(x,t)$ is the concentration of the reactant and $v(x,t)$ 
is the concentration of the autocatalyst. We suppose $(u,v)$ 
approaches the stable homogeneous steady state $(0,1)$ as 
$x \to -\infty$, and the unstable  homogeneous steady 
state $(1,0)$ as $x \to +\infty$. The diffusion parameter 
$\delta$ is the ratio of the diffusivity of the reactant 
to that of the autocatalyst and $m$ is the order of the 
autocatalytic reaction. The speed of the co-moving reference
frame is $c$. The system is globally well-posed for 
smooth initial data and any finite $\delta>0$ and $m \geq 1$.

\begin{figure}
	\includegraphics[width=9cm,height=5cm]{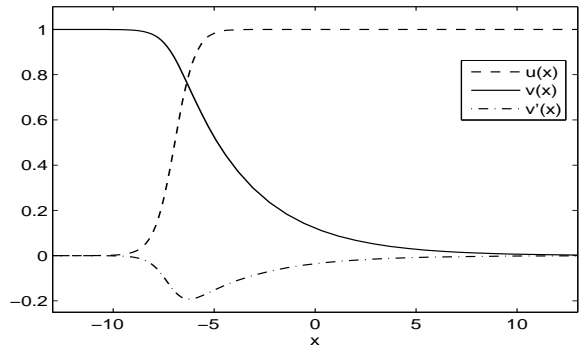}
\caption{The travelling wave solution for $\delta=0.1$ and $m=9$.}\label{figtw}
\end{figure}

From Billingham and Needham~\cite{BN} we know that a 
unique heteroclinic connection between the unstable and 
stable homogeneous steady states exists for wavespeeds $c\geq c_{\min}$. 
The unique travelling wave for $c=c_{\min}$ converges exponentially 
to the homogeneous steady states and is computed by a 
simple shooting algorithm (see Balmforth, Craster and Malham~\cite{BCM}). 
The resulting travelling wave for $\delta=0.1$ and $m=9$ 
is shown in Figure~\ref{figtw}.

\begin{figure}
	\includegraphics[width=0.45\textwidth]{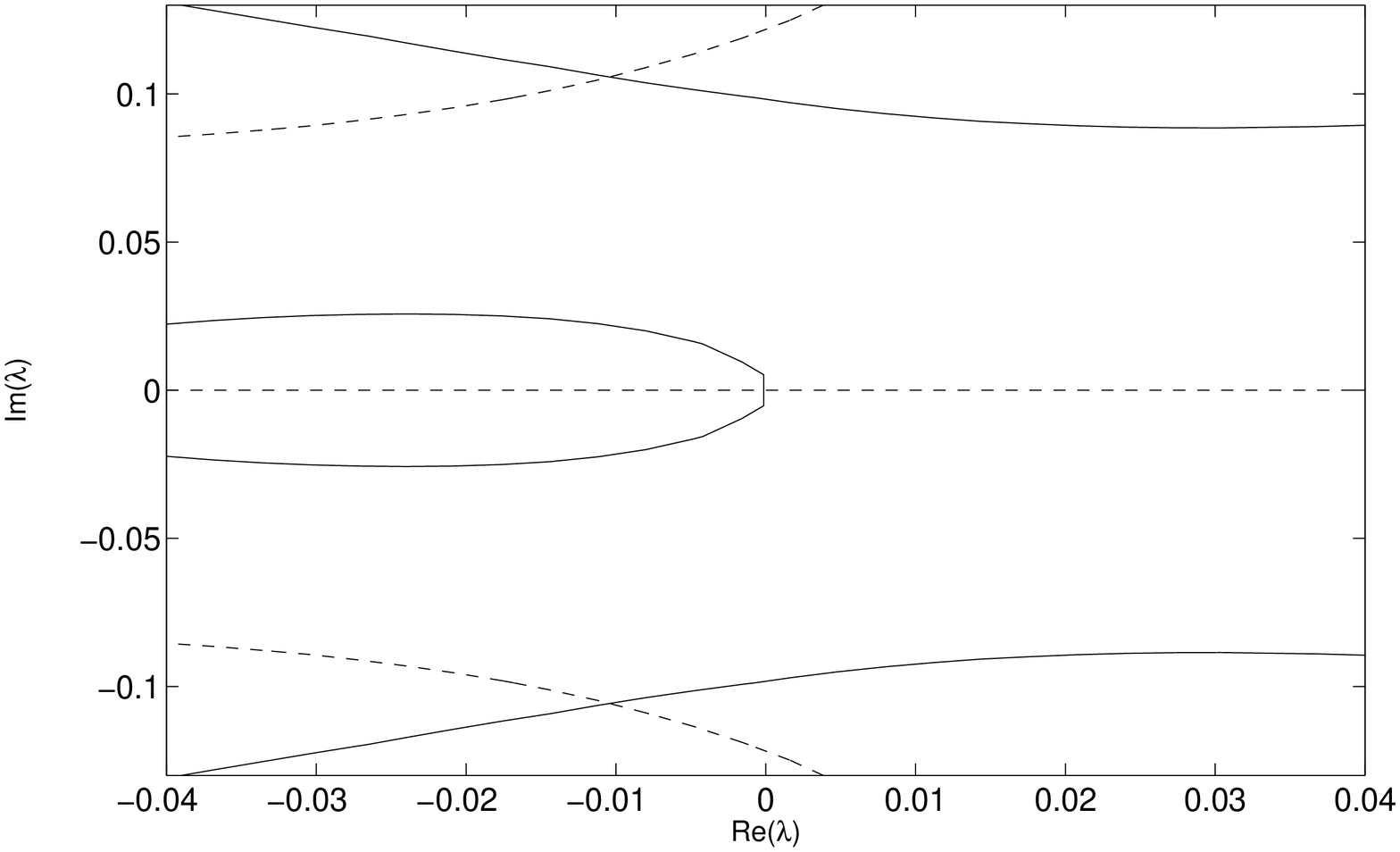}
        \includegraphics[width=0.45\textwidth]{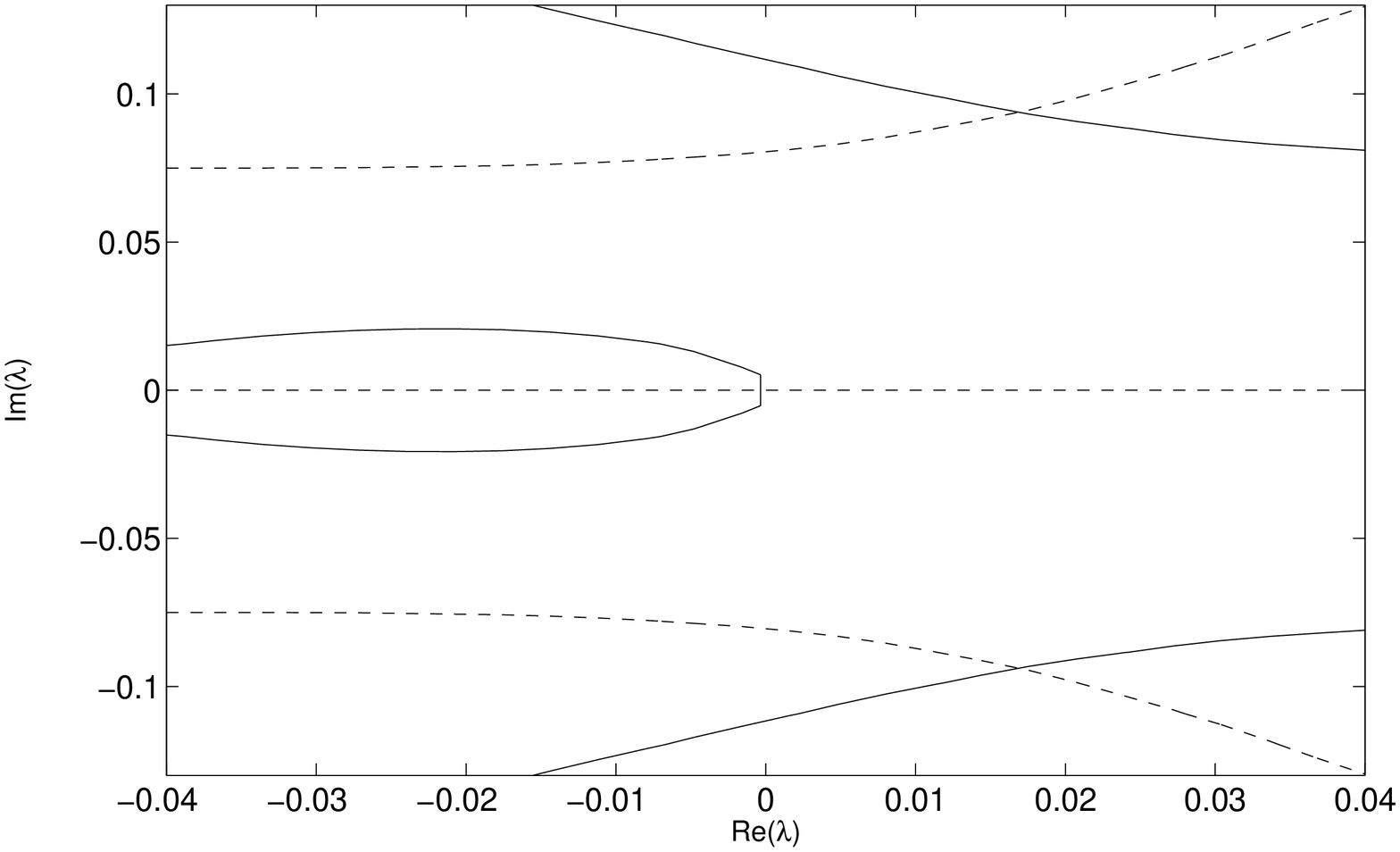}
\caption{Zero contour lines of the real (solid) and imaginary (dashed) 
parts of the Evans function for the autocatalysis problem
with $\delta=0.1$ and $m=8$ (left panel) and $m=9$ (right panel).}\label{fig4}
\end{figure}

\begin{figure}
\includegraphics[width=9cm,height=4cm]{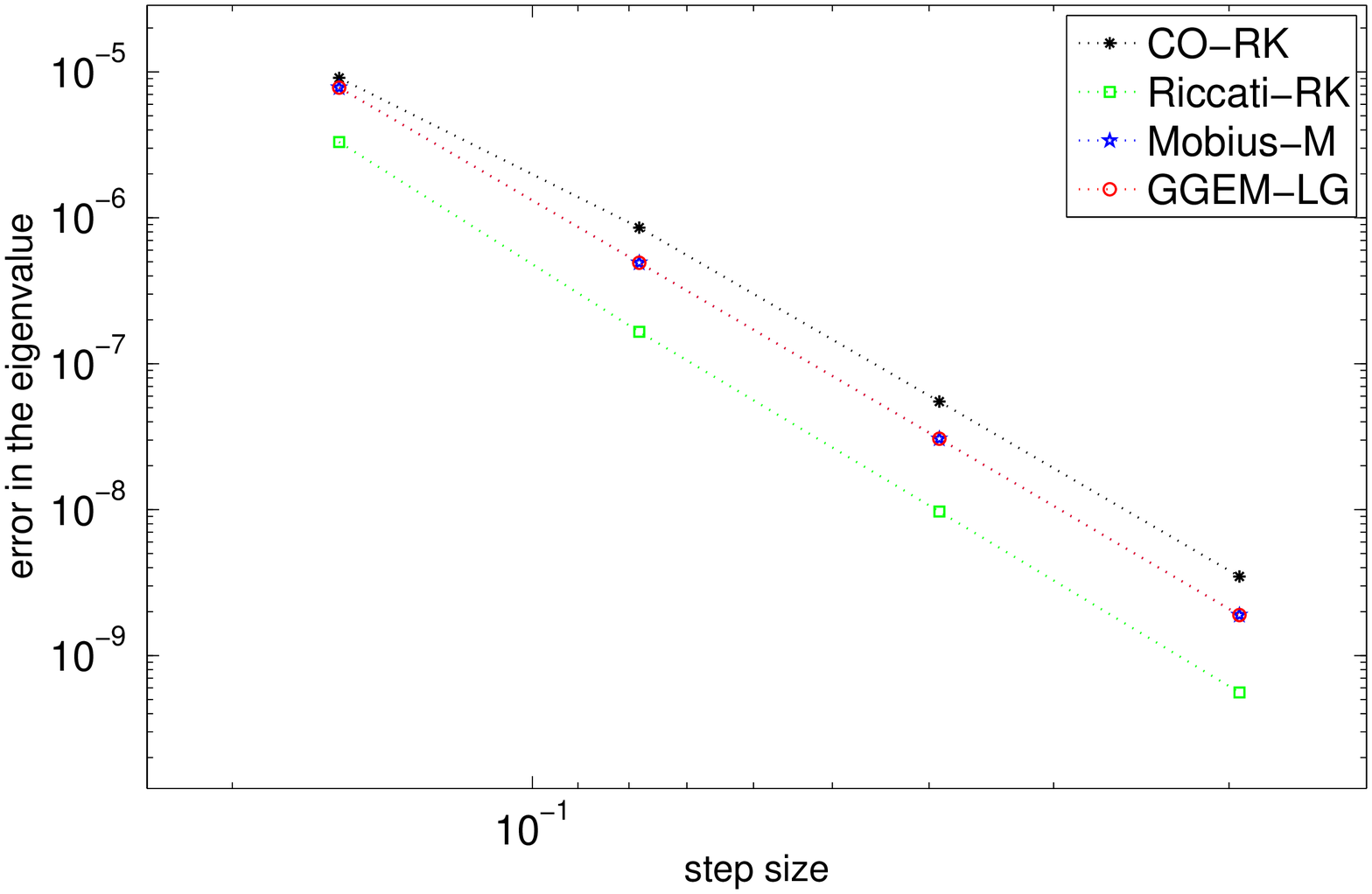}
\includegraphics[width=9cm,height=4cm]{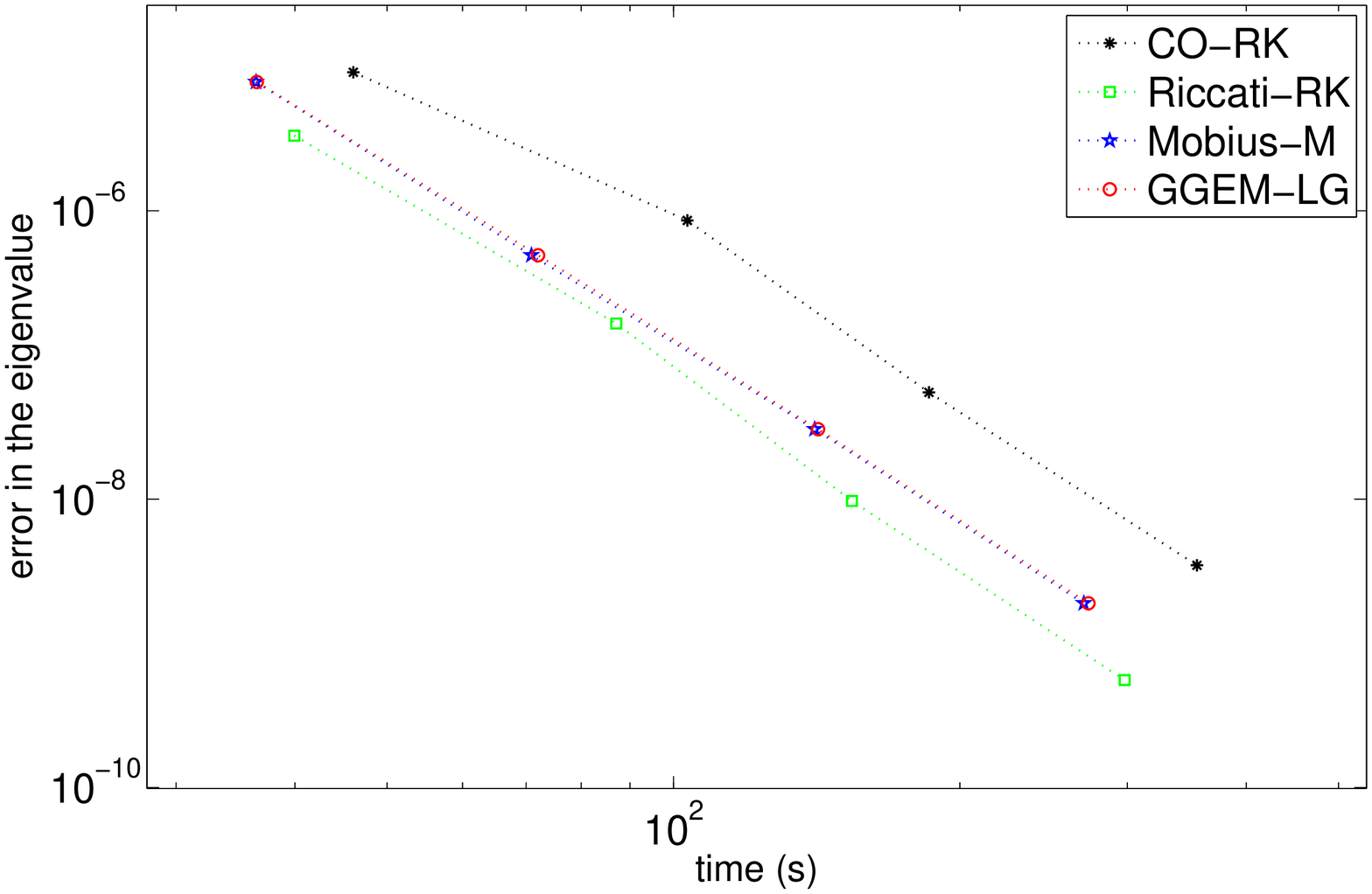}
\caption{Error in the eigenvalue in the first quadrant 
when $\delta=0.1$ and $m=9$ for different methods,
vs stepsize (upper panel) and vs cputime (lower panel), matching at $x_\ast=0$.}
\label{autocat:error}
\end{figure}

\begin{figure}
\includegraphics[width=9cm,height=5cm]{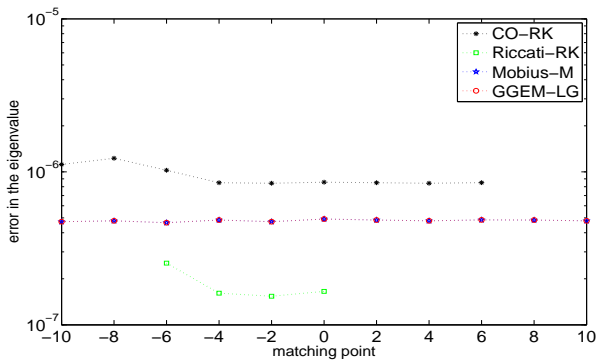}
\caption{Error in the eigenvalue in the first quadrant 
when $\delta=0.1$ and $m=9$, for different methods, 
for different matching points.
The number of steps in the equidistant mesh is $N=256$.}
\label{autocat:errormatch}
\end{figure}

The stability of the travelling wave of velocity $c$ 
can be deduced from the location of the spectrum of the 
eigenvalue problem $Y'=A(x;\lambda)Y$, where
\begin{equation*}
A(x;\lambda)=\left(\begin{matrix}0&0&1&0\\0&0&0&1\\
\lambda/\delta+\bar v^m/\delta&m\bar u\bar v^{m-1}/\delta&-c/\delta&0\\
-\bar v^m&\lambda-m\bar u\bar v^{m-1}&0&-c\end{matrix}\right),
\end{equation*}
where $\bar u$ and $\bar v$ represent the travelling wave solution.

The pulsating instability occurs when $\delta<1$ 
is sufficiently small and $m$ is sufficiently large 
(see Metcalf, Merkin and Scott~\cite{MMS} 
and Balmforth, Craster and Malham~\cite{BCM}).
For $\delta$ fixed and $m$ increasing, a complex conjugate pair of
eigenvalues crosses into the right-half $\lambda$-plane 
signifying the onset of instability
via a Hopf bifurcation. Figure~\ref{fig4} shows the onset 
of this instability for $\delta=0.1$
as $m$ is increased from $8$ to $9$
(see Aparicio, Malham and Oliver~\cite{Apar}).
The figure shows the zero contour lines of the real 
and imaginary parts of the Evans function. 
Solid lines correspond to 
zero contours of the real part of $D(\lambda)$, 
dashed lines to the imaginary part of $D(\lambda)$. 
We see that a complex-conjugate 
pair of eigenvalues crosses into the right-half plane,
indicating the onset of instability. 
Figure~\ref{fig4} was constructed using the Riccati-RK method with 
the fixed coordinate patches identified by $\ib^-=\{1,2\}$ 
from $x=-10$ to $x_\ast=-7$,  and $\ib^+=\{3,4\}$ from 
$x=10$ to $x_\ast=-7$. The matching point $x_\ast=-7$ is chosen roughly 
centred on the wavefront. 

\begin{figure}
\includegraphics[width=0.32\textwidth]{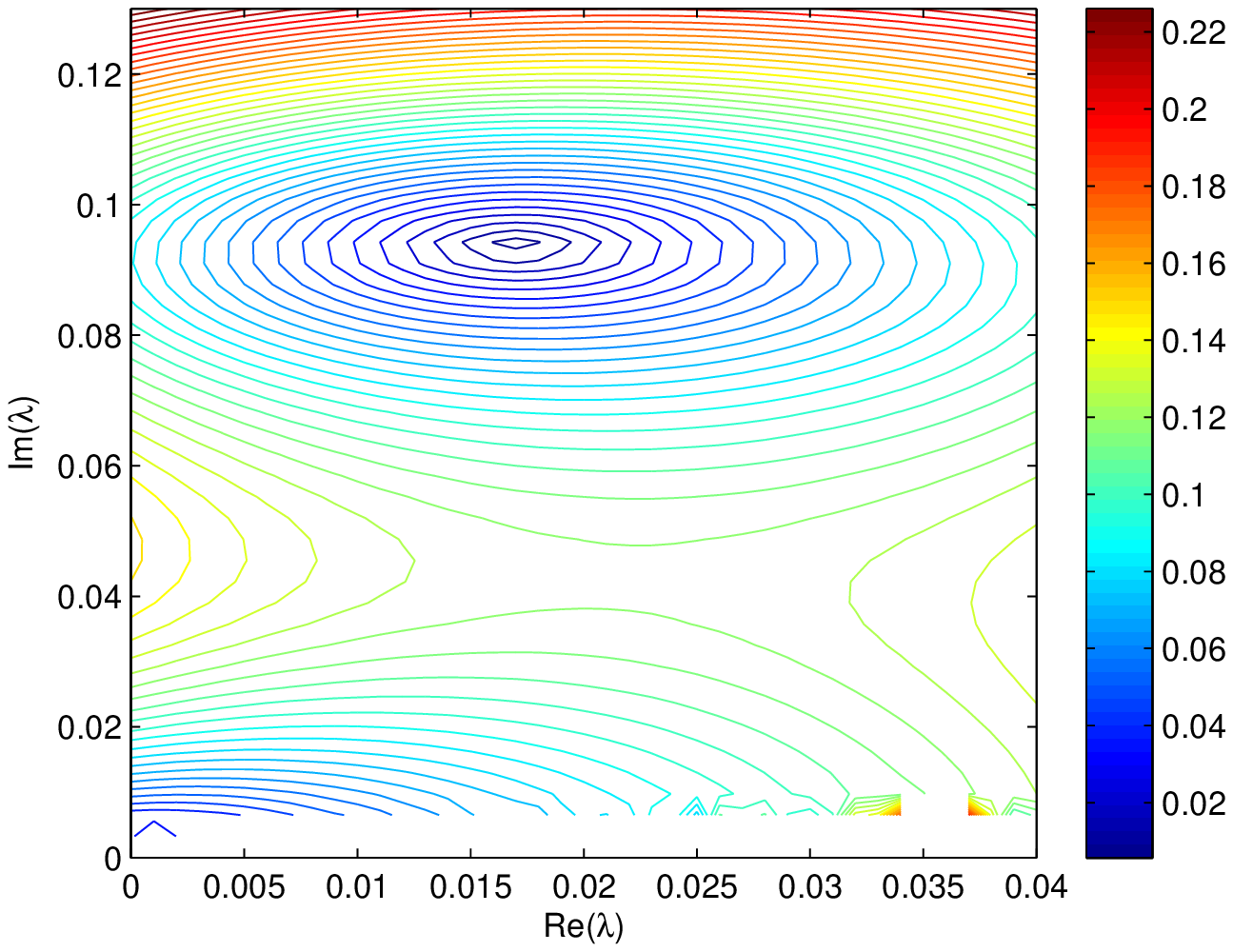}
\includegraphics[width=0.32\textwidth]{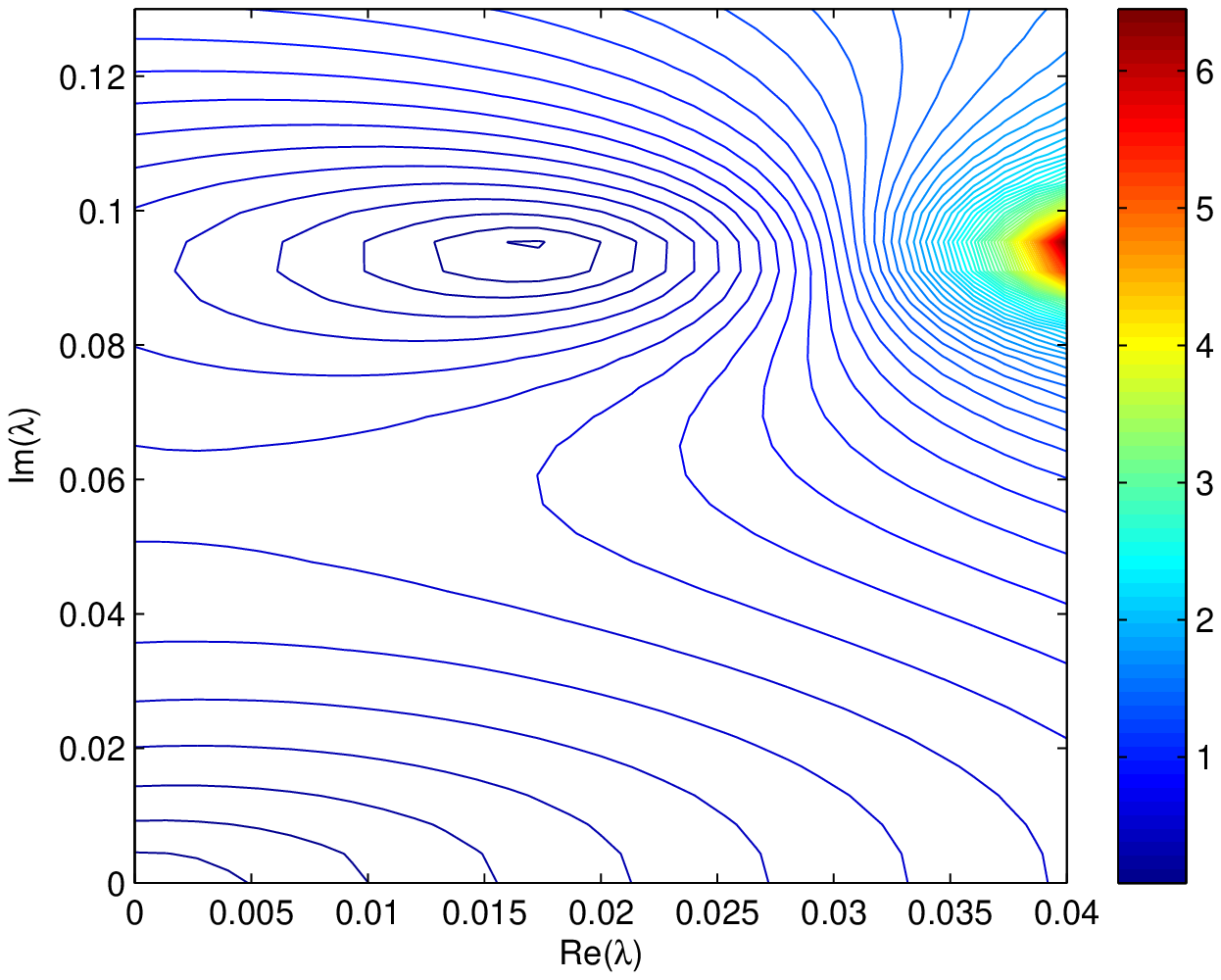}
\includegraphics[width=0.32\textwidth]{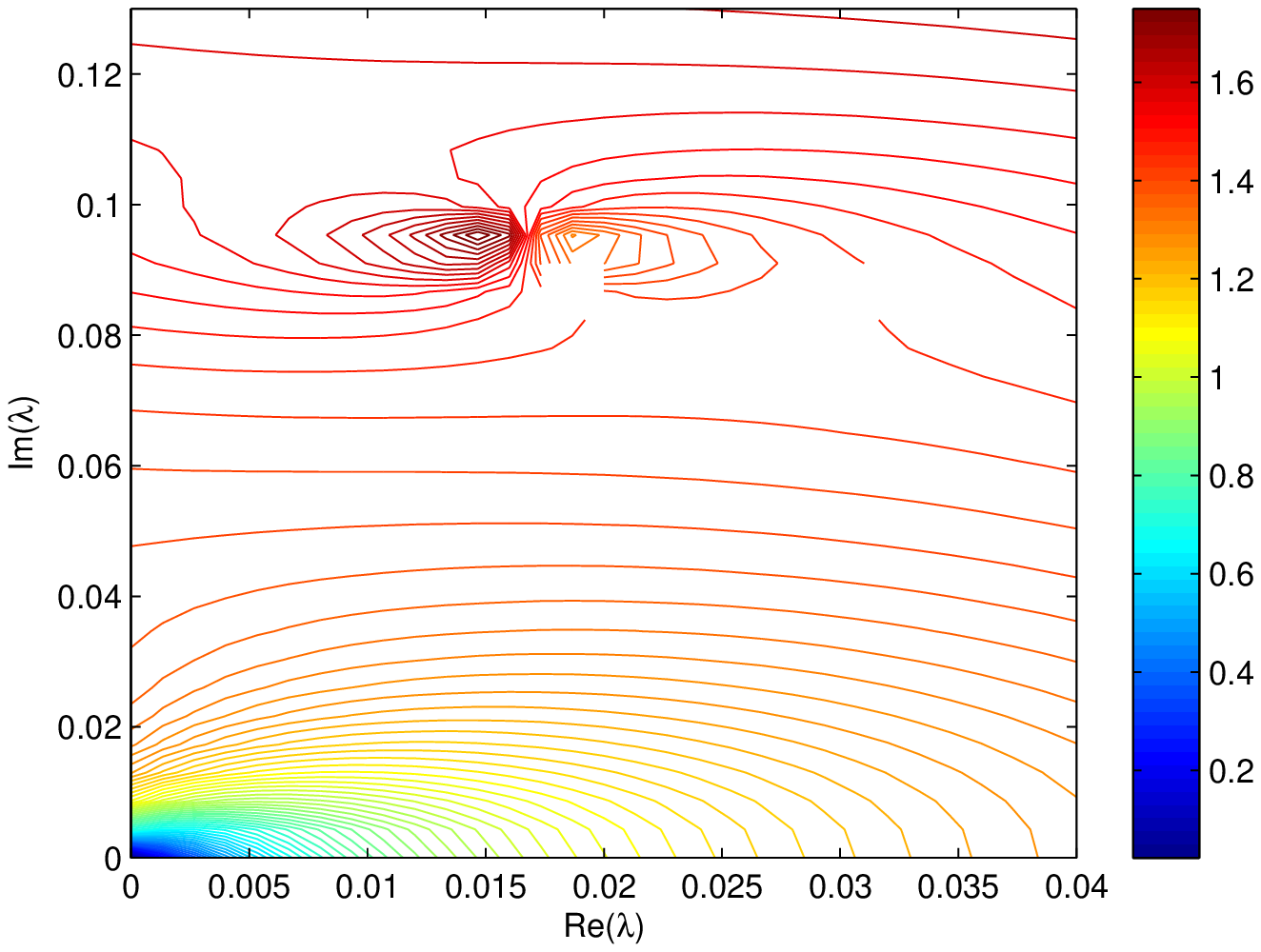}\\
\includegraphics[width=0.32\textwidth]{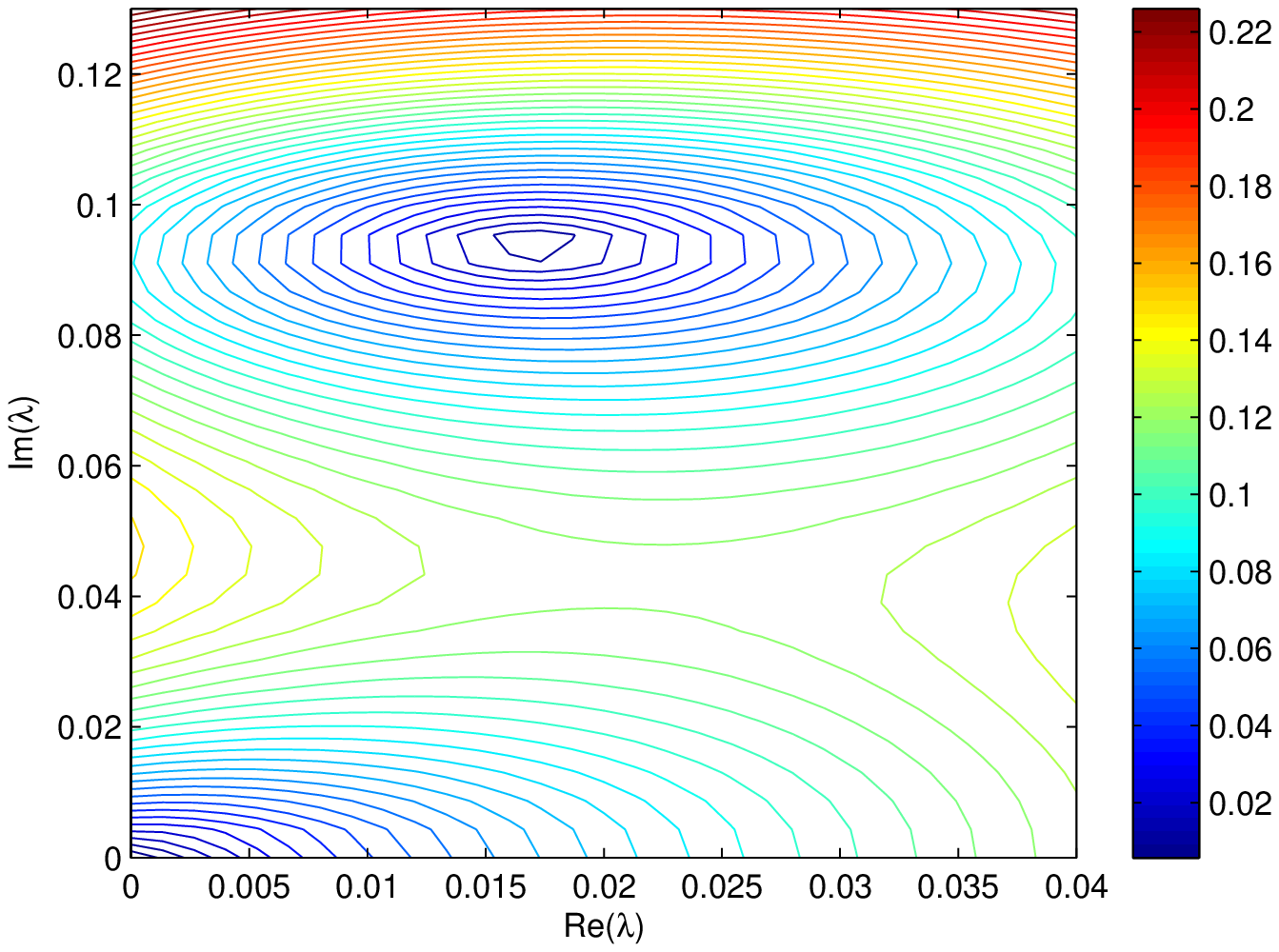}
\includegraphics[width=0.32\textwidth]{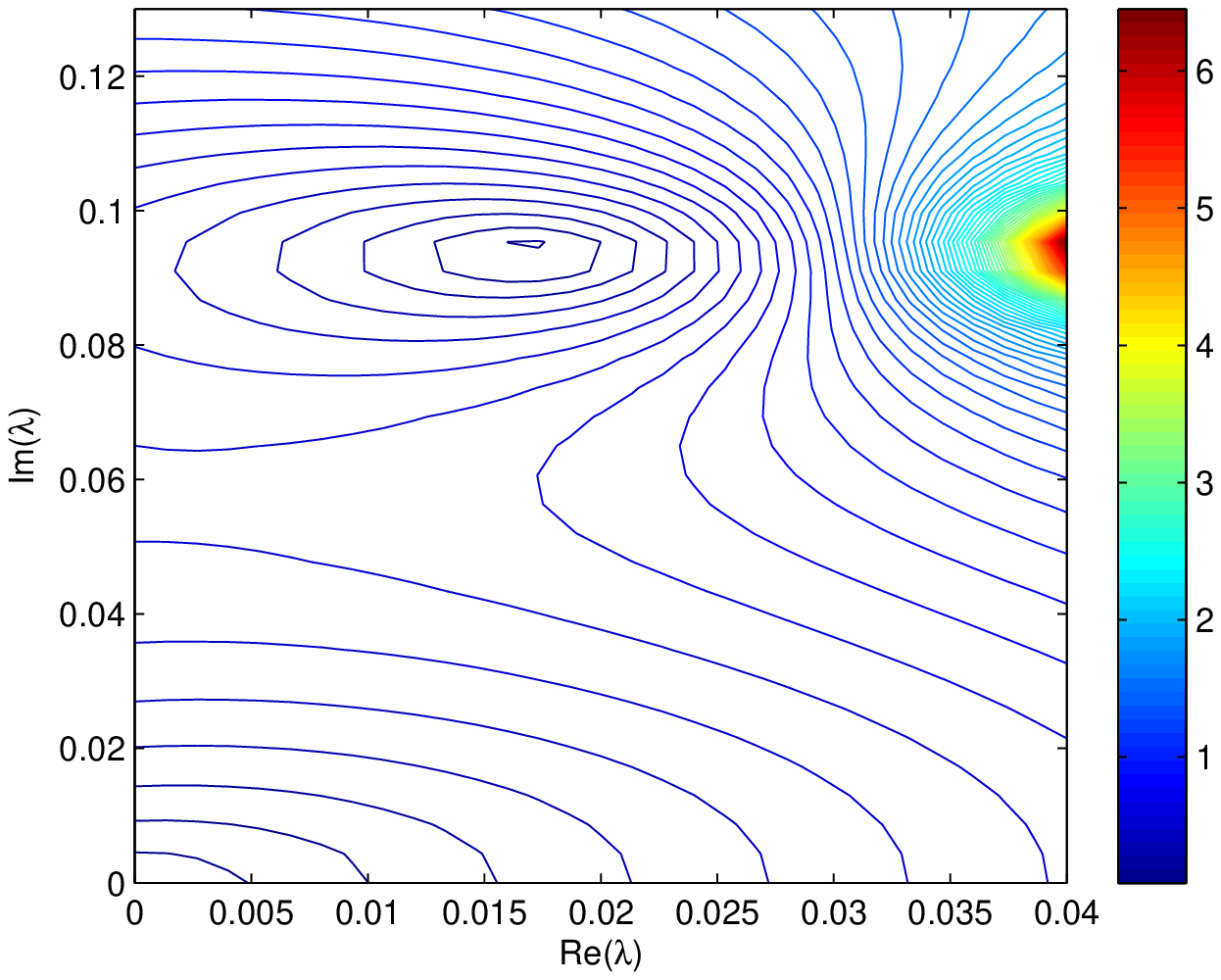}
\includegraphics[width=0.32\textwidth]{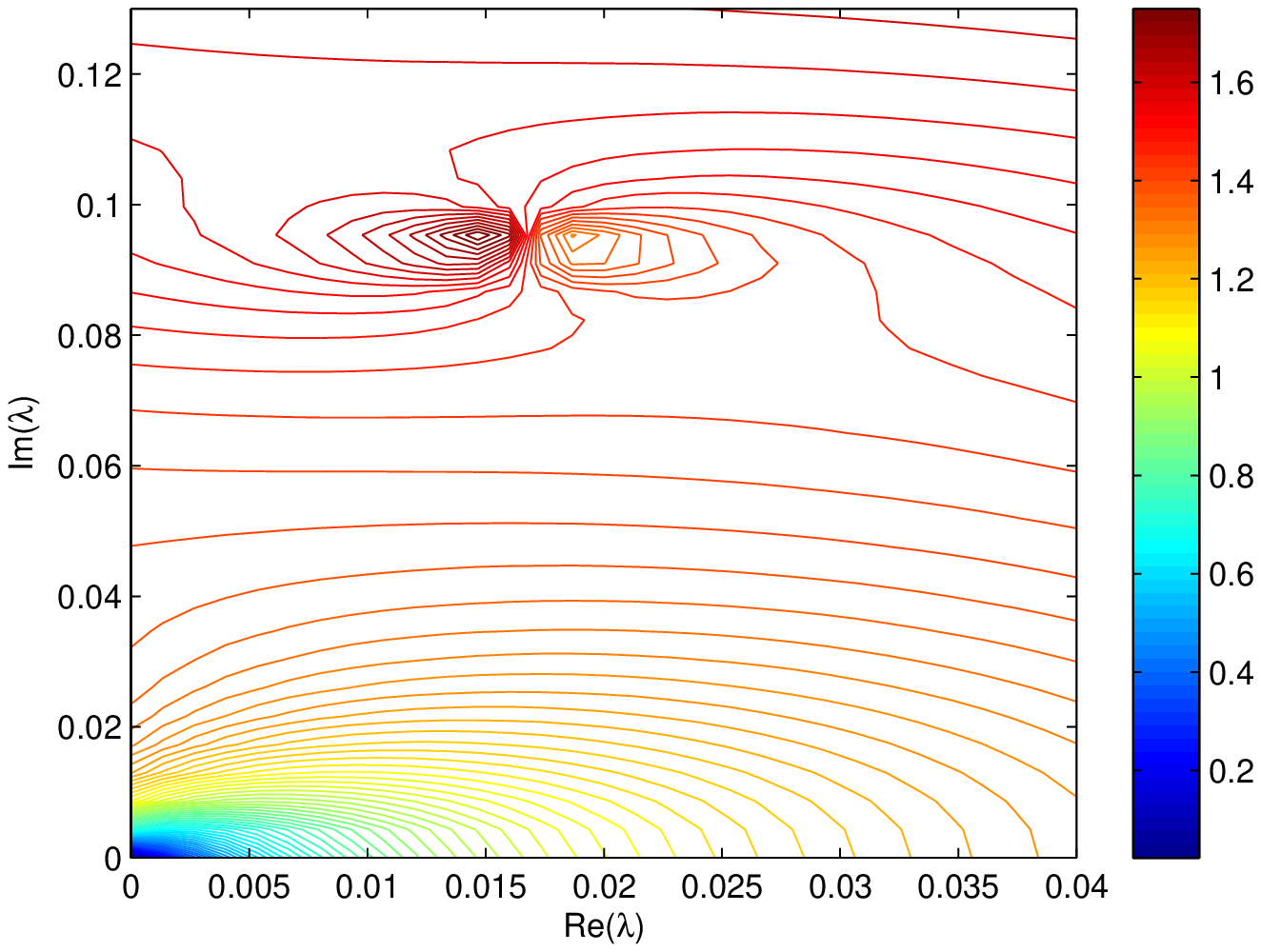}\\
\includegraphics[width=0.32\textwidth]{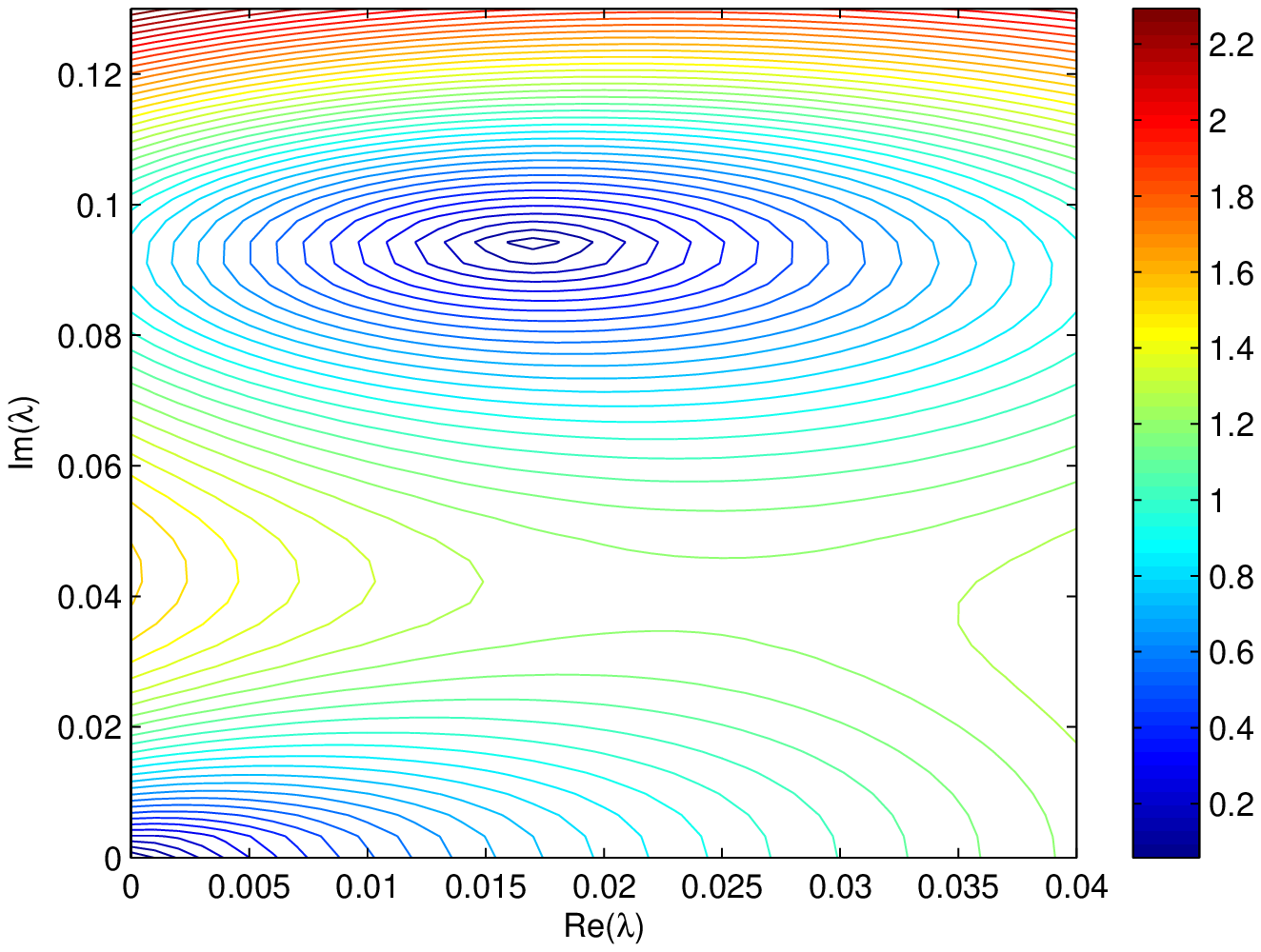}
\includegraphics[width=0.32\textwidth]{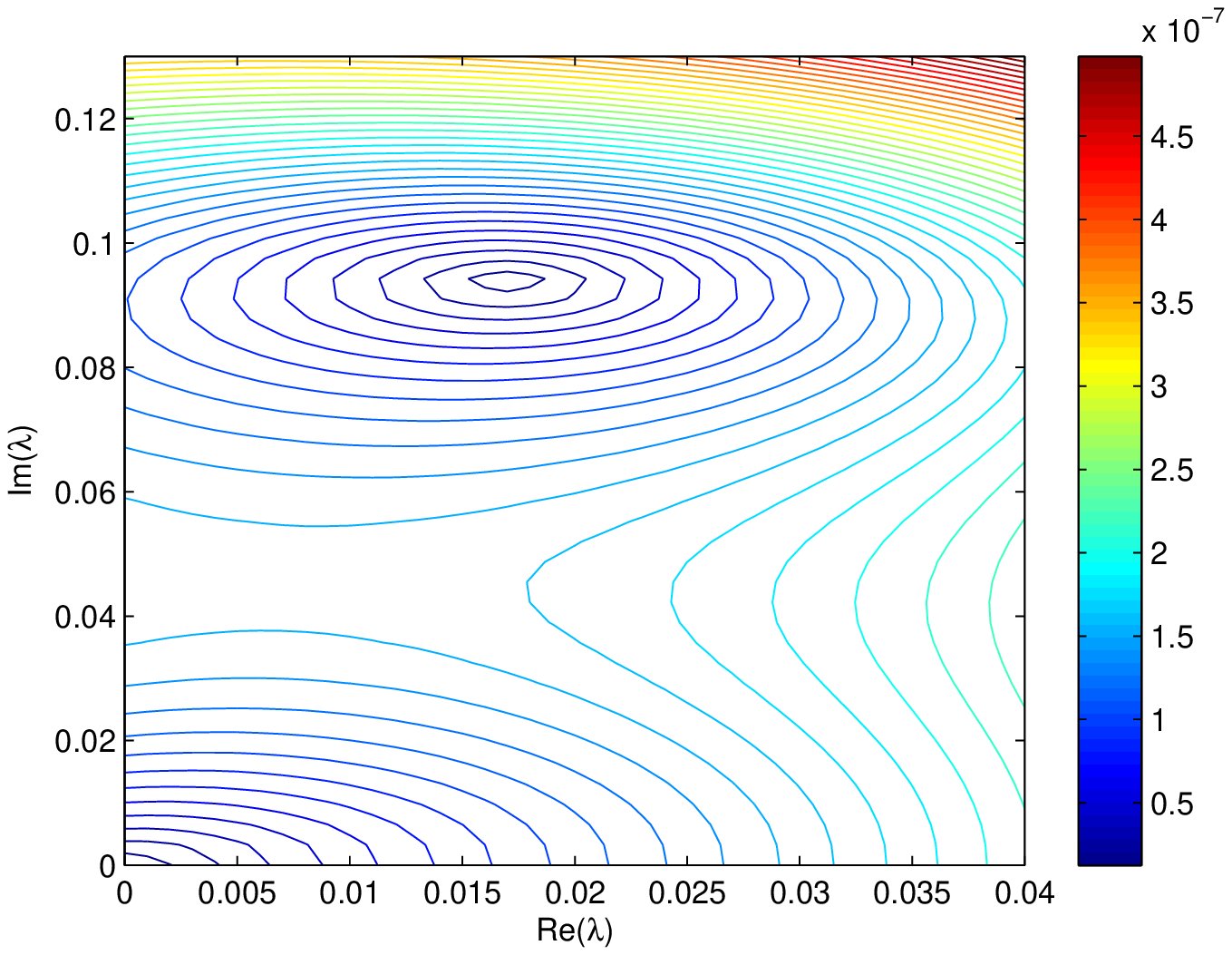}
\includegraphics[width=0.32\textwidth]{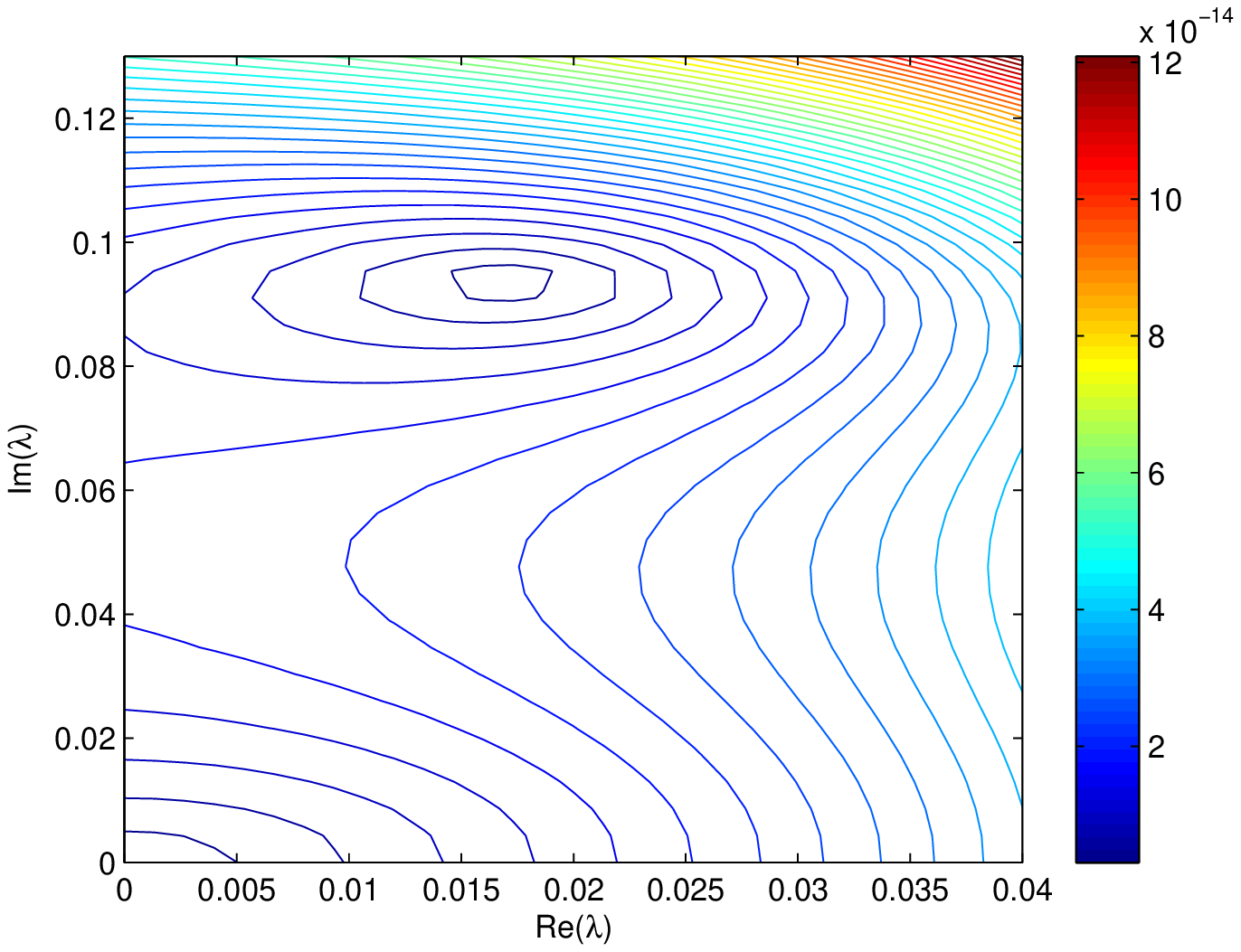}\\
\includegraphics[width=0.32\textwidth]{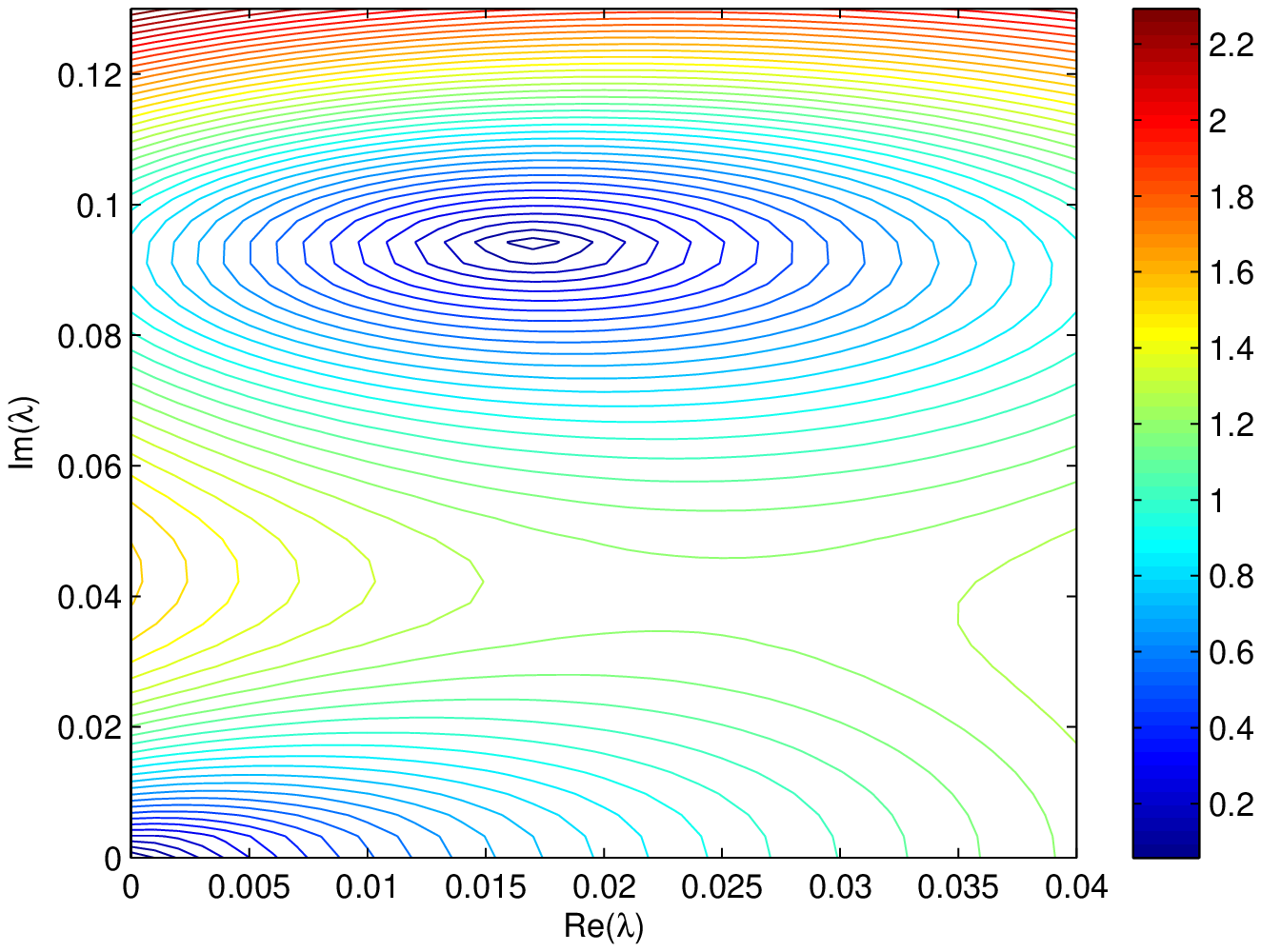}
\includegraphics[width=0.32\textwidth]{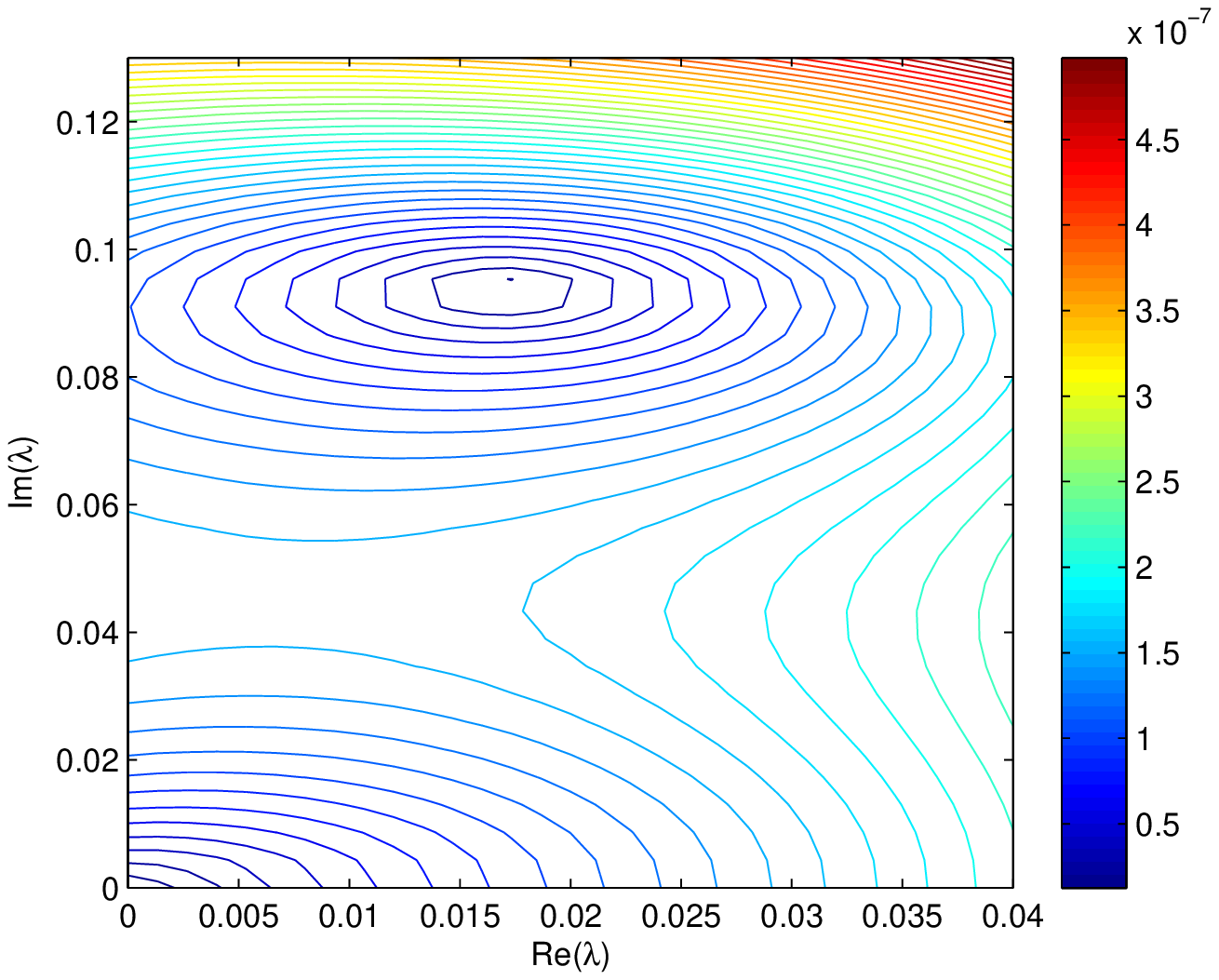}
\includegraphics[width=0.32\textwidth]{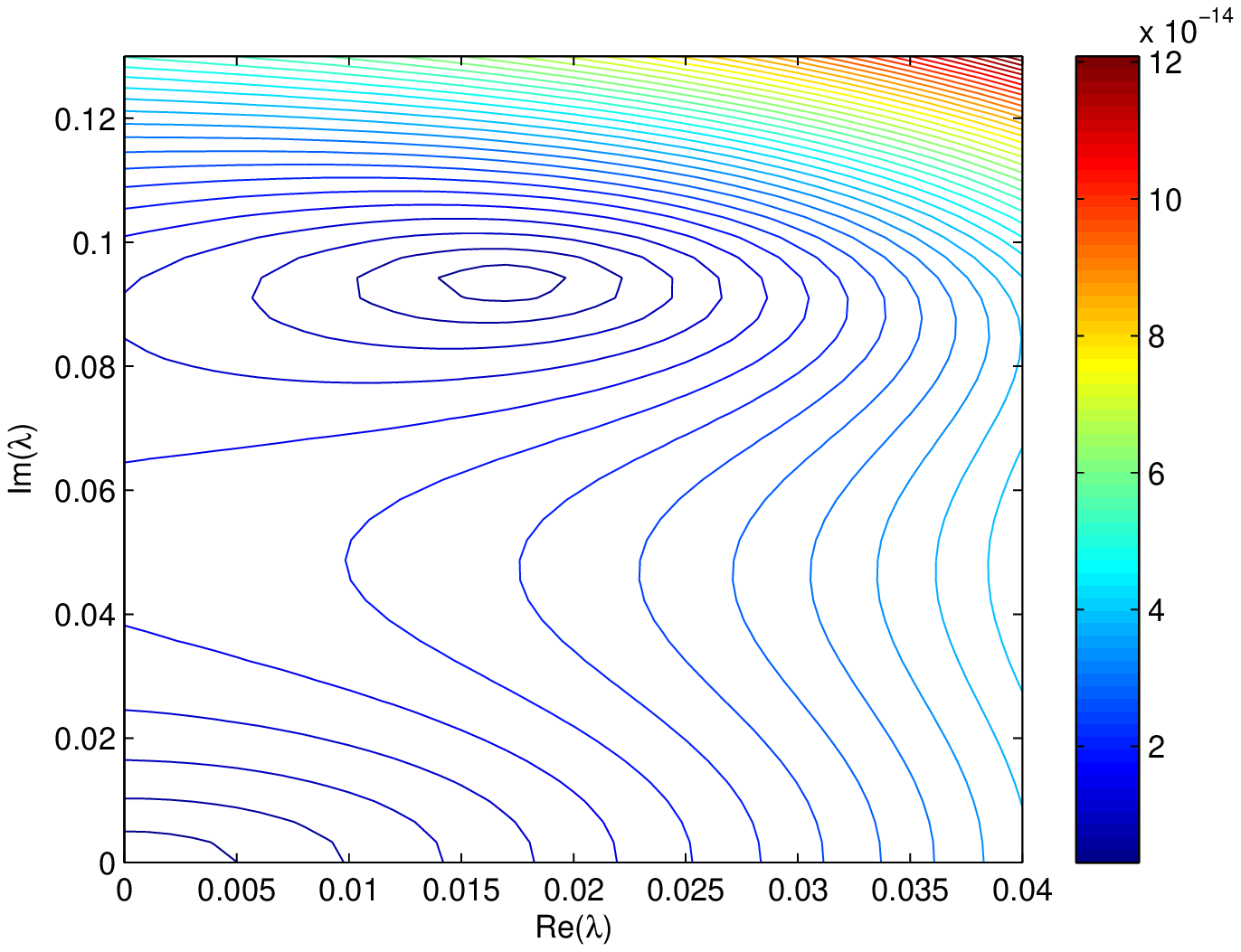}
\caption{Contour lines of $|D(\lambda)|$ for $\delta=0.1$ and $m=9$ when using 
the Riccati-RK, M\"obius--Magnus, CO-RK and GGEM-LG methods 
(order: top three down to bottom three), 
matching at positions $x_\ast=-8,0,+8$ (left to right).}
\label{autocat:contours}
\end{figure}

We compare the Riccati-RK, M\"obius--Magnus, CO-RK and GGEM-LG methods
in Figure~\ref{autocat:error} where we plot the absolute error in the eigenvalue 
vs the stepsize (upper panel) and also vs cputime (lower panel).
The eigenvalue in question is that in the first quadrant in
Figure~\ref{fig4} for $\delta=0.1$ and $m=9$.
Figure~\ref{autocat:error} was generated as follows.  
Starting with an initial guess lying within a small square 
around the eigenvalue, we iterated a standard root finding 
algorithm until we arrived in a square (containing the eigenvalue) 
which was smaller than a preset tolerance. We see in Figure~\ref{autocat:error}
that the Riccati-RK method produces a slightly better error for
a given stepsize, and is marginally more efficient than the GGEM-LG method.
The CO-RK method produces a larger error for a given computational
effort. This is not surprising, as again, the matrices in the 
Drury--Oja vector field are twice as big ($4 \times 2$) as the ones 
in the Riccati vector fields ($2 \times 2$). 

When we match at $x_\ast=-7$, there is little to distinguish the 
Riccati-RK, M\"obius--Magnus, CO-RK and GGEM-LG methods. We compare
all four methods for different matching positions $x_\ast$
in Figure~\ref{autocat:errormatch}, which was generated using the 
same root finding criteria as for Figure~\ref{autocat:error}, except 
all the methods used $N=256$ steps. Note that the errors in the 
M\"obius--Magnus and GGEM-LG methods are uniform for any matching values
in the range $[-10,10]$. The CO-RK method error doesn't vary that
much either and is slightly larger. Note that no values 
are plotted for the CO-RK method at the matching points $x_\ast=8,10$. 
In these cases the classical Runge--Kutta method applied to the 
Drury--Oja vector field is unstable for $N=256$ steps for some 
$\lambda$-values close to the eigenvalue. This problem is resolved
by increasing the number of steps to $N=512$.
For a range of matching positions
roughly in $[-10,0]$, there are no singularities of the Riccati-RK method 
in the left and right-hand integration intervals for values of the 
spectral parameter $\lambda$ close to the eigenvalue. Indeed for this
range of matching positions the Riccati-RK method delivers the best
accuracy. However for matching positions outside this range, for values of the 
spectral parameter $\lambda$ not far from the eigenvalue, the
Riccati-RK solution does have a singularity for some matching points 
(which we can see in the contour plots in Figure~\ref{autocat:contours}). 
This makes the eigenvalue-searching algorithm fail---indicated by no error points
for those matching position values. We also do not show the error
for the Riccati-RK method for the matching points $x_\ast=-8,-10$,
as there are singularities in the Evans function close to the real 
axis for these matching points (again see Figure~\ref{autocat:contours}). 
This means that we cannot for example, apply the argument principle 
in the first quadrant, though starting sufficiently close to the eigenvalue 
we can still use the Riccati-RK method as part of a root-finding algorithm 
to determine the eigenvalue.

\begin{figure}
\includegraphics[width=0.45\textwidth]{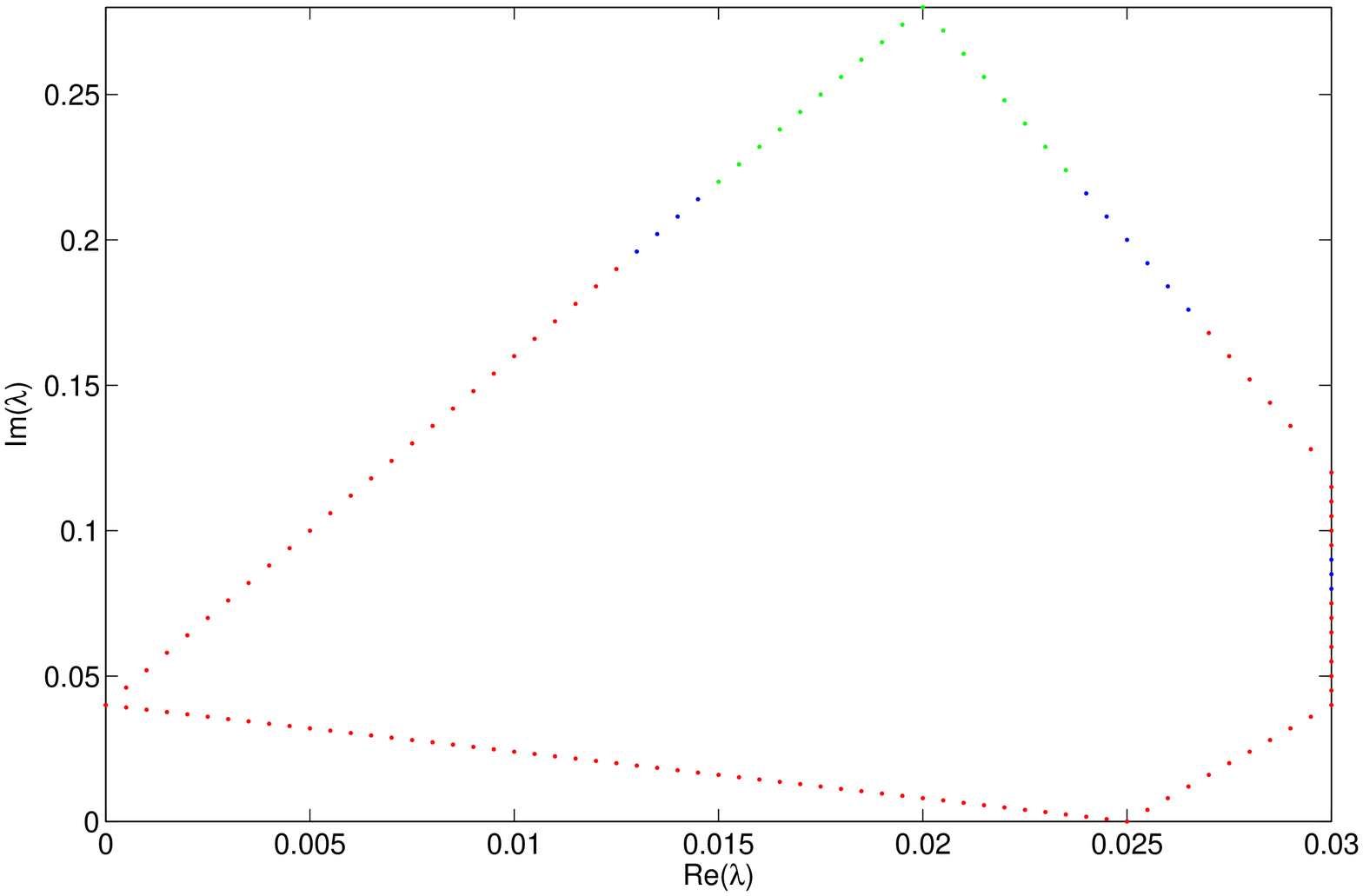}
\includegraphics[width=0.45\textwidth]{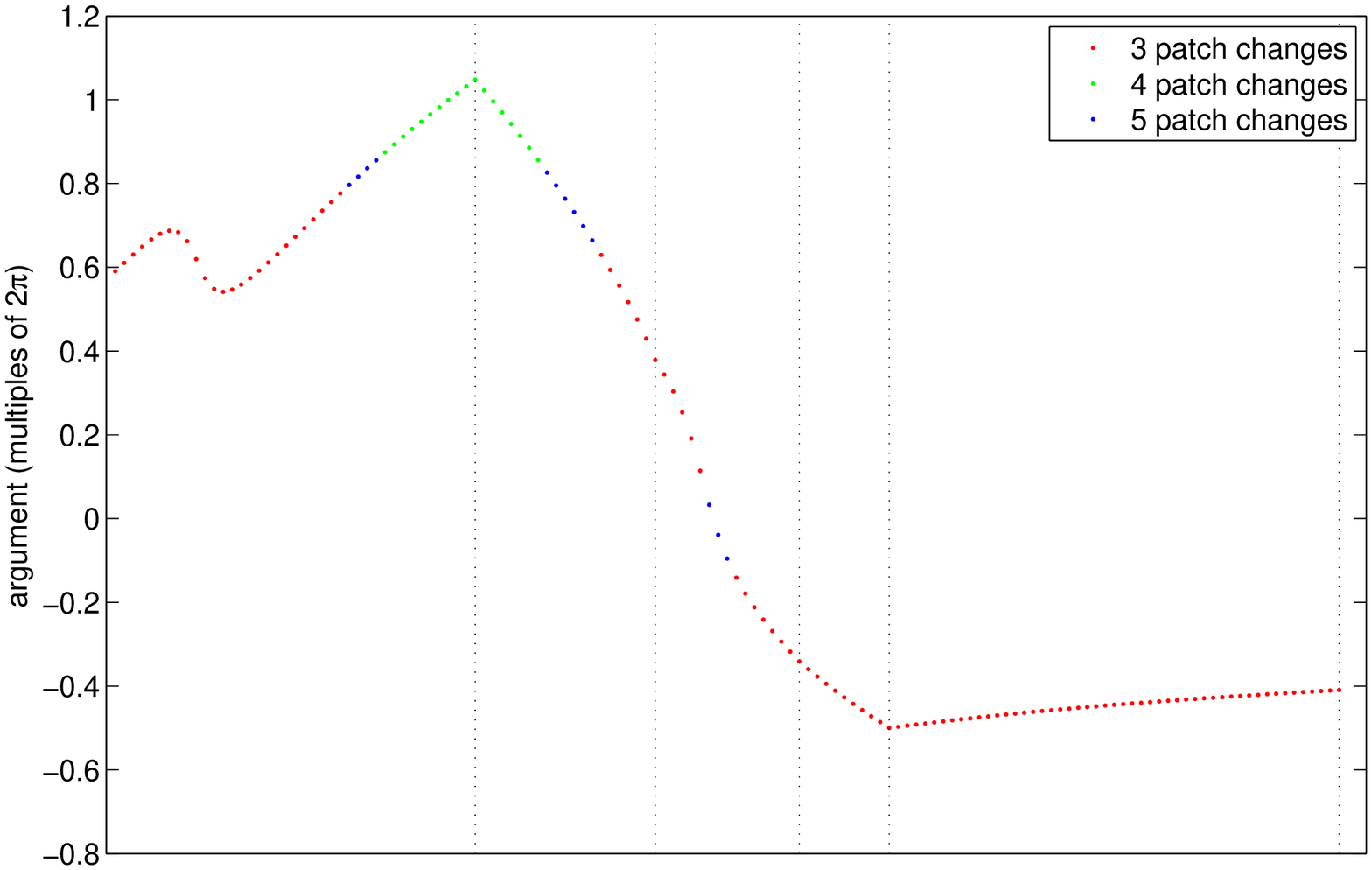}\\
\includegraphics[width=0.45\textwidth]{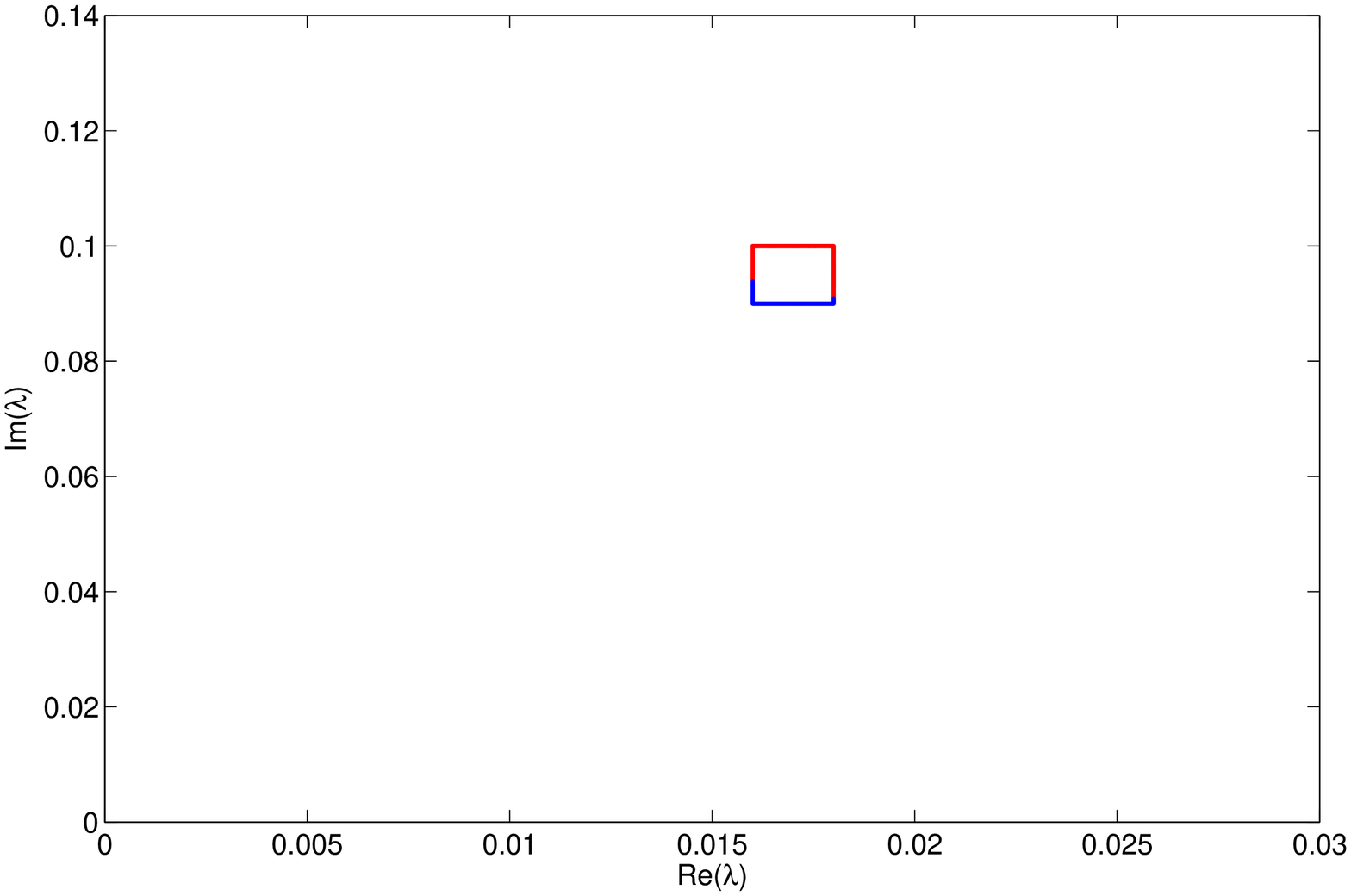}
\includegraphics[width=0.45\textwidth]{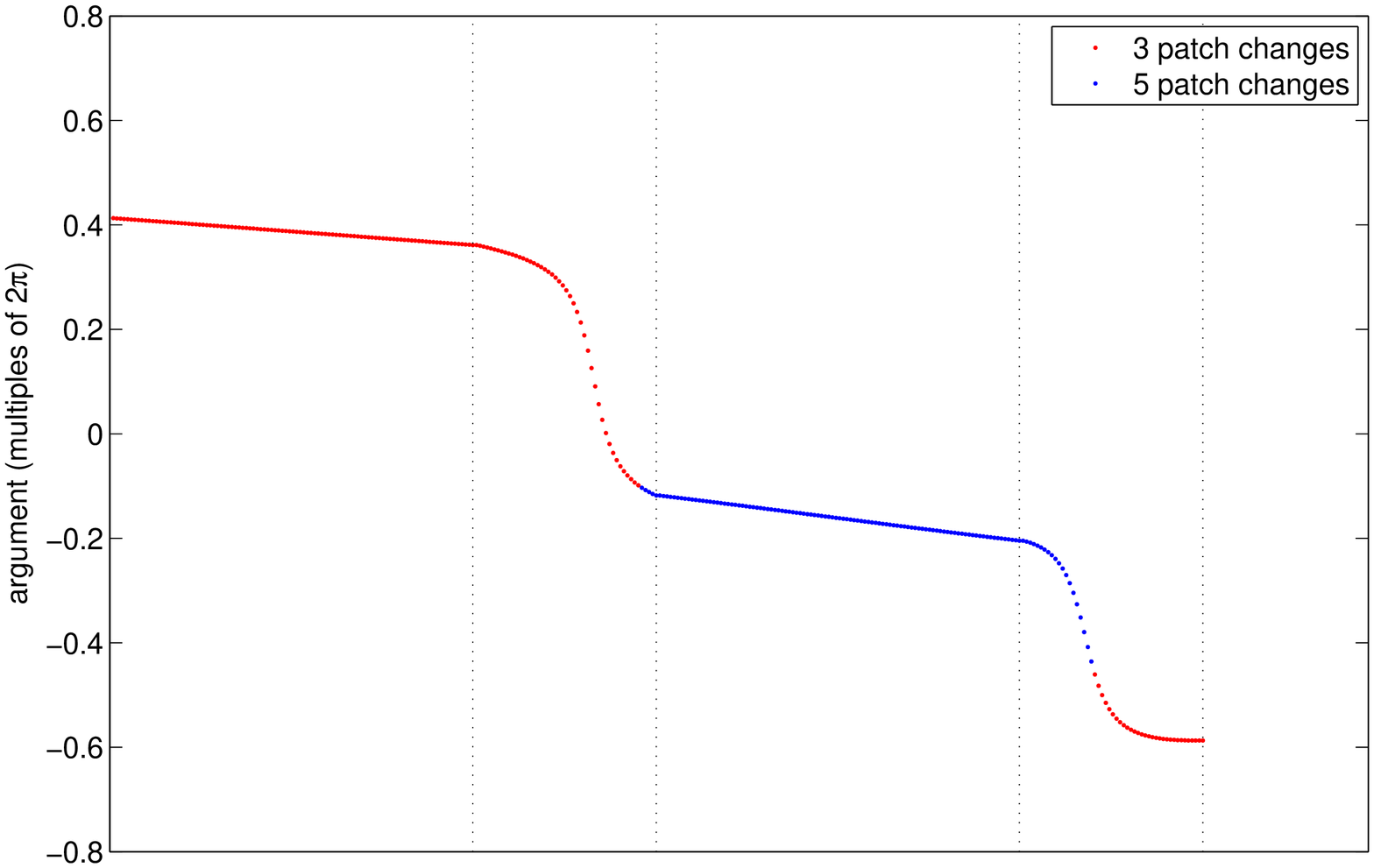}\\
\includegraphics[width=0.45\textwidth]{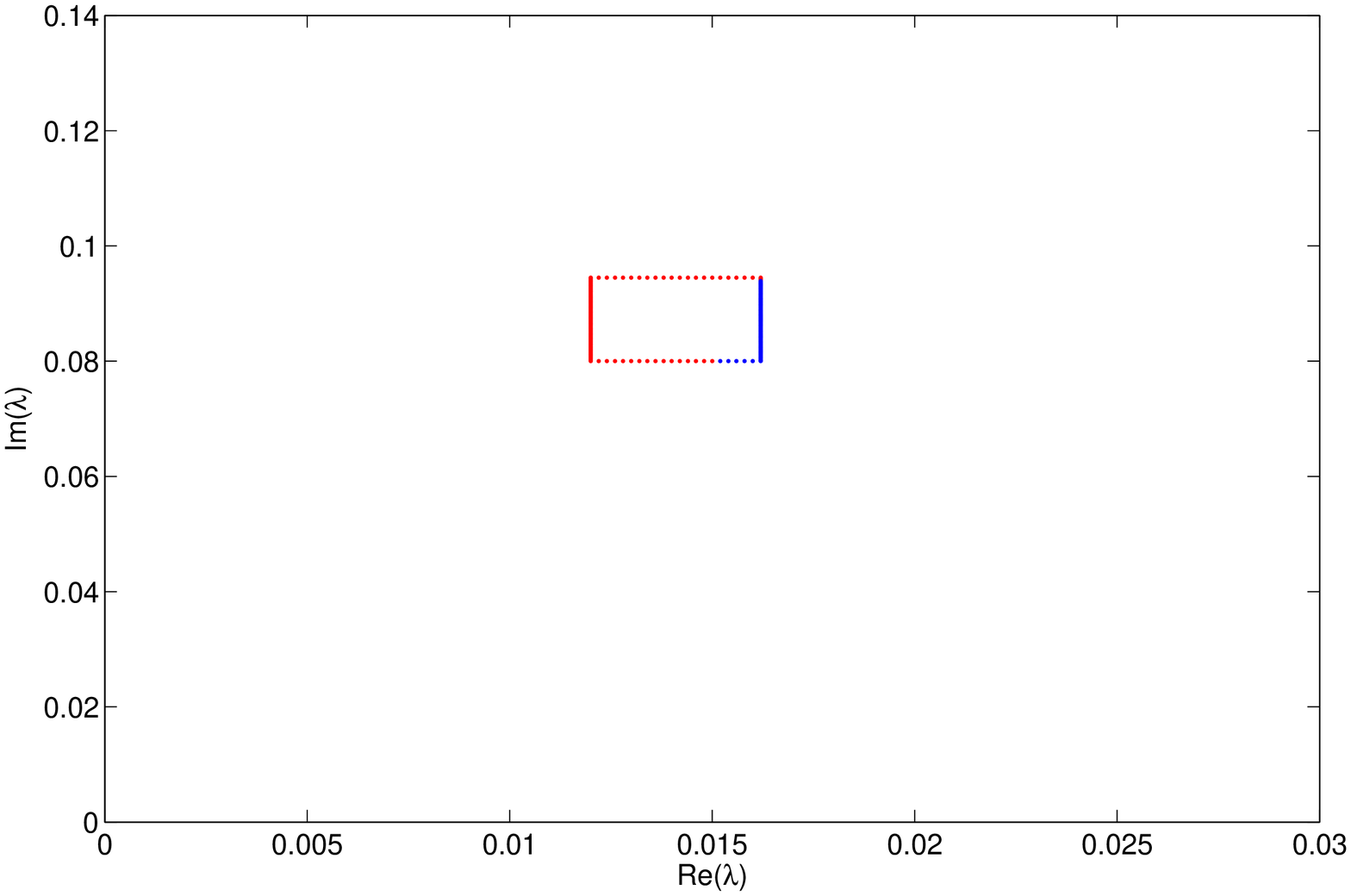}
\includegraphics[width=0.45\textwidth]{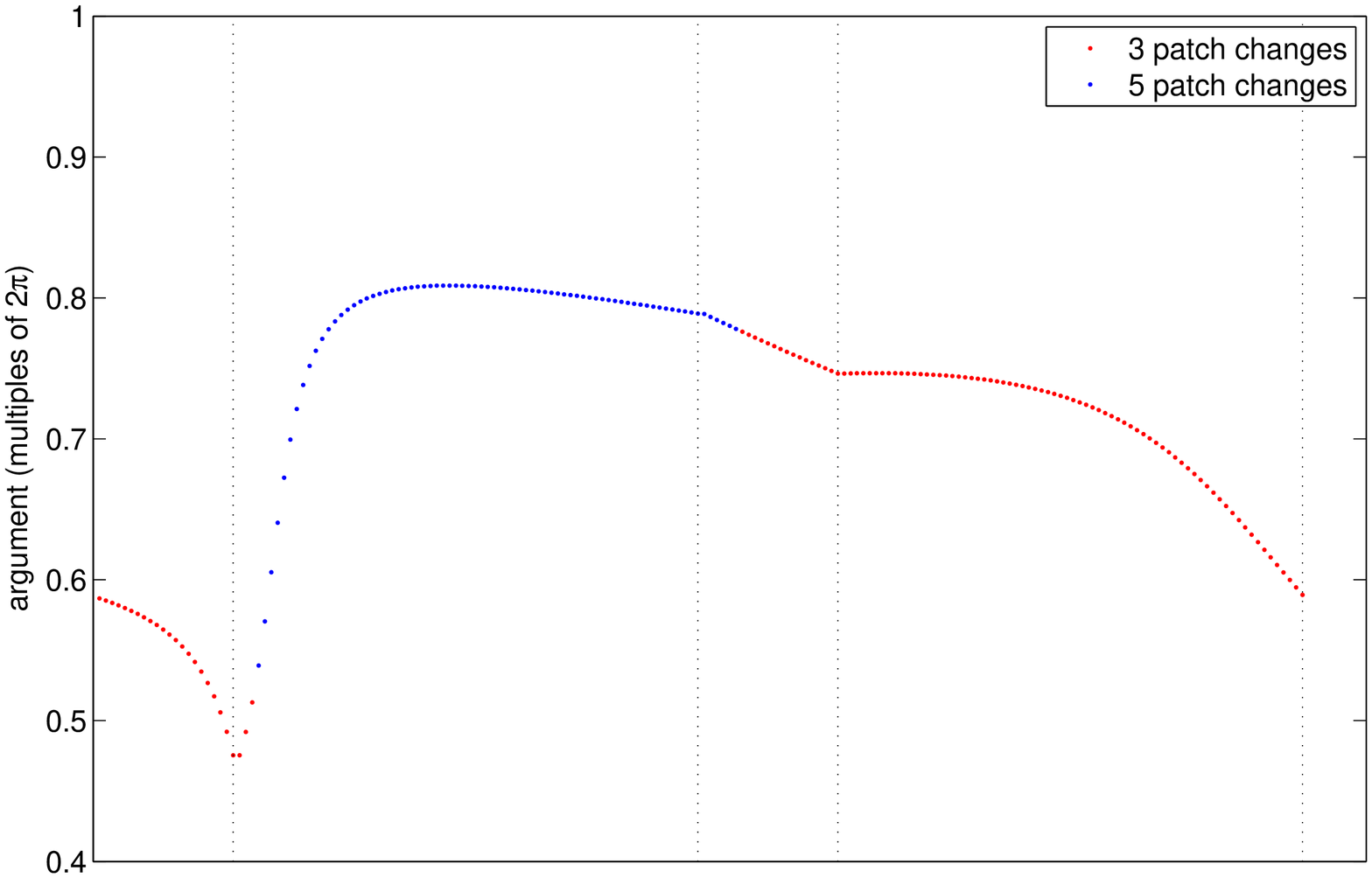}
\caption{We show three different closed contours in the
first quadrant of the complex $\lambda$-plane in 
each of the left panels. In the top two left panels,
the contour encloses the eigenvalue in that quadrant
(found in Figure~\ref{autocat:contours} for $\delta=0.1$, $m=9$).
In the corresponding three panels on the right we show how
the argument (in multiples of $2\pi$) of the 
Evans function $D(\lambda)$ changes
as we perform a complete circuit of the contour. The Evans
function was computed using the GGEM-LG method matching at
$x_\ast=+14$. The number of patch changes performed for each 
fixed $\lambda$ value are indicated.}
\label{autocat:contourint}
\end{figure}

Figure~\ref{autocat:contours} 
shows the contour lines of $|D(\lambda)|$ for $\delta=0.1$ and $m=9$,
close to the eigenvalue in the first quadrant, when using 
the Riccati-RK, M\"obius-Magnus, CO-RK and GGEM-LG methods, respectively,
and matching at three positions $x_\ast=-8,0,+8$.
We see that the CO-RK and GGEM-LG methods show the least sensitivity
to the choice of matching position and produce smooth contour plots for 
all three matching points. The contour plots across shape and scale 
look very similar for both these methods. 
By contrast the Riccati-RK and M\"obius-Magnus methods
appear to develop singularities close to the eigenvalue when the 
matching position $x_\ast$ is $0$ or $+8$.

Lastly in Figure~\ref{autocat:contourint} we demonstrate
the argument principle for counting zeros of the Evans
function inside closed contours. 
We computed the Evans function using the GGEM-LG method
and matched at $x_\ast=+14$. As expected, if the closed
contour in the complex $\lambda$-plane encloses the 
eigenvalue, then the change in the argument of the Evans
function around the complete contour is one (once we have
accounted for the $2\pi$ factor in the argument principle).
We also show, for each fixed $\lambda$ value, the number
of patch changes that occured as we integrated from $x=-14$
through to $x=+14$.

\subsection{Ekman boundary layer}
The third test system is a boundary layer flow over a flat plate 
which is infinitely extended in the $x$ and $y$ direction and 
rotates around the half infinite $z$-axis with a given rotational speed.
Linear stability of the Ekman boundary layer has been investigated in
Allen and Bridges~\cite{AB2} and Allen~\cite{A} using the compound
matrix method. The flow is governed by the continuity equation
$u_x + v_y + w_z = 0$,
and the Navier-Stokes equations in a co-rotating frame 
\begin{align*}
  u_t + u u_x + v u_y + w u_z + \tfrac{1}{\Ro} p_x - \tfrac{2}{\Ro} v &=
  \tfrac{1}{\Rey} (u_{xx} + u_{yy}) + \tfrac{1}{\Ro} u_{zz},\\
  v_t + u v_x + v v_y + w v_z + \tfrac{1}{\Ro} p_y + \tfrac{2}{\Ro} u &=
  \tfrac{1}{\Rey} (v_{xx} + v_{yy}) + \tfrac{1}{\Ro} v_{zz},\\
  w_t + u w_x + v w_y + w w_z + \tfrac{1}{\Ro \Ek} p_z &=
  \tfrac{1}{\Rey} (w_{xx} + w_{yy}) + \tfrac{1}{\Ro} w_{zz}.
\end{align*}
Here $\Rey$, $\Ro$ and $\Ek$ denote the Reynolds, Rossby and  
Ekman numbers, respectively.

After non-dimensionalization and setting $\Rey=\Ro$, 
$\Ek=1$, the linear stability of the boundary
layer is determined by the eigenvalues $\lambda$ 
of the linear problem $Y' = A(z;\lambda) Y$, 
where (see Allen~\cite[p.~176]{A})
\begin{equation*}
  A(z;\lambda) = \begin{pmatrix}
    0& 1& 0& 0& 0& 0\\
    0& 0& 1& 0& 0& 0\\
    0& 0& 0& 1& 0& 0\\
    -a(z,\lambda) & 0& b(z, \lambda) & 0& 0& -2\\
    0& 0& 0& 0& 0& 1\\
    i \gamma \Rey \U_z(z) & 2& 0 & 0 & b(z,\lambda)-\gamma^2& 0
  \end{pmatrix}
\end{equation*}
and
\begin{align*}
  \U(z)&=\cos(\epsilon)\bigl(1-\exp(-z)\cos(z)\bigr) 
         +\sin(\epsilon)\exp(-z)\sin(z),\\  
  \U_z(z)&=\exp(-z)\bigl(\sin(z+\epsilon)+\cos(z+\epsilon)\bigr) \,,\\
  \V(z)&=-\sin(\epsilon)\bigl(1-\exp(-z)\cos(z)\bigr)+\cos(\epsilon)\exp(-z)\sin(z)\,, \\
  \V_{zz}(z)&=-2\exp(-z)\cos(z+\epsilon), \\
  a(z,\lambda)&=\gamma^4+i\Rey\gamma^2\bigl(\gamma \V(z)-i \lambda\bigr)
  +i \gamma \mathrm{\Rey} \V_{zz}(z),\\
  b(z, \lambda) &= 2 \gamma^2+\Rey \bigl(i \gamma \V(z)+\lambda\bigr).
\end{align*}
Here the parameters $\gamma$ and $\epsilon$ represent the radial and
angle components, respectively, of a polar coordinate parameterization
of horizontal wavenumbers associated with the $x$ and $y$ directions---see
Allen and Bridges~\cite{AB2} for more details.

We choose the fixed coordinate patch identified by $\ib^+=\{1,2,3\}$
to integrate the corresponding Riccati equation from $z=\ell_+$ to $z_\ast=0$.
The boundary condition for the rigid wall at $z_\ast=0$ 
as given in Allen and Bridges is
\begin{equation*}
  Y_1(0;\lambda)=Y_2(0;\lambda)=Y_5(0;\lambda)=0.
\end{equation*}
If 
\begin{equation*}
Y^-=\begin{pmatrix}
  1 & 0 & 0 & 0 & 0 & 0\\ 0 & 1 & 0 & 0 & 0 & 0 \\ 0 & 0 & 0 & 0 & 1 & 0
\end{pmatrix}\,,
\end{equation*}
then the boundary conditions are equivalent to 
\begin{equation*}
  \det\bigl(Y^-\cdot Y^+(0;\lambda)\bigr)=0\,.
\end{equation*}
We thus compute the Evans function: 
\begin{equation*} 
D(\lambda;z_\ast)= \det\Biggl(
 Y^-\cdot\begin{pmatrix} I \\  \hat y \end{pmatrix}\Biggr)
  = \det \begin{pmatrix}
    1 & 0 & 0 \\ 0 & 1 & 0 \\ \hat y_{21} & \hat y_{22} & \hat y_{23} 
  \end{pmatrix} \equiv \hat y_{23}(z_\ast;\lambda).
\end{equation*}

\begin{figure}
  \includegraphics[width=0.45\textwidth]{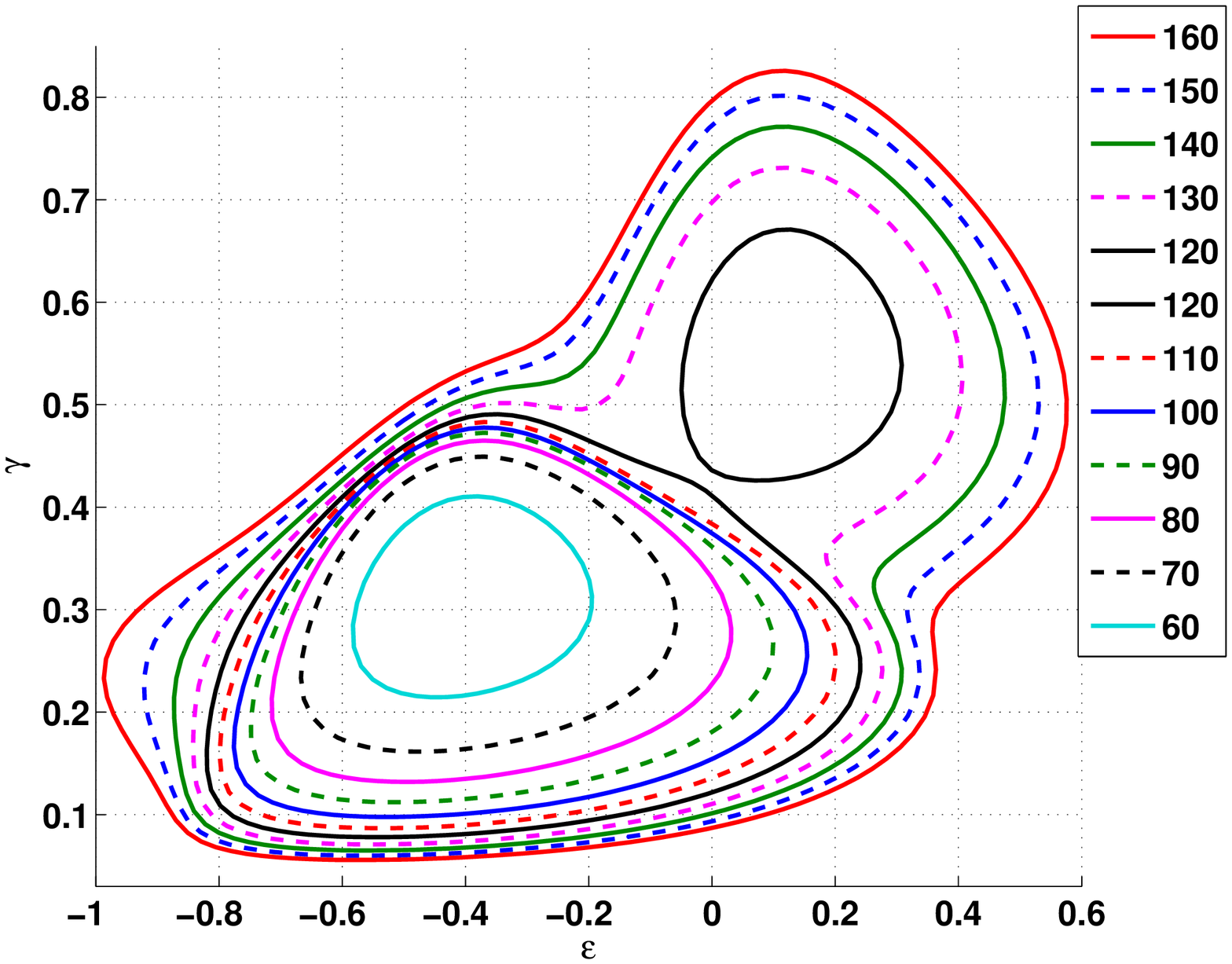}
  \includegraphics[width=0.45\textwidth]{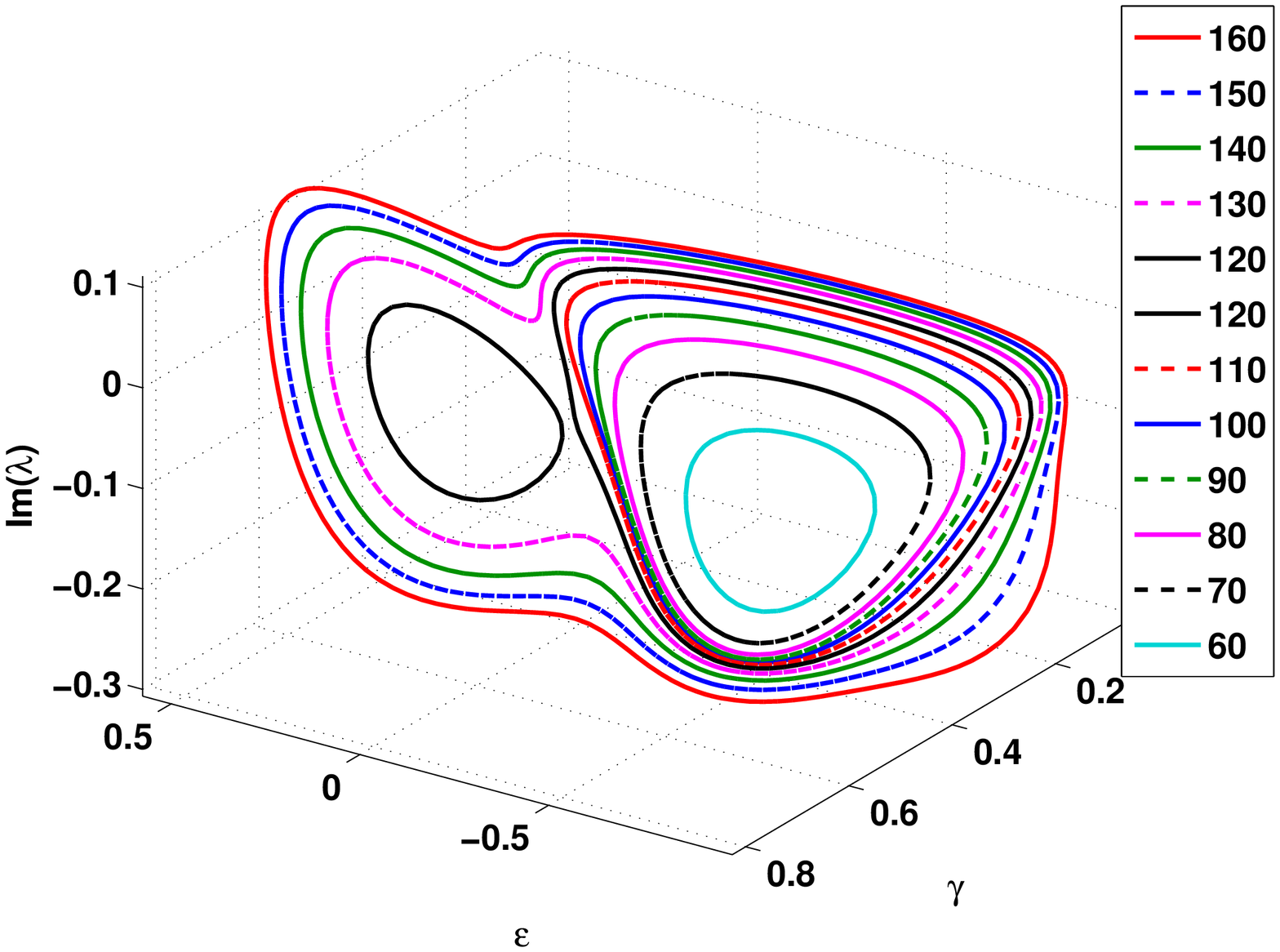}
  \caption{Neutral curves for rigid wall, $\Rey$ fixed (values indicated).}
\label{fig:ncurves_Re}
\end{figure}

\begin{figure}
  \centering
    \includegraphics[width=0.45\textwidth]{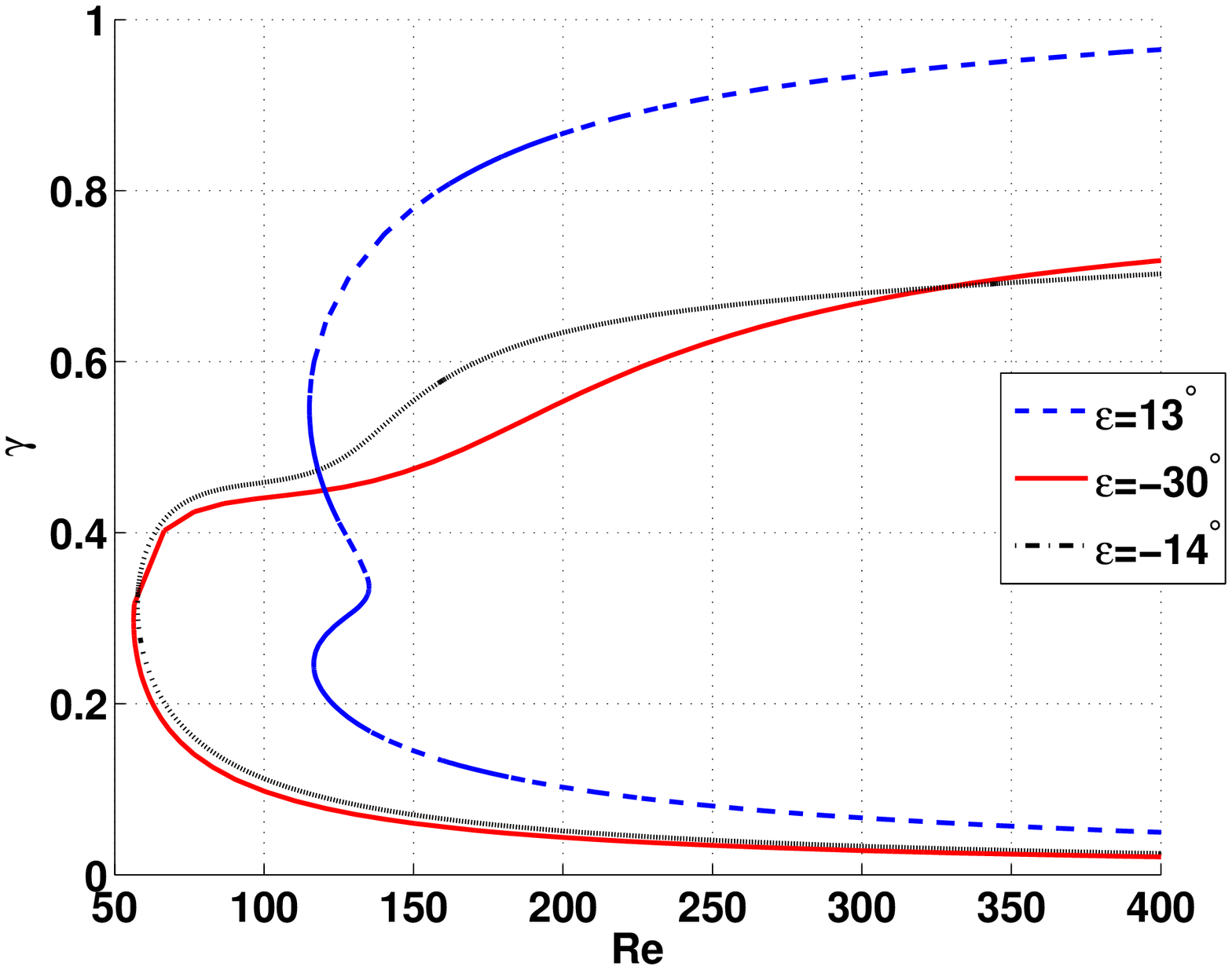}
    \includegraphics[width=0.45\textwidth]{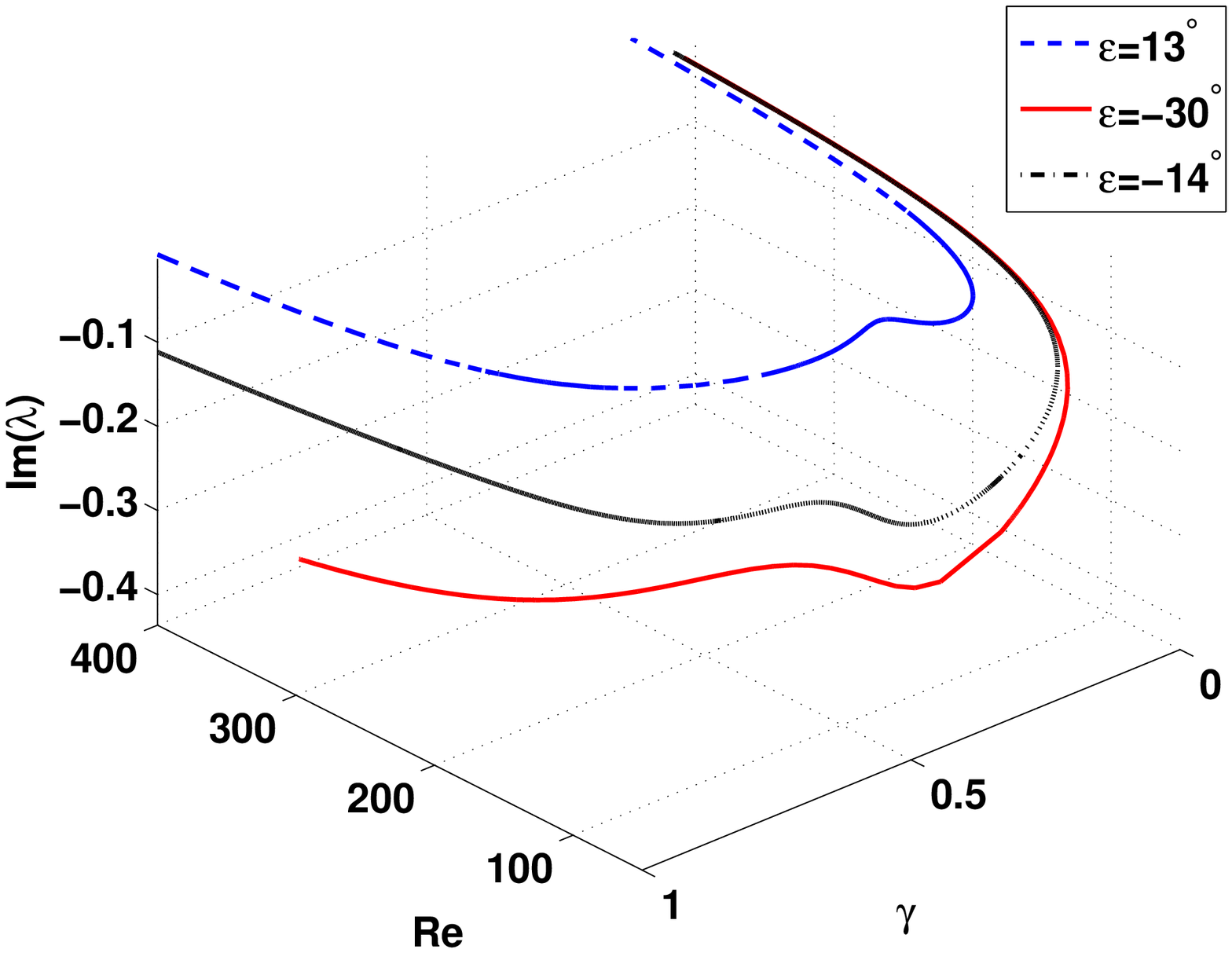}
    \caption{Neutral curves for rigid wall, $\epsilon$ fixed.}
\label{fig:ncurves_eps}
\end{figure}

\begin{figure}
    \includegraphics[width=0.32\textwidth]{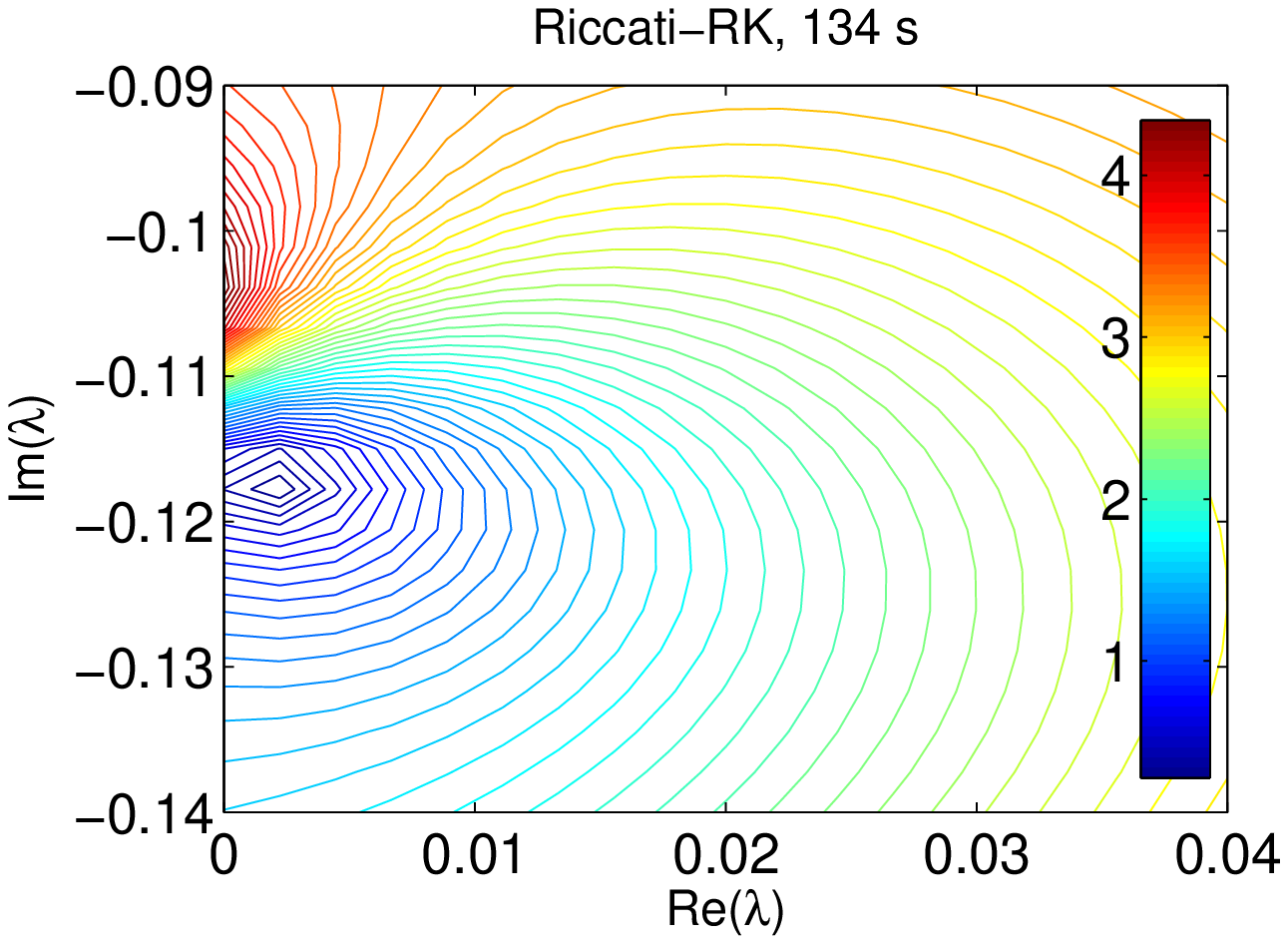}
    \includegraphics[width=0.32\textwidth]{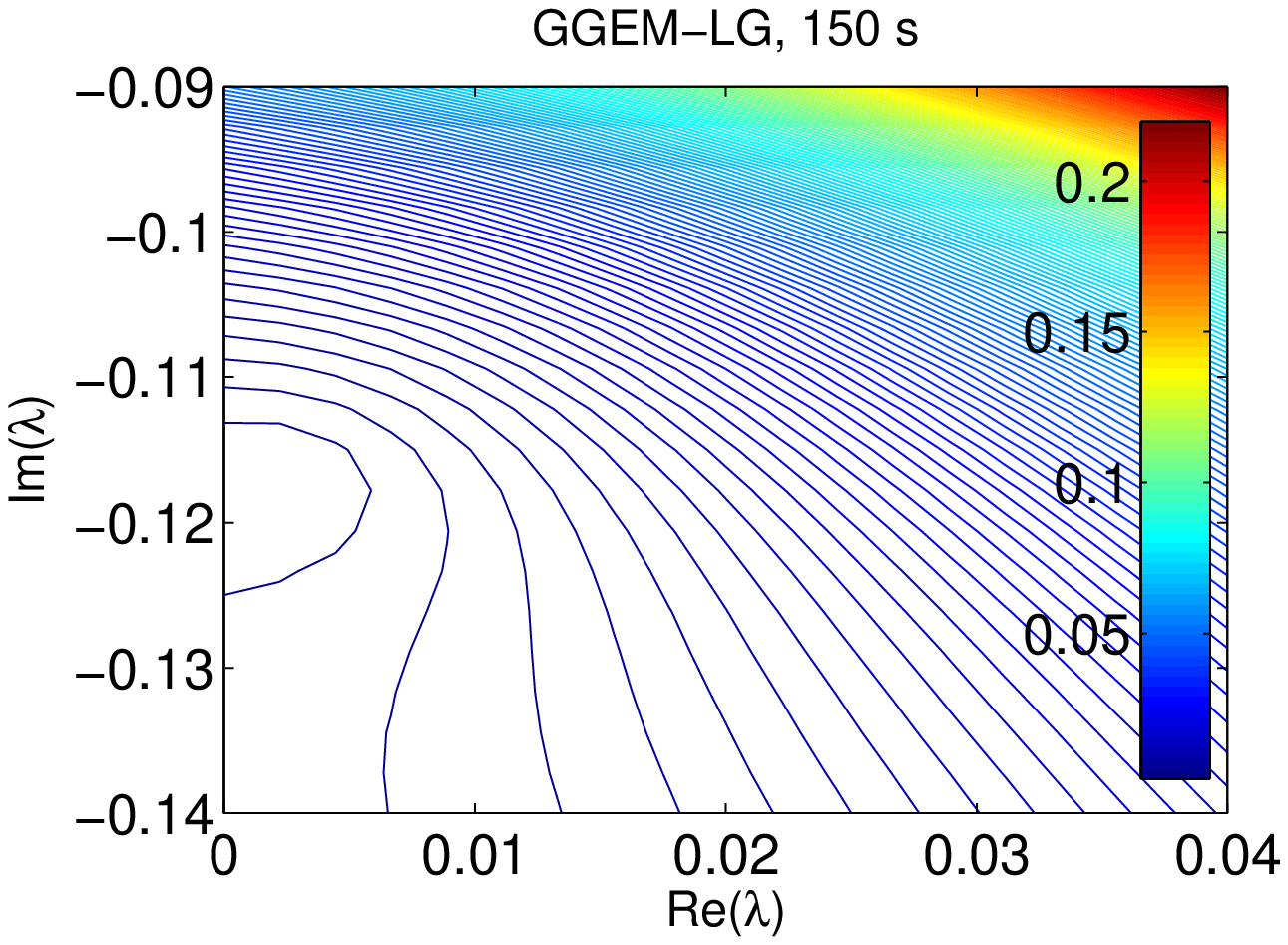}
    \includegraphics[width=0.32\textwidth]{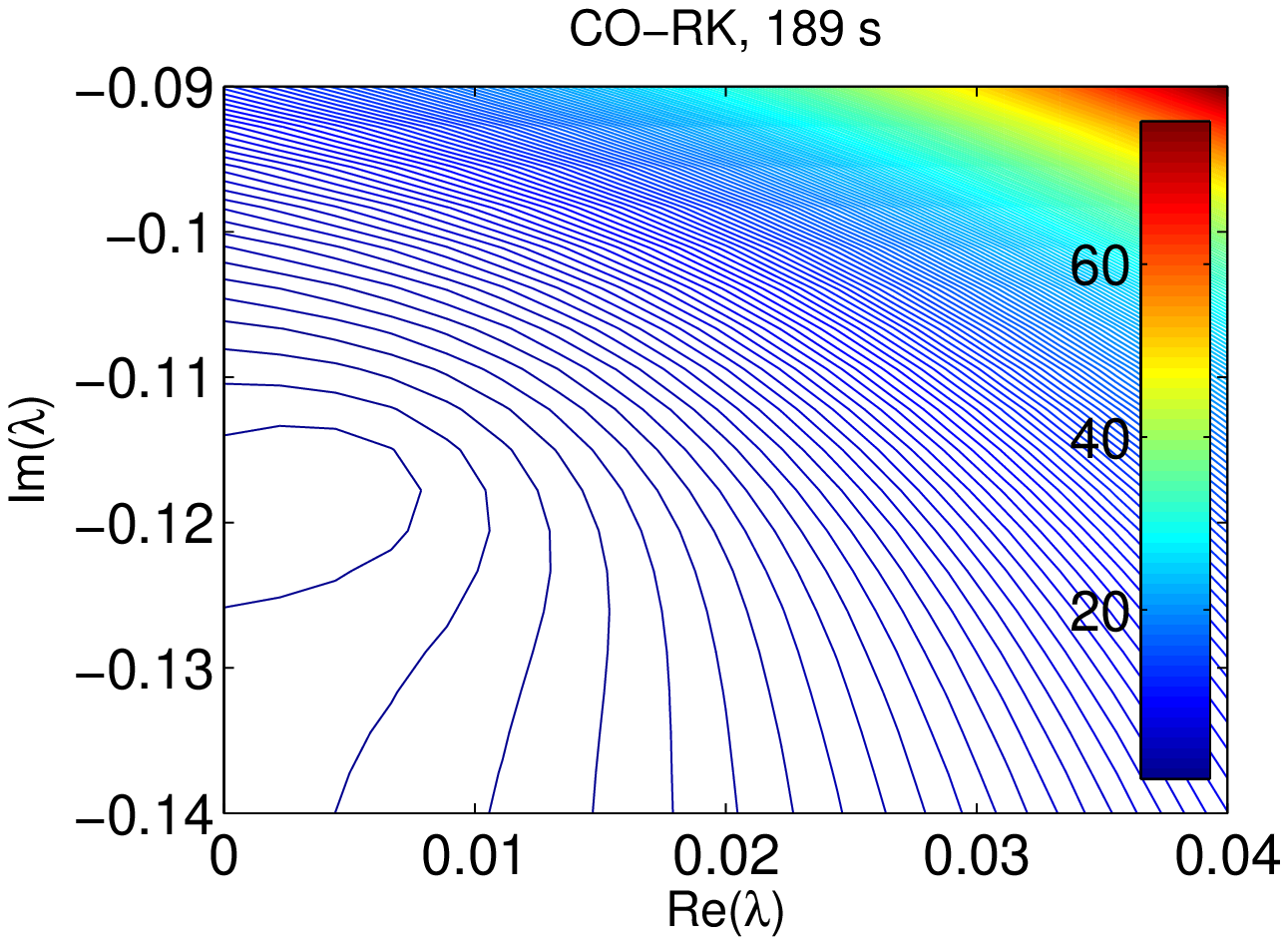}
\caption{Contour plots of $|D(\lambda;0)|$ for $\Rey=140$, 
$\epsilon=0.014156$, $\gamma=0.70575$. 
An equidistant mesh was used with $500$ intervals  
and values computed on a $20\times20$ grid in the complex $\lambda$-plane.}
\label{fig:EkmanGGEM}
\end{figure}

We computed neutral curves, i.e.\ curves in the 
$\epsilon$--$\gamma$ plane where
$Re(\lambda)=0$, using the Riccati-RK method with $\ib^+=\{1,2,3\}$ 
to compute $\hat y(0;\lambda)\in\mathbb C^{3\times3}$ 
and consequently the Evans function $D(\lambda;z_\ast)$. 
For continuation of the curves we used the Matlab package MatCont 
which uses pseudo-arclength continuation
(Dhooge, Govaerts and Kuznetsov~\cite{DGK03}).
Figures~\ref{fig:ncurves_Re} and \ref{fig:ncurves_eps}
show the neutral curves which match those 
in Allen and Bridges~\cite{AB2} and Allen~\cite{A}.
The integration of the Riccati system has been done 
with the Matlab ODE-solver \texttt{ode23s} from
$z=10$ to $z_\ast=0$ (as in Allen and Bridges)
with absolute and relative tolerances $10^{-6}$
and $10^{-4}$. The stable subspace of $A(+\infty;\lambda)$
was constructed using the Matlab eigenvalue-solver \texttt{eig} 
(we also used direct formulae for the eigenvectors
to construct analytic bases for the stable subspace but this did 
not significantly change the overall performance).

For comparison we also implemented the CO-RK and GGEM-LG methods.
We tested the performance of all three methods, in each case evaluating the 
Evans function on a $20\times20$ grid for $\lambda$ in the complex plane.
In Figures~\ref{fig:EkmanGGEM} we present contour plots of $|D(\lambda; 0)|$, 
and see that all methods find a root at $\lambda\approx 0.002-0.117\mathrm{i}$.
The computation times for a $2.4$Ghz machine were:
$134$ seconds for Riccati-RK, $150$ seconds for GGEM-LG and 
$189$ seconds for CO-RK. As a comprehensive check, we also 
implemented the compound matrix method (i.e.\ Pl\"ucker coordinates,
of which there are $20$), described in Allen and Bridges, for this 
performance test. As expected, since this 
method involves integrating a linear system of order $20$, it
was an order of magnitude slower (while giving the same results).

\section{Concluding discussion}\label{sec:conclu}
We have shown that the new scaled Grassmann Gaussian elimination method
as well as the Riccati method with quasi-optimal patch swapping, compete with 
the continuous orthogonalization method for computing the Evans function.
Both new methods deliver superior accuracy for the same computational
cost when combined with Lie group Magnus integration to advance the 
solution. Moreover, as hoped, numerically these new methods appear to
be robust in the sense that they are insensitive to the choice of 
the matching position in the computational domain. We now outline several
directions in which we plan to use and extend these methods. 

One of the main goals we have had in mind in this paper is that
of large scale spectral problems, in particular the stability of travelling
waves with a multi-dimensional structure. 
There is recent research extending the Evans function
approach in this direction---see Deng and Nii~\cite{DN}, 
Gesztesy, Latushkin and Makarov~\cite{GLM}
and Gesztesy, Latushkin and Zumbrun~\cite{GLZ}. 
From a numerical perspective we have, together with Niesen, 
implemented some of the methods we propose in this paper
in a multi-dimensional context.
In particular, it is well known in autocatalysis and combustion 
that planar travelling fronts can be unstable to transverse perturbations
and develop into steadily propagating travelling fronts with wrinkles.
In Ledoux, Malham, Niesen and Th\"ummler~\cite{LMNT} we show
that the wrinkled fronts themselves develop an instability
as a diffusion parameter is further increased.

For large scale problems the Lie group methods we propose 
using the Magnus expansion may become prohibitive. 
This is because of the effort required to compute the 
matrix exponential---see Moler and Van Loan~\cite{MV}, 
Celledoni and Iserles~\cite{CI},
Munthe-Kaas and Zanna~\cite{MKZ} and Iserles and Zanna~\cite{IZ}.
For the examples we considered this was not an issue. However
it remains to be seen if such Lie group methods will be cost effective
for larger problems---the methods we proposed based on Runge--Kutta
integration such as GGEM-RK can be used as they scale 
favourably with system size.

The constructs and Grassmannian reductions we have considered
in this paper, it turns out, have their origins in the 
control theory literature dating back to the early seventies, 
in particular in the pioneering
papers of Hermann and Martin~\cite{HMc,HM1,HM2,HM3,HM4},
Martin and Hermann~\cite{MH} and Brockett and Byrnes~\cite{BB}.
We also found Bittanti, Laub and Willems~\cite{BLW}, 
Lafortune and Winternitz~\cite{LW},
Rosenthal~\cite{Ro1}, Shayman~\cite{Sh} and Zelikin~\cite{Ze}
particularly useful resources. A future direction we would
like to explore is whether there are any applications of
the numerical methods we have outlined here to practical 
non-autonomous control problems?

Riccati methods in particular also have their origins
in the quantum chemistry literature also dating back to
the early seventies---a recent survey of
these numerical methods can be found in Chou and Wyatt~\cite{CW2}.
However also see Light and Walker~\cite{LightWalker},
Johnson~\cite{Johnson}, Hutson~\cite{Hu} 
and Gray and Manopoulous~\cite{GM}.
In particular the log-derivative and $R$-propagation methods
correspond to special choices of Grassmannian patch in the 
Riccati methods we mention above. Pr\"ufer methods, for which
we can think of the patch evolving, originate even further back;
see Pr\"ufer~\cite{Prufer} and Pryce~\cite{Pryce}.

Of course, our quasi-optimal Gaussian elimination process
for choosing a suitable representative patch was inspired by 
the Schubert cell decomposition of the Grassmann manifold;
see for example Billey~\cite{Bi}, Griffiths and Harris~\cite{GH},
Kleiman and Laksov~\cite{KL}, Kresch~\cite{Kre}, Postnikov~\cite{Po},
Sottile~\cite{So} and, in a somewhat different vein, Kodama~\cite{Ko}. 
Since the Grassmann manifold is the disjoint union of Schubert cells,
the question is, can we express the flow on the Grassmann manifold 
as a flow on Schubert cells 
(see Griffiths and Harris and also Ravi, Rosenthal and Wang~\cite{RRW})?
Can we construct the corresponding flow on the 
cohomological ring of Schubert cycles (Chern~\cite{C}; Fulton~\cite{Fu})?

\section*{Acknowledgements}
We would especially like to thank the anonymous 
referee \#2, who coined the following phrase 
for the Riccati flow in an earlier version of this manuscript:
\emph{this approach is more of Gaussian elimination type 
(including the question of pivoting, which is not discussed here)}.
This comment eventually lead us to the Grassmann
Gaussian elimination method we investigate in this paper.
We also thank this referee for pointing out our 
incomplete analyticity arguments in the second draft.
We would also like to thank referee \#1 for useful
background on the history of the Evans function,
and referee \#3 for making us aware of the 
control theory literature associated with this topic.
Chris Jones, Yuri Latushkin, Bob Pego, 
Bjorn Sandstede and Arnd Scheel 
organised a workshop at AIM in Palo Alto in May 2005
on \emph{Stability Criteria for Multi-Dimensional 
Waves and Patterns}, which instigated 
the topic of this paper. All three authors were 
visiting the Isaac Newton Institute in the Spring 
of 2007 when this research was initiated.
We would like to thank Arieh Iserles and Ernst Hairer
for inviting us and providing so much support and 
enthusiasm. We are also indebted to the facilities
at the Isaac Newton Institute which were invaluable.
We also thank Tom Bridges, Jitse Niesen, Jacques Vanneste 
and Antonella Zanna for stimulating discussions 
on this work. Veerle Ledoux is a postdoctoral
fellow of the Fund of Scientific Research---Flanders
(F.W.O.---Vlaanderen). Vera Th\"ummler was supported by
CRC 701: Spectral Structures and Topological Methods 
in Mathematics.

\label{lastpage}

\end{document}